\documentclass[10pt]{article}
\usepackage{latexsym}
\usepackage{epic}
\usepackage{eepic}
\usepackage{color}
\usepackage{amsmath,multirow,amsfonts}
\usepackage{tikz-cd}
\usetikzlibrary{matrix}
\usetikzlibrary{arrows}
\usepackage{array}
\usepackage[english]{babel}
\usepackage{times}
\usepackage[T1]{fontenc}
\usepackage{amssymb}
\usepackage{graphics}
\usepackage{geometry}
\usepackage{graphicx}
\usepackage{epsfig}
\usepackage{cite}
\usepackage{caption}
 \usepackage{adjustbox}
 \usepackage{tikz}
 \usepackage{lscape}
 \tikzstyle{block}=[draw opacity=.7,line width=1.5cm]
\usetikzlibrary{arrows}
\begin{document}
\newtheorem{def1}{Definition}[section]
\newtheorem{lem}{Lemma}[section]
\newtheorem{exa}{Example}[section]
\newtheorem{thm}{Theorem}[section]
\newtheorem{pro}{Proposition}[section]
\newtheorem{cor}{Corollary}[section]
\newtheorem{rem}{Remark}[section]
\newtheorem{exam}{Example}[section]
\title{Granular $F$-transform and its application}
\author{Abha Tripathi\thanks{abhatripathi745@yahoo.in},~ S.P. Tiwari\thanks{sptiwarimaths@gmail.com}\\
Department of Applied Mathematics\\ Indian Institute of Technology (ISM)\\
Dhanbad-826004, India\\\\ and\\\\\ J. kavikumar\thanks{kavi@uthm.edu.my}\\
Department of Mathematics and Statistics,\\ Universiti Tun Hussein Onn Malaysia,\\
Campus Pagoh 84600, Malaysia }
\date{}
\maketitle
\begin{abstract} This contribution introduces the concept of granular $F$-transform  and investigates its basic properties by using the theory of fuzzy numbers and horizontal membership functions. Further, we present a numerical method based on granular $F$-transform to solve a fuzzy prey-predator model consisting of two prey and one predator due to its natural variability and investigate the existence of the equilibrium points and their stability. 
\end{abstract}
\section{Introduction}
The concept of fuzzy transform ($F$-transform) was firstly introduced by Perfilieva \cite{per}, a theory that attracted the interest of {many researchers}. It has now been greatly expanded upon, and a new chapter in the theory of semi-linear spaces has been opened. The main idea of the $F$-transform is to factorize (or fuzzify) the precise values of independent variables by using a proximity relation, and to average the precise values of dependent variables to an approximation value. The theory of $F$-transform has already been developed and used to real-valued to lattice-valued functions (cf., \cite{{per},{irin1}}), from fuzzy sets to parametrized fuzzy sets \cite{st} and from the single variable to the two (or more variables) (cf., \cite{{mar2},{mar}, {mar1}, {step}}). Recently, several studies have begun to look into $F$-transforms based on any $L$-fuzzy partition of an arbitrary universe (cf., \cite{kh1,jir,mo, mocko, anan,spt1,rus,spt}), where $L$ is a complete residuated lattice. Among these researches, the relationships between $F$-transforms and {semimodule homomorphisms} were introduced in \cite{jir}; a categorical approach of $L$-fuzzy partitions was studied in \cite{mo}; while, the relationships between $F$-transforms and similarity { relations were} discussed in \cite{mocko}. Further, in \cite{anan}, an interesting link among $F$-transforms, $L$-fuzzy topologies/co-topologies and $L$-fuzzy approximation operators (which are concepts used in the study of an operator-oriented view of fuzzy rough set theory) was established, while in \cite{spt1}, the relationship between fuzzy pretopological spaces and spaces with $L$-fuzzy partition was shown. Also, in a different direction, a generalization of $F$-transforms was presented in \cite{rus} by considering the so-called $Q$-module transforms, where $Q$ stands for an unital quantale, while $F$-transforms based on a generalized residuated lattice {were studied} in \cite{spt}. Further, classes of $F$-transforms taking into account the well-known classes of implicators, namely {$R-,S-, QL-$implicators} were discussed in \cite{tri}. The several researches carried out in the application fields of $F$-transforms, e.g., trend-cycle estimation \cite{holc}, {data compression \cite{hut}}, numerical solution of partial differential equations \cite{kh}, scheduling \cite{li}, time series \cite{vil}, data analysis \cite{no}, denoising \cite{Ir}, face recognition \cite{roh}, neural network {approaches \cite{ste} and trading} \cite{to}.\\\\
{It is easy to see how the dynamics of the species would be impacted when they interact. Many studies have been conducted on the problem of the food chain. If the growth rate of one species increases while that of the other decreases during their interaction, we say they are in a predator-prey situation. Such a situation arises when one species (predator) feeds on another species (prey). The Lotka-Volterra model is the first fundamental system representing the interaction between prey and predator species. During the first world war, the Lotka-Volterra model was developed to explain the oscillatory levels of certain fish in the Adriatic sea in \cite{elet,murray}. Several features of predator-prey models have been the subject of many mathematical and ecological studies. There are many factors such as functional response \cite{liu}, competition \cite{cush,gak}, cooperation \cite{elet} affecting dynamics of predator-prey model. The stability and other dynamical behavior of predator-prey models could be found in \cite{aziz,bell,liou,hsu,gakk,brau,oat,past,jp}. The parameters in all the above-cited models are crisp in nature. However, in real-world ecosystems, many parameters may oscillate simultaneously with periodically varying environments. They also change due to natural and human-caused events, such as fire, earthquakes, climate warming, financial crisis, etc. As a result, environmental variables significantly impact the interaction process between the species and its dynamics. Thus the fuzzy mathematical model is more effective than the crisp model. Therefore we have considered the fuzzy set theory to create the prey-predator model. Specifically, the imprecise parameters are replaced by fuzzy numbers in the fuzzy approach. The analysis of the behavior of most phenomena is often based on mathematical models in the form of differential equations. Fuzzy differential equations are equations in which uncertainties are modeled by fuzzy sets (possibility). In recent years, analyzing the dynamical behavior of prey-predator systems whose mathematical models have been considered fuzzy differential equations. Obtaining the solution of the fuzzy differential equation has been investigated under the concept of H-derivative, SGH-derivative \cite{bede}, gH-derivative and g-derivative \cite{bede1}, H$_2$-differentiability \cite{maz}, and gr-derivative \cite{mazan}, it has been researched how to solve fuzzy differential equations. Also, in \cite{hoa,hoa1,long,long1,long2}, fractional calculus has been used to discuss fuzzy differentiable equations. In \cite{sun1}, it has been investigated how stable fuzzy differentiable equations using the second kind of Hukuhara derivative are in the application.\\\\
As long as the underlying fuzzy functions are highly generalized Hukuhara differentiable, the approach described in \cite{sun} could be used for the stability analysis. Thus the limitation of the method as mentioned above relates to the existence of a highly generalized Hukuhara derivative. As the method presented in \cite{sun} is based on what is known as Fuzzy Standard Interval Arithmetic (FSIA), then it has a flaw known as the UBM phenomenon (see \cite{mazan} for more details). A novel fuzzy derivative idea termed granular derivative terms of relative-distance-measure fuzzy interval arithmetic (RDM-FIA) was presented in \cite{mazan} to overcome the drawbacks of the FSIA-based approach. A new idea of the conventional membership functions called the horizontal membership functions (HMFs) proposed in \cite{piegat} was used to construct RDM-FIA. Based on the findings in \cite{land,piegat1,piegat2,son}, it has been established that the RDM-FIA is a more helpful application tool than the FSIA.} \\\\ 
Differential equations cannot always be solved analytically, requiring numerical methods. Therefore in scientific research, numerical methods for solving differential equations have been elaborated frequently. In this connection, differential equations are successfully solved using fuzzy techniques. The fuzzy transform ($F$-transform) introduced by Perfilieva \cite{per} is one of the fuzzy techniques that has been introduced in the literature. An approximation method based on $F$-transform for second order differential equation was introduced in \cite{chen}. In \cite{ali,kh,stev,p}, numerical methods based on $F$-transform to solve initial value problem and boundary value problem were introduced. Also, a numerical method based on $F$-transform to solve a class of delay differential equations by $F$-transform was introduced in \cite{tom}.\\\\
It is to be pointed out here that the numerical solution of a differential equation or fuzzy differential  based on $F$-transform (by using the concept of level sets) was studied, but the numerical solution based on granular $F$-transform of a fuzzy mathematical model, which is represented by fuzzy differential equations under granular differentiability, is yet to be done. In which all parameters and initial conditions can be uncertain. Specifically,
\begin{itemize}
\item[$\bullet$] we introduce the concepts of granular $F$-transform and granular inverse $F$-transform associated with the fuzzy function and discuss some basic results by using the concept of granular metric;
\item[$\bullet$] we formulate a fuzzy prey-predator model and investigate the equilibrium points and their stability in terms of fuzzy numbers;
\item[$\bullet$] we establish a numerical method based on granular $F$-transform for the fuzzy prey-predator model; and
\item[$\bullet$] we present a comparison between two numerical solutions with the exact solution.
\end{itemize}
\section{Preliminaries} Herein, the ideas associated with fuzzy number, horizontal membership function, granular differentiability and fuzzy partition (cf., \cite{maz,mazan,naj,naja,per,kh}, for details). Throughout this chapter, $E^1$ denotes the collection of fuzzy numbers defined on the real number $\mathbb{R}$ and $E^m=E^1\times E^1\times...\times E^1$. The $\alpha$-level sets of $\widehat{a}\in E^1$ is $\widehat{a}^\alpha=[\underline{a}^\alpha,\overline{a}^\alpha]$, where $\underline{a}^\alpha $ and $ \overline{a}^\alpha$ are the left and right end points of $\widehat{a}$.
\begin{def1}\label{grhor}
For a fuzzy number $\widehat{a}: [a_1,a_2] \subseteq \mathbb{R} \rightarrow [0, 1]$ and $\widehat{a}^\alpha=[\underline{a}^\alpha,\overline{a}^\alpha]$, the
{\bf horizontal membership function} is a function $ a^{gr} : [0, 1] \times [0, 1] \rightarrow [a, b]$ such that $a^{gr}(\alpha, \mu_a)=\underline{a}^\alpha+(\overline{a}^\alpha-
\underline{a}^\alpha)\mu_a$, where $\mu_a\in [0, 1]$ is called the relative-distance-measure (RDM) variable.
\end{def1}
\begin{rem}\label{grr1} (i) The horizontal membership function of $\widehat{a}\in E^1$, i.e., $a^{gr}(\alpha, \mu_{a})$ is also denoted by $\mathbb{K}(\widehat{a})$. Also, the $\alpha$-level set of $\widehat{a}$ can be given by \[\mathbb{K}^{-1}(a^{gr}(\alpha, \mu_{a}))=\widehat{a}^\alpha=[\inf\limits_{\beta\geq\alpha}\min\limits_{\mu_{a}}a^{gr}(\beta, \mu_{a}),\sup\limits_{\beta\geq\alpha}\max\limits_{\mu_{a}}a^{gr}(\beta, \mu_{a})].\]
\begin{itemize}
 \item[(ii)] For the fuzzy numbers $\widehat{a_1},\widehat{a_2}\in E^1$, $\widehat{a_1}=\widehat{a_2}$ iff $\mathbb{K}(\widehat{a_1})=
 \mathbb{K}(\widehat{a_2})$ and $\widehat{a_1}\geq \widehat{a_2}$ if $\mathbb{K}(\widehat{a_1})\geq
 \mathbb{K}(\widehat{a_2}),\,
 \forall\,\mu_{a_1}=\mu_{a_2}\in [0,1]$. 
 \item[(iii)] Let each of addition, subtraction, multiplication and division operations between fuzzy numbers $\widehat{a_1},\widehat{a_2}\in E^1$ be represented by $\odot$. Therefore $\widehat{a_1}\odot\widehat{a_2}=\widehat{m}\in E^1$ iff $\mathbb{K}(\widehat{m})=
 \mathbb{K}(\widehat{a_1})\odot\mathbb{K}(\widehat{a_2})$.
 \item[(iv)] Let $\widehat{a_1},\widehat{a_2},\widehat{a_3}\in E^1$. Then we have
 \item[$\bullet$] $\widehat{a_1}-\widehat{a_2}=-(\widehat{a_1}+\widehat{a_2})$,
 \item[$\bullet$] $\widehat{a_1}-\widehat{a_1}=0$,
 \item[$\bullet$] $\widehat{a_1}\div\widehat{a_1}=1$, and
 \item[$\bullet$] $(\widehat{a_1}+\widehat{a_2})\widehat{a_3}=\widehat{a_1}\widehat{a_3}+\widehat{a_2}\widehat{a_3}$.
 \item[(v)] A fuzzy function is a generalization of a classical function in which the domain or range, or both, is a subset of the fuzzy numbers set.
 \end{itemize}
 \end{rem}
\begin{def1}
 The horizontal membership function of $\widehat{g}(\widehat{p_{1}}(u),\widehat{p_{2}}(u),...,\widehat{p}_m(u))$ is defined by $\mathbb{K}(\widehat{g}(\mathbb{K}(\widehat{p}_1(u)),\mathbb{K}(\widehat{p_{2}}(u)),...,\mathbb{K}(\widehat{p}_m(u))))$, where $\widehat{g}:E^m\rightarrow E^1$ and $\widehat{p_i}:[a_1,a_2]\subseteq\mathbb{R}\rightarrow E^1,\,i=1,2,...,m$, are fuzzy functions.
\end{def1}
\begin{def1}\label{grmetric}
A fuzzy function $d^{gr}:E^1\times E^1\rightarrow \mathbb{R}^+\cup \{0\}$ 
is called the {\bf granular metric} if \[d^{gr}(\widehat{a_1},\widehat{a_2})=\sup\limits_{\alpha\in[0,1]}\max\limits_{\mu_{a_1},\mu_{a_2}\in[0,1]}|a_1^{gr}(\alpha,\mu_{a_1})-a_2^{gr}(\alpha,\mu_{a_2})|,\,\forall\,\widehat{a_1},\widehat{a_2}\in E^1.\]
\end{def1}
\begin{def1}\label{grcont}
A fuzzy function $\widehat{g}: [a_1,a_2] \subseteq \mathbb{R} \rightarrow E^1$ is called the {\bf granular continuous} (gr-continuous) if for all $u_0\in[a_1,a_2],\epsilon>0$ there is a $\delta>0$ such that $$d^{gr}(\widehat{g}(u),\widehat{g}(u_0))<\epsilon, \,\,\mbox{whenever}\, |u-u_0|<\delta.$$
\end{def1}
\begin{def1} \label{grgrdif}
A fuzzy function $\widehat{g}: [a_1,a_2] \subseteq \mathbb{R} \rightarrow E^1$ is called the {\bf granular differentiable} (gr-differentiable) at the point
$u\in [a_1,a_2]$ if there is a fuzzy number $\mathcal{D}^{gr} \widehat{g} (u)\in E^1$ such that the following limit exists:
\[\mbox{lim}_{h\to 0}\frac{\widehat{g}(u + h)-\widehat{g}(u)}{h}=\mathcal{D}^{gr} \widehat{g} (u).\]
\end{def1}
\begin{pro}\label{grgrdif1}
At any point $u\in [a_1,a_2]$, a fuzzy function $\widehat{g} : [a_1,a_2] \subseteq \mathbb{R} \rightarrow E^1$ is gr-differentiable iff $\mathbb{K}(\widehat{g})$ is differentiable w.r.t. $u$ at that point. In addition, 
$$\mathbb{K}(\mathcal{D}^{gr} \widehat{g} (u))=\frac{\partial}{\partial t}\mathbb{K}(\widehat{g} (u)).$$
\end{pro}
\begin{pro}\label{grgrpdif}
The fuzzy function $\widehat{g}(\widehat{p}_1(u),\widehat{p_{2}}(u),...,\widehat{p}_m(u))$ is gr-partial differentiable w.r.t. $\widehat{p}_i(u)$ iff its horizontal membership function is differentiable w.r.t. the horizontal membership function of $\widehat{p}_i(u)$, where $\widehat{g} : E^m \rightarrow E^1$ and $\widehat{p}_i : [a_1,a_2] \subseteq \mathbb{R} \rightarrow E^1,\,i = 1, 2, . . . ,m$, are fuzzy functions.. Moreover,
\[\mathbb{K}(\frac{\partial^{gr}}{\partial \widehat{p}_i} \widehat{g}(\widehat{p}_1(u),\widehat{p_{2}}(u),...,\widehat{p}_m(u)))=\frac{\partial}{
\partial \mathbb{K}(\widehat{p}_i)}\mathbb{K}(\widehat{g}(\mathbb{K}(\widehat{p}_1(u)),\mathbb{K}(\widehat{p_{2}}(u)),...,\mathbb{K}(\widehat{p}_m(u))))
.\]
\end{pro}
\begin{def1}\label{grgrint}
Let $g^{gr}(u, \alpha, \mu_g )$ is integrable on $u\in [a_1,a_2]$, where $g^{gr}(u, \alpha, \mu_g )$ is a horizontal membership
function of a gr-continuous fuzzy function $\widehat{g} : [a_1,a_2]\subseteq \mathbb{R}\rightarrow E^1$. In addition, let integral of $\widehat{g} $ on $[a_1,a_2]$ be represented by $\oint^{a_2}_{a_1} \widehat{g} (u)du$. If there exists a fuzzy number $\widehat{m} = \oint^{a_2}_{a_1} \widehat{g} (u)du$ such that $\mathbb{K}(\widehat{m})=\int^{a_2}_{a_1} \mathbb{K}(\widehat{g} (u))du$, then the fuzzy function $\widehat{g}$ is called the {\bf granular fuzzy integrable} on $[a_1,a_2]$.
\end{def1}
\begin{def1}\label{grgrpol}
A {\bf granular fuzzy polynomials} is an expression consisting of fuzzy variables and fuzzy coefficients that involves only the granular operations of addition, subtraction, multiplication and nom-negative with integer exponents of fuzzy variables. 
\end{def1}
For example- $\widehat{g}(\widehat{p})=\widehat{a}_m\widehat{p}^n+\widehat{a}_{m-1}\widehat{p}^{m-1}+...+\widehat{a}_1\widehat{p}+\widehat{a}_0,\,\forall\,\widehat{a}_i\in E^1,\,i=1,...,m$.
\begin{def1}\label{grroot}
A {\bf fuzzy root} of a granular fuzzy polynomial $\widehat{g}(\widehat{p})$ is a fuzzy number $\widehat{p}_i$ such that $\widehat{g}(\widehat{p}_i)=0$.
\end{def1}
\begin{rem}\label{grr2}
It is easy to check that if $\widehat{p}_k$ is a fuzzy root of $\widehat{g}(\widehat{p})$, then $\mathbb{K}(\widehat{p}_i)$ is a root of $\mathbb{K}(\widehat{g}(\mathbb{K}(\widehat{p})))$, i.e., $\widehat{g}(\widehat{p}_i)=0\Rightarrow \mathbb{K}(\widehat{g}(\mathbb{K}(\widehat{p_i}) ))=0$.
\end{rem}
Next, the concepts of fuzzy partition and $F$-transform introduced by \cite{per} are recalled. 
\begin{def1}\label{grFP} Let $u_1 <...<u_m$ be fixed nodes within $[a_1,a_2]$, where $u_1 = a_1, u_m = a_2$ and $m \geq 2$. Then fuzzy sets $P_1, . . . , P_m$ (are called {\bf basic functions}) identified with their membership functions $P_1(u), . . . , P_m(u)$ defined on $[a_1,a_2]$, form a {\bf fuzzy partition} of $[a_1,a_2]$ if they satisfy the following properties for $i = 1, . . . ,m$,
\begin{itemize}
  \item[(i)] $P_i :[a_1,a_2]\rightarrow [0,1],\, P_i(u_i)=1$,
   \item[(ii)] $P_i(u) = 0$ if $t\notin (u_{i-1},u_{i+1})$, where for the uniformity of the notation, we put $x_0 = a_1$ and $u_{m+1} = a_2$,
   \item[(iii)] $P_i(u)$ is continuous,
   \item[(iv)] $P_i(u), i= 2, . . . ,m$, strictly increases on $[u_{i-1}, u_{i}]$ and $P_i(u), i = 1, . . . ,m-1$, strictly decreases on $[u_i, u_{i+1}]$, and
 \item[(v)] $\sum_{i=1}^{m}P_i(u)=1,\,\forall\,u\in[a_1,a_2],\,$.
   \end{itemize}
\end{def1}
If $u_1, . . . , u_m,\,m \geq 2$,
are equidistant, then the fuzzy partition of $[a_1,a_2]$ is $h$-uniform. That is to say that $u_i= a_1+ h (i -1) , i = 1, . . . ,m$, where $h = \frac{a_2-a_1}{m-1 },m\geq2$ and the two
additional properties are satisfied:
\begin{itemize}
  \item[(i)] $P_i(u_i -u) = P_i(u_i +u ), i= 2, . . . ,m -1, \forall\, u\in [0, h]$, and
  \item[(ii)] $P_i(u) = P_{i-1}(u -h)$ and $P_{i+1}(u) = P_i(u -h), i= 2, . . . ,m -1,\,\forall u\in [u_k, u_{i+1}]$.
\end{itemize} 
\begin{lem}\label{grlem}
Let $P_1,...,P_m,\,m\geq 3$, be basic functions that form the uniform fuzzy partition of $[a_1,a_2]$. Then
\[\int_{u_1}^{u_2}P_1(u)du=\int_{u_{m-1}}^{u_{m}}P_m(u)du=\frac{h}{2},\,\,and\]
\[\int_{u_{i-1}}^{u_{i+1}}P_i(u)du=h,\, i=2,...,m-1\]
where $h$ is the distance between two adjacent nodes.
\end{lem}
\begin{def1}\label{grft} Let $P_1, ...,P_m$ be basic functions that form a fuzzy partition of $[a_1,a_2]$ and ${g}$ be the continuous function on $[a_1,a_2]$. Then the {\bf $F$-transform} of ${g}$ w.r.t. $P_1,...,P_m$ is the $m$-tuple
of real numbers (components) $F[{g}] = (F_1[{g}], ..., F_m[{g}])$, where
\[F_i[{g}]=\dfrac{\int^{a_2}_{a_1} {g}(u)P_i(u)du}{\int^{a_2}_{a_1} P_i(u)du},\,i=1,2,...,m,\]
where $F_i[{g}]$ (or simply $F_i)$) is $i$-th component of the $F$-transform.
\end{def1}
\begin{def1}\label{grift} Let $F[\widehat{g}] = (F_1, ..., F_m)$ be the $F$-transform of ${g}$ w.r.t. $P_1, ...,P_m$. Then the function 
$$\hat{g}_m(u)=\sum_{i=1}^{m}F_i[\widehat{g}]P_i(u),\,\forall\,u\in[a_1,a_2],$$
is called the {\bf inverse $F$-transform}.
\end{def1}
\section{Granular fuzzy transform} 
This section introduces the concept of the granular fuzzy transform ($F$-transform), which creates a relationship between a set of gr-continuous fuzzy functions from $[a_1,a_2]$ to $E^1$ and the collection of fuzzy numbers. The formula, which will refer to as a granular inverse $F$-transform (inversion formula), converts a fuzzy number into another gr-continuous fuzzy function that approximates the original one. Further, we investigate their basic properties. Now, we initiate with the following.
\begin{def1}\label{grgrft} Let $P_1, ...,P_m$ be basic functions that form a fuzzy partition of $[a_1,a_2]$ and $\widehat{g}: [a_1,a_2]\rightarrow E^1$ be the gr-continuous fuzzy function. Further, let $\widehat{g}(u)P_i(u)$ be the granular fuzzy integrable on $[a_1,a_2]$. Then the {\bf granular $F$-transform} of $\widehat{g}$ w.r.t. $P_1,...,P_m$ is the $m$-tuple of fuzzy numbers (components) $F^{gr}[\widehat{g}] = (F^{gr}_1[\widehat{g}], ..., F^{gr}_m[\widehat{g}])$, where
\[F^{gr}_i[\widehat{g}]=\dfrac{\oint^{a_2}_{a_1}\widehat{g}(u)P_i(u)du}{\oint^{a_2}_{a_1} P_i(u)du},\,i=1,2,...,m,\]
where $F^{gr}_i[\widehat{g}]$ (or simply $F^{gr}_i)$) is $i$-th component of the granular $F$-transform.
\end{def1}
Further, let $P_1, . . . ,P_m$ be basic functions that form a uniform fuzzy partition of $[a_1,a_2]$. Then from Lemma \ref{grlem}, the components of granular $F$-transform can be written as 
\begin{eqnarray*}
\nonumber F^{gr}_1[\widehat{g}]&=&\dfrac{2}{h}\oint^{u_2}_{u_1}\widehat{g}(u)P_1(u)du,\\
\nonumber F^{gr}_m[\widehat{g}]&=&\dfrac{2}{h}\oint^{u_m}_{u_{m-1}}\widehat{g}(u)P_m(u)du,\\
\nonumber F^{gr}_i[\widehat{g}]&=&\dfrac{1}{h}\oint^{u_{i+1}}_{u_{i-1}}\widehat{g}(u)P_i(u)du,\,i=2,...,m-1.\\
\end{eqnarray*}
\begin{rem}\label{grr3}
(i) From Definition \ref{grgrint}, it is clear that $\mathbb{K}(F^{gr}[\widehat{g}])=F[\mathbb{K}(\widehat{g})]$, where $F[.]$ is $F$-transform in \cite{per}; and\\\\
(ii) For any crisp function $g:[a_1,a_2]\rightarrow \mathbb{R},\,F^{gr}[g]=F[g]$, i.e., granular $F$-transform becomes $F$-transform in \cite{per}.
\end{rem}
The following are towards the some fundamental properties of the granular $F$-transform associated with the fuzzy function. 
\begin{pro}\label{grp1}
Let $\widehat{g_1},\widehat{g_2}:[a_1,a_2]\rightarrow E^1$ be gr-continuous fuzzy functions. Then the granular $F$-transform is linear, i.e.,
\[F^{gr}[\widehat{a_1}\widehat{g_1}+\widehat{a_2}\widehat{g_2}]=\widehat{a_1}F^{gr}[\widehat{g_1}]+\widehat{a_2}F^{gr}[\widehat{g_2}],\,\forall\,\widehat{a_1},\widehat{a_2}\in E^1.\]
\end{pro}
\textbf{Proof:} Let $\widehat{a_1},\widehat{a_2}\in E^1$. Then from Remarks \ref{grr1} and \ref{grr3}, we have
\begin{eqnarray*}
\mathbb{K}(F^{gr}[\widehat{a_1}\widehat{g_1}+\widehat{a_2}\widehat{g_2}])&=&F[\mathbb{K}(\widehat{a_1}\widehat{g_1}+\widehat{a_2}\widehat{g_2})]\\
&=&F[\mathbb{K}(\widehat{a_1})\mathbb{K}(\widehat{g_1})+\mathbb{K}(\widehat{a_2})\mathbb{K}(\widehat{g_2})]\\
&=&\mathbb{K}(\widehat{a_1})F[\mathbb{K}(\widehat{g_1})]+\mathbb{K}(\widehat{a_2})F[\mathbb{K}(\widehat{g_2})]\\
&=&\mathbb{K}(\widehat{a_1})\mathbb{K}(F^{gr}[\widehat{g_1}])+\mathbb{K}(\widehat{a_2})\mathbb{K}(F^{gr}[\widehat{g_2}])\\
&=&\mathbb{K}(\widehat{a_1}F^{gr}[\widehat{g_1}]+\widehat{a_2}F^{gr}[\widehat{g_2}]).
\end{eqnarray*}
Thus $F^{gr}[\widehat{a_1}\widehat{g_1}+\widehat{a_2}\widehat{g_2}]=\widehat{a_1}F^{gr}[\widehat{g_1}]+\widehat{a_2}F^{gr}[\widehat{g_2}]$.
\begin{pro}\label{grp2} Let $P_1, . . . ,P_m,m\geq3$, be basic functions that form a uniform fuzzy partition of $[a_1,a_2]$ and $\widehat{g}:[a_1,a_2]\rightarrow E^1$ be the gr-differentiable fuzzy function. Moreover, $\mathcal{D}^{gr}\widehat{g}$ is the gr-continuous fuzzy function on $[a_1,a_2]$. Then $F^{gr}[\widehat{g}]$ gives minimum to the fuzzy function defined on $S\subseteq E^1$ such that
\[\widehat{\phi}(\widehat{a})={\oint^{a_2}_{a_1}(\widehat{g}(u)-\widehat{a})^2P_i(u)du},\,i=1,2,...,m,\]
where $S=\{\widehat{a}\in E^1:\widehat{g}(a_1)\leq \widehat{a}\leq\widehat{g}(a_2)\}.$
\end{pro}
\textbf{Proof:} Since the fuzzy function $(\widehat{g}(u)-\widehat{a})^2P_i(u)$ is the gr-continuously differentiable w.r.t. $\widehat{a}$ on $S$. Therefore we can write 
\begin{eqnarray}\label{grmin}
\nonumber\mathcal{D}^{gr}\widehat{\phi}(\widehat{a})&=&\mathcal{D}^{gr}\left(\oint^{a_2}_{a_1}(\widehat{g}(u)-\widehat{a})^2 P_i(u)du\right).
\end{eqnarray}
Now, from Remark \ref{grr1} and Proposition \ref{grgrdif1}, we have
\begin{eqnarray}\label{grmin1}
\nonumber\mathbb{K}(\mathcal{D}^{gr}\widehat{\phi}(\widehat{a}))&=&\mathbb{K}\left(\mathcal{D}^{gr}\left(\oint^{a_2}_{a}(\widehat{g}(u)-\widehat{a})^2P_i(u)du\right)\right)\\
\frac{\partial }{\partial \nonumber \mathbb{K}(\widehat{a})}\mathbb{K}(\widehat{\phi}(\mathbb{K}(\widehat{a})))&=&-2\int^{a_2}_{a_1}(\mathbb{K}(\widehat{g}(u))-\mathbb{K}(\widehat{a}))P_i(u)du.
\end{eqnarray}
Also, it can be easily check that $\mathbb{K}(\widehat{\phi}(\mathbb{K}(\widehat{a})))$ reaches its minimum at that time, resulting in a solution to the equation $\frac{\partial }{\partial \mathbb{K}(\widehat{q})}\mathbb{K}(\widehat{\phi}(\mathbb{K}(\widehat{a})))=0$, i.e.
\[ \mathbb{K}(\widehat{a})=\dfrac{\oint^{a_2}_{a_1}\mathbb{K}(\widehat{g}(u))P_i(u)du}{\oint^{a_2}_{a_1} P_i(u)du}=\mathbb{K}(F^{gr}_i[\widehat{g}]).\]
Thus we have $\widehat{a}=F^{gr}_i[\widehat{g}]$.\\\\
In the following, we demonstrate that several assumptions about the smoothness of $\widehat{g}$ can be used to estimate each granular $F$-transform component $F^{gr}_i,\,i = 1,...,m$.
\begin{pro}\label{grp3}
Let $P_1, . . . ,P_m,\,m\geq3$, be basic functions that form a uniform fuzzy partition of $[a_1,a_2]$ and $\widehat{g}:[a_1,a_2]\rightarrow E^1$ be the gr-continuous fuzzy function. In addition, let $F^{gr}_i[\widehat{g}], i= 1, . . . ,m$, be the granular $F$-transform components of $\widehat{g}$ w.r.t. $P_1, . . . ,P_m$. Then
\begin{eqnarray*}
 d^{gr}(\widehat{g}(u),F^{gr}_i[\widehat{g}])&\leq& \omega(2h,\widehat{g}),\\
d^{gr}(\widehat{g}(u),F^{gr}_{i+1}[\widehat{g}])&\leq& \omega(2h,\widehat{g}),\,\forall\,i= 1, . . . ,m-1,\, u\in [u_i, u_{i+1}],
\end{eqnarray*}
where $h=\frac{a_2-a_1}{m-1}$ and $\omega(2h,\widehat{g})=\max\limits_{|\delta|\leq 2h}\max\limits_{t'\in[a_1,a_2-\delta]}d^{gr}(\widehat{g}(u'+\delta),\widehat{g}(u'))$ is the modulus of gr-continuity of $\widehat{g}$ on $[a_1,a_2]$.
\end{pro}
\textbf{Proof:} Let $u\in[u_i,u_{i+1}]$ and $1\leq i\leq m-1$. Then 
\begin{eqnarray*}
 d^{gr}(\widehat{g}(u),F^{gr}_i[\widehat{g}])&=&\sup\limits_{\alpha\in [0,1]}\max\limits_{\mu_g\in [0,1]}|g^{gr}(u,\alpha,\mu_g)-F_i[\mathbb{K}(\widehat{g})]|\\
 &=&\sup\limits_{\alpha\in [0,1]}\max\limits_{\mu_g\in [0,1]}|g^{gr}(u,\alpha,\mu_g)-\frac{1}{h}{\int^{u_{i+1}}_{u_{i-1}}g^{gr}(u',\alpha,\mu_g)P_i(u')du'}|\\
 &=&\sup\limits_{\alpha\in [0,1]}\max\limits_{\mu_g\in [0,1]}|\frac{1}{h}{\int^{u_{i+1}}_{u_{i-1}}(g^{gr}(u,\alpha,\mu_g)-g^{gr}(u',\alpha,\mu_g))P_i(u')du'}|\\
 &\leq&\sup\limits_{\alpha\in [0,1]}\max\limits_{\mu_g\in [0,1]}\frac{1}{h}{\int^{u_{i+1}}_{u_{i-1}}|g^{gr}(u,\alpha,\mu_g)-g^{gr}(u',\alpha,\mu_g)|P_i(u')du'}\\
 &\leq&\sup\limits_{\alpha\in [0,1]}\max\limits_{\mu_g\in [0,1]}\max\limits_{|\delta|\leq 2h}\max\limits_{t'\in[a_1,a_2-\delta]}|g^{gr}(u'+\delta,\alpha,\mu_g)-g^{gr}(u',\alpha,\mu_g)|\\
 &=&\max\limits_{|\delta|\leq 2h}\max\limits_{t'\in[a_1,a_2-\delta]}d^{gr}(\widehat{g}(u'+\delta),\widehat{g}(u'))\\
 &=&\omega(2h,\widehat{g}).
\end{eqnarray*}
Similarly, we can show that $d^{gr}(\widehat{g}(u),F^{gr}_{i+1}[\widehat{g}])\leq \omega(2h,\widehat{g})$.
\begin{pro}\label{grp4}
Let $P_1, . . . ,P_m,\,m\geq3$, be basic functions which form a uniform fuzzy partition of $[a_1,a_2]$ and $\widehat{g}:[a_1,a_2]\rightarrow E^1$ be the gr-continuous fuzzy function. Then for $c_{i1} \in [u_{i-1}, u_{i}]$ and $c_{i2} \in [u_{i}, u_{i+1}]$, the components of granular $F$-transform satisfy the following condition:
\[F^{gr}_i[\widehat{g}]=\frac{1}{h}{\oint^{c_{i2}}_{c_{i1}}\widehat{g}(u)du},\,i= 2, . . . ,m-1.\]
Moreover, if $i=1$ and $i=m$, then there exists $c\in [u_{1},u_{2}]$ \mbox{and} $c\in [u_{m-1}, u_{m}]$ such that
\[F^{gr}_1[\widehat{g}]=\frac{2}{h}{\oint^{c}_{u_1}}\widehat{g}(u)du,\]
\[F^{gr}_m[\widehat{g}]=\frac{2}{h}{\oint^{u_{m}}_{c}}\widehat{g}(u)du,\,\,\mbox{respectively}.\]
\end{pro}
\textbf{Proof:} Let $2\leq i \leq m-1$. Then from Definition \ref{grFP} , $P_i(u)$ increases monotonically on $[u_{i-1}, u_{i}]$ and decreases monotonically on $[u_{i}, u_{i+1}]$,
we find
\[F^{gr}_i[\widehat{g}]=\frac{1}{h}{\oint^{u_{i+1}}_{u_{i-1}}\widehat{g}(u)P_i(u)du}.\]
Now, from Remark \ref{grr1}, we can write
\begin{eqnarray}
\nonumber\mathbb{K}(F^{gr}_i[\widehat{g}])&=&\mathbb{K}\left(\frac{1}{h}{\oint^{u_{i+1}}_{u_{i-1}}\widehat{g}(u)P_i(u)du}\right)\\
\nonumber&=& \frac{1}{h}\int^{u_{i+1}}_{u_{i-1}}\mathbb{K}(\widehat{g}(u))P_i(u)du\\
\nonumber&=& \frac{1}{h}\int^{u_{i}}_{u_{i-1}}\mathbb{K}(\widehat{g}(u))P_i(u)du+\frac{1}{h}\int^{u_{i+1}}_{u_{i}}\mathbb{K}(\widehat{g}(u))P_i(u)du\\
\nonumber&=&\frac{1}{h}{\int^{u_{i}}_{c_{i1}}}\mathbb{K}(\widehat{g}(u))du+\frac{1}{h}\int^{c_{i2}}_{u_{i}}\mathbb{K}(\widehat{g}(u))du\\
\nonumber&=&\frac{1}{h}{\int^{c_{i2}}_{c_{i1}}}\mathbb{K}(\widehat{g}(u))du\\
\nonumber&=& \mathbb{K}\left(\frac{1}{h}{\oint^{c_{i2}}_{c_{i1}}}\widehat{g}(u)du\right).
\end{eqnarray}
{Thus $F^{gr}_i[\widehat{g}]=\frac{1}{h}{\oint^{c_{i2}}_{c_{i1}}}\widehat{g}(u)du$. Similarly, we can obtain} $${F^{gr}_1[\widehat{g}]=\frac{2}{h}{\oint^{c}_{u_1}}\widehat{g}(u)du\left(F^{gr}_m[\widehat{g}]=\frac{2}{h}{\oint^{u_{m}}_{c}}\widehat{g}(u)du\right).}$$
From the above result, $F^{gr}_i$ can be defined as an integral mean value of $\widehat{g}$ within the interval $[c_{i1}, c_{i2}]$ that accumulates information about the fuzzy function $\widehat{g}$. However, for the given fuzzy function and nodes of the partition, this interval cannot be specified precisely. We can demonstrate the closeness between the components of the granular $F$-transform and fuzzy function at the corresponding nodes. The estimation of closeness is given below.
\begin{pro}\label{grp5}
Let the conditions of Proposition \ref{grp4} be hold and the fuzzy function $\widehat{g}$ be twice gr-continuously differentiable in $(a_1,a_2)$. Then for all $i= 1, . . . ,m$
\[F^{gr}_i[\widehat{g}] = \widehat{g}(u_k) + O(h^2).\]
\end{pro}
\textbf{Proof:} Let $2\leq i\leq m-1$. Then
\[\vspace{-1mm}{F^{gr}_i[\widehat{g}]=\frac{1}{h}{\oint^{u_{i+1}}_{u_{i-1}}\widehat{g}(u)P_i(u)du}.}\]
From Remark \ref{grr1}, the above expression can be written as
\begin{eqnarray}
\nonumber\mathbb{K}(F^{gr}_i[\widehat{g}])&=&\mathbb{K}\left(\frac{1}{h}{\oint^{u_{i+1}}_{u_{i-1}}\widehat{g}(u)P_i(u)du}\right)\\
\nonumber&=& \frac{1}{h}\int^{u_{i+1}}_{u_{i-1}}\mathbb{K}(\widehat{g}(u))P_i(u)du.
\end{eqnarray}
Now, by using trapezoidal rule with nodes $u_{i-1},u_i,u_{i+1}$, we have 
\begin{eqnarray}
\nonumber \mathbb{K}(F^{gr}_i[\widehat{g}]) &=&\frac{1}{h}\int^{u_{i+1}}_{u_{i-1}}\mathbb{K}(\widehat{g}(u))P_i(u)du\\
\nonumber&=&\frac{1}{h}\frac{h}{2}\left[\mathbb{K}(\widehat{g}(u_{i-1}))P_i(u_{i-1})+\mathbb{K}(\widehat{g}(u_{i+1}))P_i(u_{i+1})+ 2\mathbb{K}(\widehat{g}(u_{i}))P_i(u_{k})\right]+O(h^2)\\
\nonumber&=&\mathbb{K}(\widehat{g}(u_{i}))+ O(h^2)\\
\nonumber&=&\mathbb{K}(\widehat{g}(u_{i})+ O(h^2)).
\end{eqnarray}
Thus $F^{gr}_i[\widehat{g}] = \widehat{g} (u_i) + O(h^2)$. Similarly, for $i=1\,(i=m)$, $F^{gr}_1[\widehat{g}] = \widehat{g} (u_1) + O(h^2)\, (F^{gr}_m[\widehat{g}] = \widehat{g} (u_m) + O(h^2))$.\\\\
The following is towards the notion of the granular inverse $F$-transform associated with the fuzzy function. Now, we initiate with the following.
\begin{def1}\label{grgrift} Let $F^{gr}[\widehat{g}] = (F^{gr}_1, ..., F^{gr}_m)$ be the granular $F$-transform of $\widehat{g}$ w.r.t. $P_1, ...,P_m$. Then the function 
$$\widehat{g}^{gr}_m(u)=\sum_{i=1}^{m}F^{gr}_i[\widehat{g}]P_i(u),\,\forall\,u\in[a_1,a_2],$$
is called the {\bf granular inverse $F$-transform}.
\end{def1}
\begin{rem}
(i) From Definition \ref{grgrift}, it is easy to say that $\mathbb{K}(\widehat{g}^{gr}_m(u))=\sum_{i=1}^{m}\mathbb{K}(F^{gr}_i[\widehat{g}])P_i(u)=\sum_{i=1}^{m}F_i[\mathbb{K}(\widehat{g})]P_i(u)=\hat{g}_m(u)$, where $F[.],\hat{g}_m(u)$ are $F$-transform and inverse $F$-transform in \cite{per}, respectively; and\\\\
(ii) For any crisp function $g:[a_1,a_2]\rightarrow \mathbb{R},\widehat{g}^{gr}_m(u)=\hat{g}_m(u)$, i.e., granular $F$-transform and granular inverse $F$-transform become $F$-transform and inverse $F$-transform in \cite{per}, respectively.
\end{rem}
Below, we demonstrate that the granular inverse $F$-transform $\widehat{g}^{gr}_m$ can approximate the original gr-continuous fuzzy function $\widehat{g}$ with an arbitrary precision. To do this, we need the following remark.
\begin{rem}\label{grremark} (i) Let $\widehat{g}:[a_1,a_2]\rightarrow E^1$ be the gr-continuous fuzzy function. Then for all $\epsilon > 0$ there is a $\delta>0$ such that for all $u\in[a_1,a_2]$, $d^{gr}(\widehat{g}(u),\widehat{g}(u))<\epsilon$, whenever $|t-t'|<\delta$. As $d^{gr}(\widehat{g}(u),\widehat{g}(u'))=\sup\limits_{\alpha\in [0,1]}\max\limits_{\mu_g\in [0,1]}|g^{gr}(u,\alpha,\mu_g)-g^{gr}(u',\alpha,\mu_g)|$. So, by using the fact that a function continuous on $[a_1,a_2]$, is uniformly continuous, we have $$\sup\limits_{\alpha\in [0,1]}\max\limits_{\mu_g\in [0,1]}|g^{gr}(u,\alpha,\mu_g)-g^{gr}(u',\alpha,\mu_g)|<\epsilon.$$
(ii) \vspace{-1mm}{By using (i), we can show that for all $u\in [a_1,a_2]$ and $i=1,...,m-1,$ $$d^{gr}(\widehat{g}(u),F_i^{gr}[\widehat{g}])\leq\epsilon,\,d^{gr}(\widehat{g}(u),F_{i+1}^{gr}[\widehat{g}])\leq\epsilon.$$
Now, from Proposition \ref{grp3}, we have}
\begin{eqnarray}
\nonumber d^{gr}(\widehat{g}(u),F_i^{gr}[\widehat{g}])&\leq&\sup\limits_{\alpha\in [0,1]}\max\limits_{\mu_g\in [0,1]}\frac{1}{h}{\int^{u_{i+1}}_{u_{i-1}}|g^{gr}(u,\alpha,\mu_g)-g^{gr}(u',\alpha,\mu_g)|P_i(u')du'}\\
\nonumber&<&\frac{\epsilon}{h}{\int^{u_{i+1}}_{u_{i-1}}P_i(u')du'}\\
\nonumber&=&\epsilon.
\end{eqnarray}
{Similarly, we can show that $d^{gr}(\widehat{g}(u),F_{i+1}^{gr}[\widehat{g}])\leq\epsilon$.}
\end{rem}
\begin{pro}\label{grp6}
Let $\widehat{g}:[a_1,a_2]\rightarrow E^1$ be the gr-continuous fuzzy function. Then for all $\epsilon > 0$ there is $m_\epsilon$ and the fuzzy partition
$P_1, . . . ,P_{m_\epsilon}$ of $[a_1,a_2]$ such that for all $u\in[a_1,a_2]$,
\[d^{gr}(\widehat{{g}}_m(u) , \widehat{g}^{gr}_{m_\epsilon}(u))\leq\epsilon,\]
where $\widehat{g}^{gr}_{m_\epsilon}$ is the granular inverse $F$-transform of $\widehat{g}$.
\end{pro}
\textbf{Proof:} Let $u\in[a_1,a_2]$ and $\widehat{g}:[a_1,a_2]\rightarrow E^1$ be the gr-continuous fuzzy function. Then from Definition \ref{grmetric} and Remark \ref{grremark} (ii),
\begin{eqnarray*}
\nonumber d^{gr}(\widehat{{g}}(u), \widehat{g}^{gr}_m(u))&=&\sup\limits_{\alpha\in [0,1]}\max\limits_{\mu_g\in [0,1]}|g^{gr}(u,\alpha,\mu_g)-{g}^{gr}_m(u,\alpha,\mu_g)|\\
\nonumber&=&\sup\limits_{\alpha\in [0,1]}\max\limits_{\mu_g\in [0,1]}|g^{gr}(u,\alpha,\mu_g)-\sum_{i=1}^{m}F_i[\mathbb{K}(\widehat{g})]P_i(u)|\\
\nonumber&=&\sup\limits_{\alpha\in [0,1]}\max\limits_{\mu_g\in [0,1]}|g^{gr}(u,\alpha,\mu_g)\sum_{i=1}^{m}P_i(u)-\sum_{i=1}^{m}F_i[\mathbb{K}(\widehat{g})]P_i(u)|\\
\nonumber&=&\sup\limits_{\alpha\in [0,1]}\max\limits_{\mu_g\in [0,1]}|\sum_{i=1}^{m}(g^{gr}(u,\alpha,\mu_g)-F_i[\mathbb{K}(\widehat{g})])P_i(u)|\\
\nonumber&\leq&\sup\limits_{\alpha\in [0,1]}\max\limits_{\mu_g\in [0,1]}\sum_{i=1}^{m}|g^{gr}(u,\alpha,\mu_g)-F_i[\mathbb{K}(\widehat{g})]|P_i(u)\\
\nonumber&\leq&\epsilon\sum_{i=1}^{m}P_i(u)\\
\nonumber&=&\epsilon.
\end{eqnarray*}
Proposition \ref{grp6} can be formulated for the uniform fuzzy partitions of $[a_1,a_2]$.
\begin{cor}\label{grco1}
Let $\widehat{g}:[a_1,a_2]\rightarrow E^1$ be the gr-continuous fuzzy function and $\{(P^{(m)}_1, . . . ,P^{(m)}_{m})_m\}$ be a sequence of uniform fuzzy partitions of $[a_1,a_2]$ for all $m$ . In addition, let $\{\widehat{g}^{gr}_m(u)\}$ be the sequence of granular inverse $F$-transforms w.r.t. $\{(P^{(m)}_1, . . . ,P^{(m)}_{m})_m\}$, respectively. Then for all $\epsilon > 0$ there exist $m_\epsilon$ such that $m>m_\epsilon$ and
\[d^{gr}(\widehat{{g}}_m(u) , \widehat{g}^{gr}_{m}(u))\leq\epsilon,\,\forall\,u\in [a_1,a_2].\]
\end{cor}
\textbf{Proof:} Follows from Proposition \ref{grp6}.
\begin{cor}
Let the assumptions of Corollary \ref{grco1} be hold. Then the sequence of granular inverse $F$-transforms $\{\widehat{g}^{gr}_m(u)\}$ uniformly
converges to fuzzy function $\widehat{{g}}$ .
\end{cor}
\textbf{Proof:} Follows from Corollary \ref{grco1}.\\\\
The following result describes how to estimate the difference between any two approximations of a given fuzzy function by the granular inverse $F$-transforms based on different sets of basic functions. As can be observed, it depends on the original fuzzy function's smoothness behavior, as defined by its modulus of continuity.
\begin{pro}\label{grp7} 
Let $P_1,...,P_m$, $P'_1,...,P'_m,\,m\geq 3$, be basic functions that form different uniform fuzzy partitions of $[a_1,a_2]$ and $\widehat{g}:[a_1,a_2]\rightarrow E^1$ be the gr-continuous fuzzy function. Further, let $\widehat{g}^{gr}_m,\widehat{g}^{gr'}_m$ be the granular inverse $F$-transforms of
$\widehat{g}$ w.r.t. different sets of basic functions. Then
\[d^{gr}(\widehat{g}^{gr}_m(u),\widehat{g}^{gr'}_m(u))\leq 2 \omega(2h,\widehat{g}), \,\forall\,u\in [a_1,a_2],\]
where $h=\frac{a_2-a_1}{m-1}$ and $2 \omega(2h,\widehat{g})$ is the modulus of continuity of $\widehat{g}$ on $[a_1,a_2]$.
\end{pro} 
\textbf{Proof:} Follows from Proposition \ref{grp3} and Definition \ref{grgrift}.\\\\
The following is towards the graphical representation of horizontal membership functions and level sets of the granular $F$-transform and the granular inverse $F$-transform associated with a fuzzy function are presented.
\begin{exa}\label{grexa}
\end{exa}
Let $P_1 = (0, 0, 1), P_2= (0, 1, 2), P_3 = (1, 2, 3), P_4 = (2, 3, 3)$ be the triangular fuzzy numbers that form a fuzzy partition of $[0,3]$. To calculate the components of the granular $F$-transform and the granular inverse $F$-transform, we assume a fuzzy function $\widehat{g}(u)=(\frac{u^3}{3},\frac{u^3}{3}+u+3,\frac{2u^3}{3}+4)$ and its horizontal membership function is given by 
\begin{eqnarray*}
\mathbb{K}(\widehat{g}(u))&=&g^{gr}(u,\alpha,\mu_g)\\
&=&\frac{u^3}{3}+\alpha(u+3)+\mu_g(1-\alpha)(\frac{u^3}{3}+4),\,u\in[0,3],\alpha,\mu_g\in\{0,0.5,1\}.
\end{eqnarray*}
Next, the horizontal membership functions of the components of the granular $F$-transform and the granular inverse $F$-transform associated with the fuzzy function $\widehat{g}(u)$ corresponding to different values of $\alpha,\mu_g$, are presented in Figure \ref{grfig:0}. In which, at $\alpha=1$, the fuzzy function $\widehat{g}(u)$ shows crisp behavior and the granular $F$-transform and the granular inverse $F$-transform become $F$-transform and inverse $F$-transform as given in \cite{per}, respectively. Also, the $\alpha(=0)$-level sets of the fuzzy function $\widehat{g}(u)$, granular $F$-transform and granular inverse $F$-transform are given in Figure \ref{grfig:00}.
\begin{figure}
    \centering
    \includegraphics[scale=0.5]{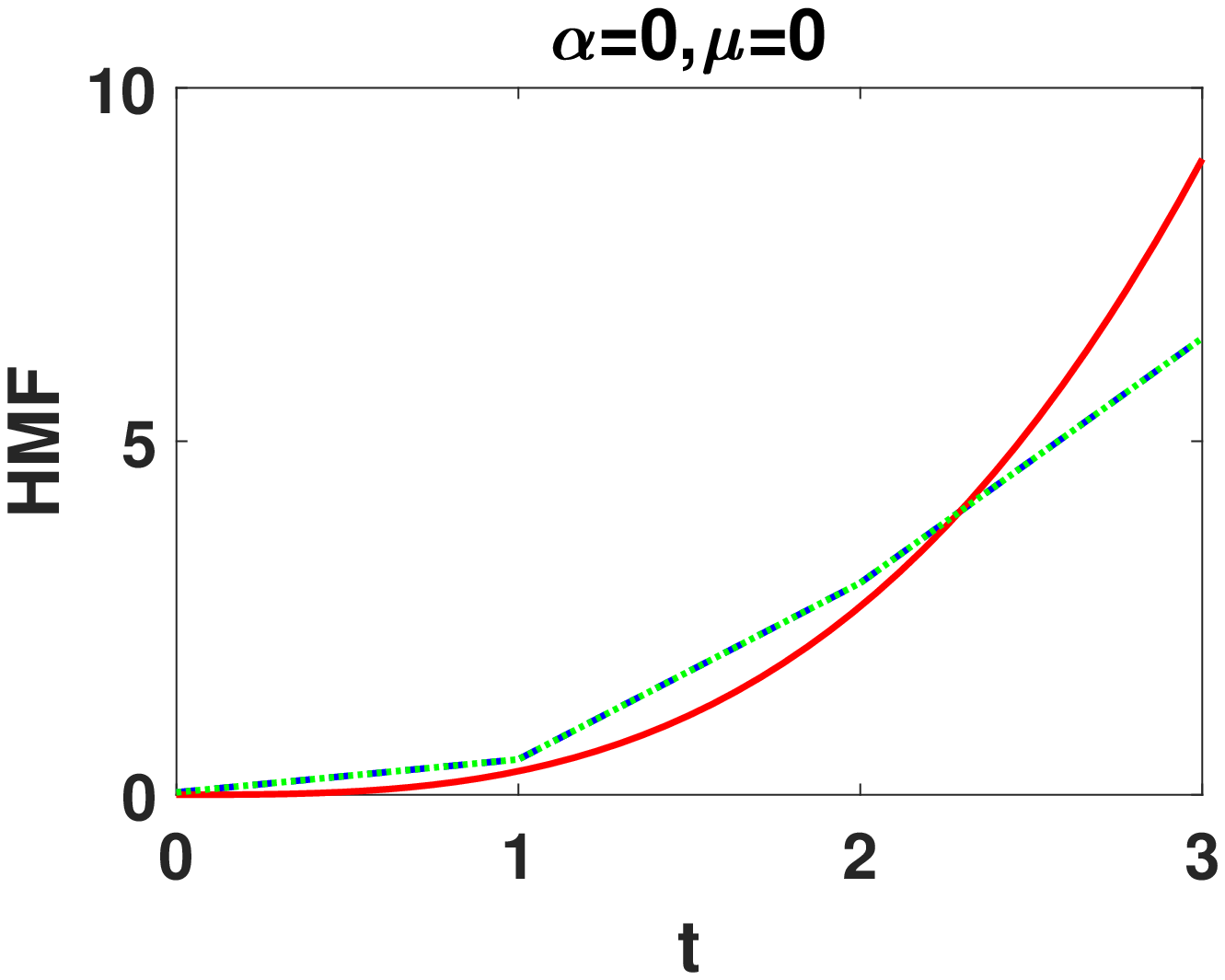} 
    \includegraphics[scale=0.5]{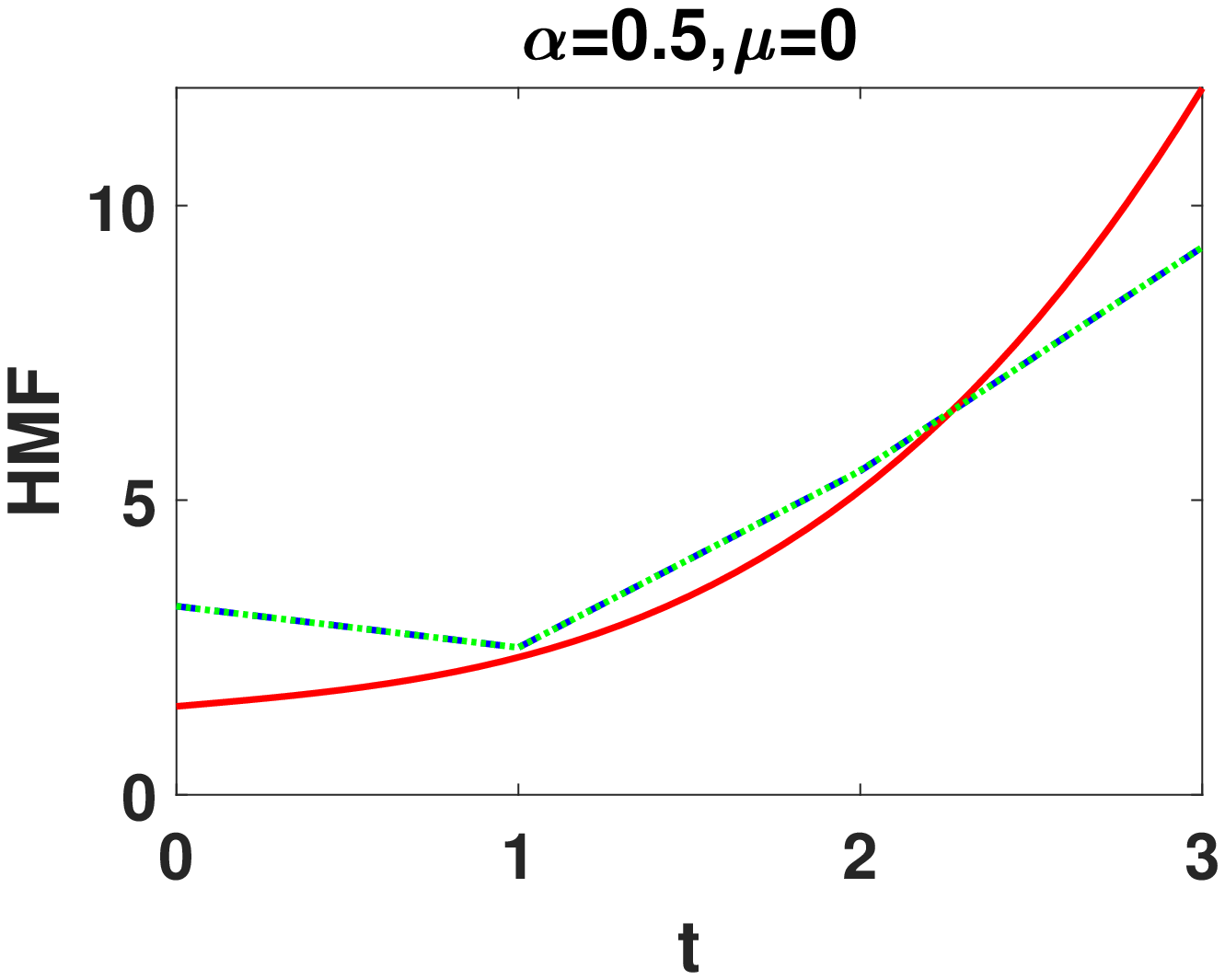}
      \includegraphics[scale=0.5]{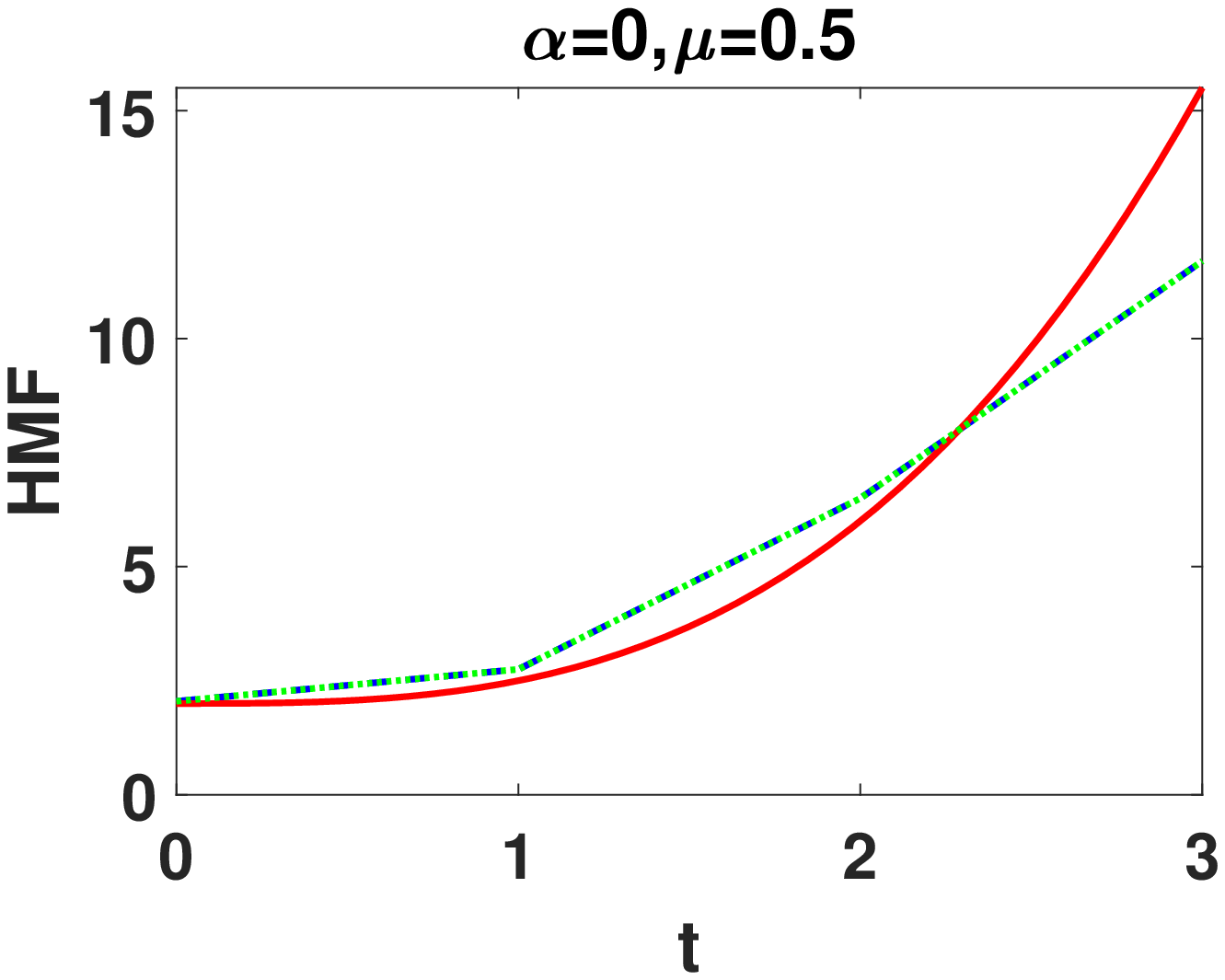}
        \includegraphics[scale=0.5]{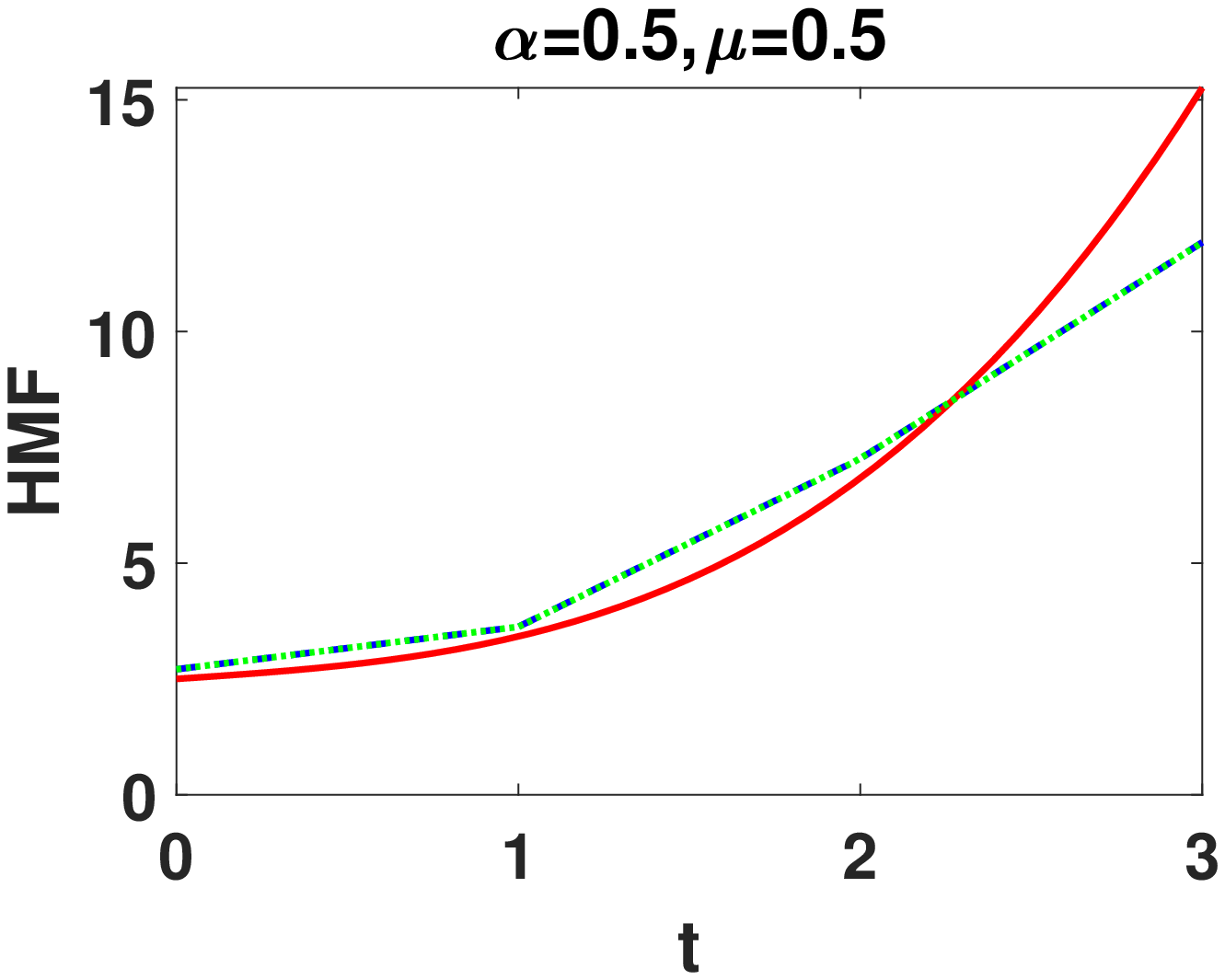}
       \includegraphics[scale=0.5]{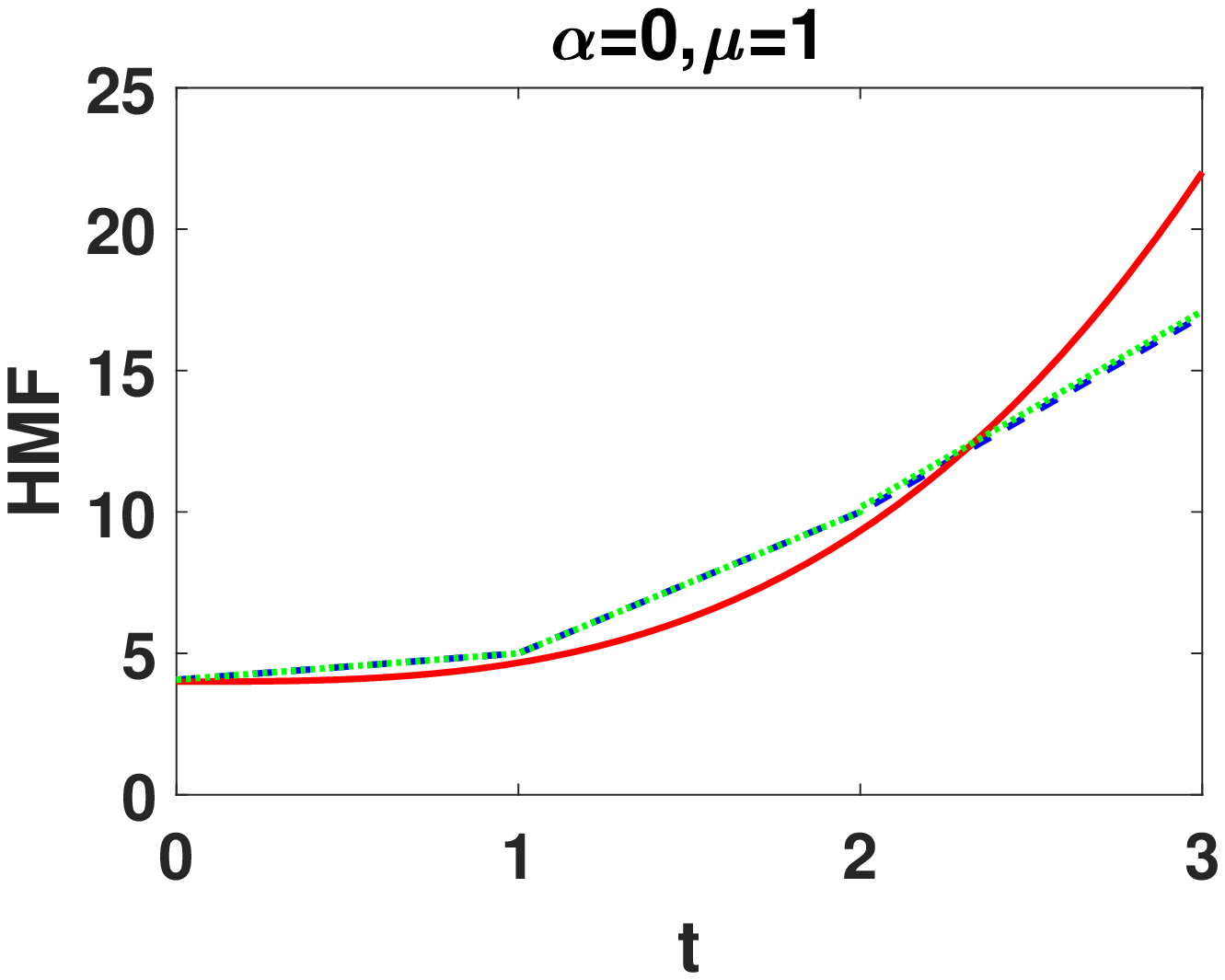}
         \includegraphics[scale=0.5]{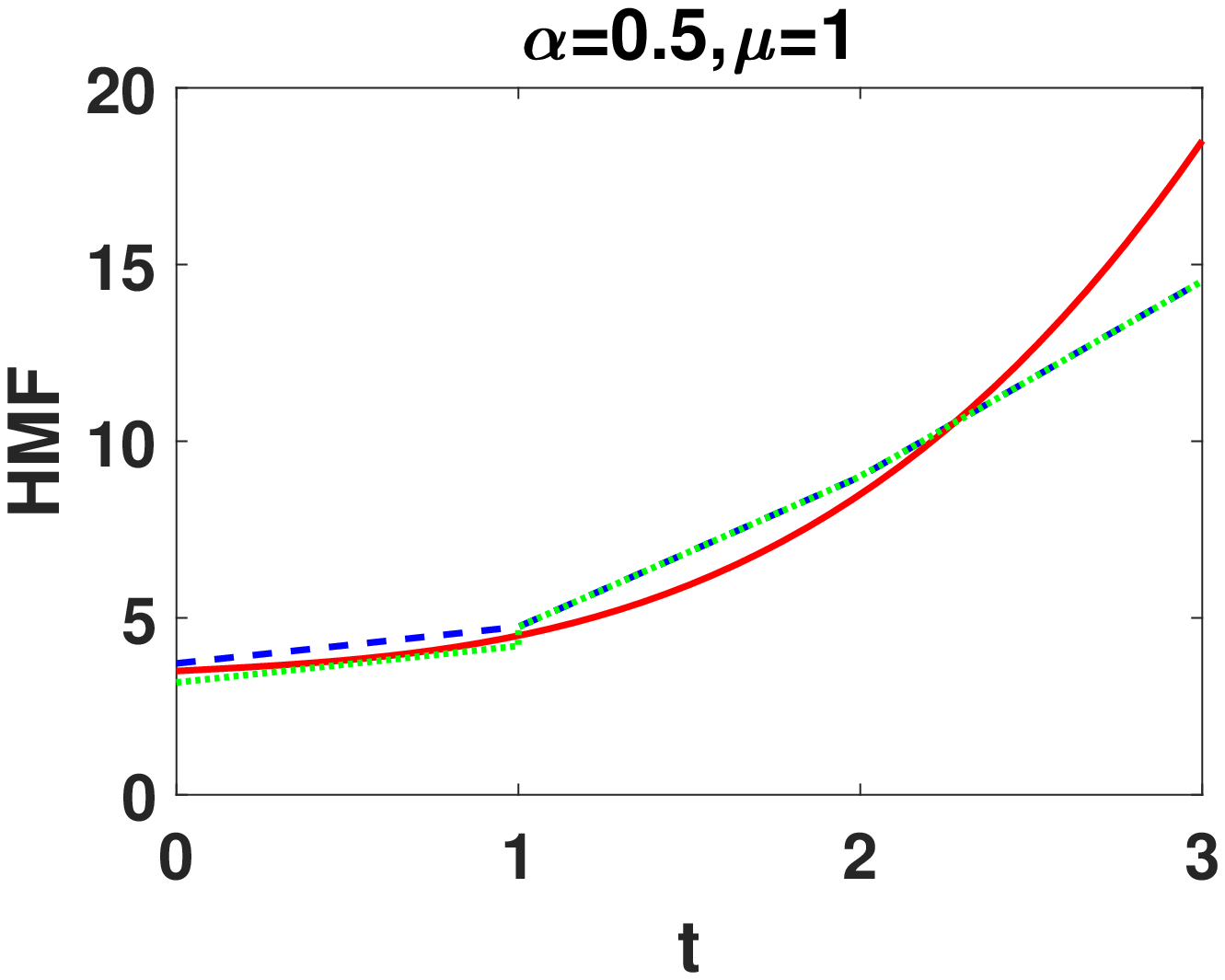}
     \includegraphics[scale=0.5]{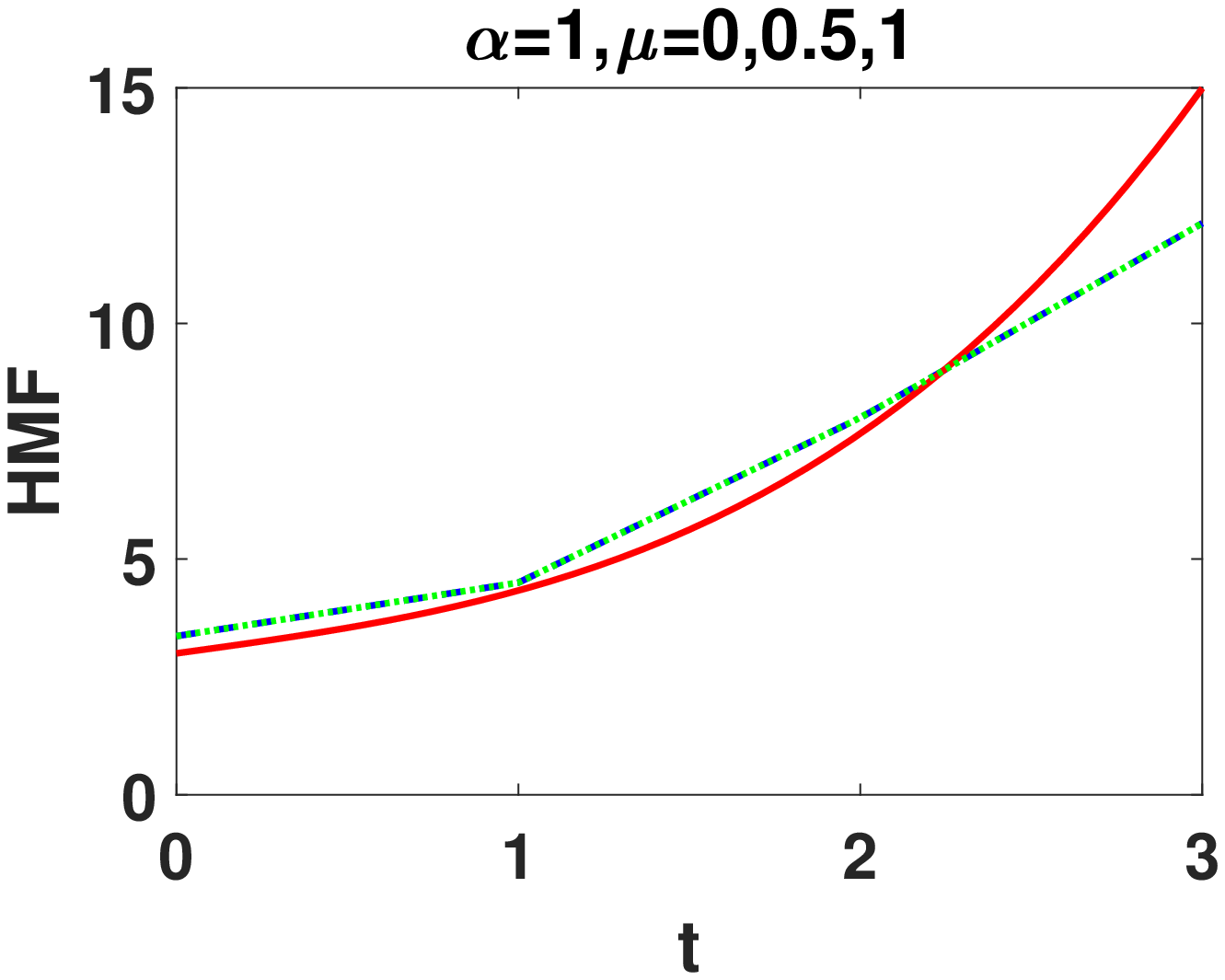}
     \caption{Red, blue, green curves show horizontal membership function (HMF) of fuzzy function $\widehat{f}(t)$, granular $F$-transform and granular inverse $F$-transform presented in Example \ref{exa} corresponding to different values of $\alpha,\mu$, respectively.}
    \label{fig:0}
\end{figure}
\begin{figure}
    \centering
    \includegraphics[scale=0.55]{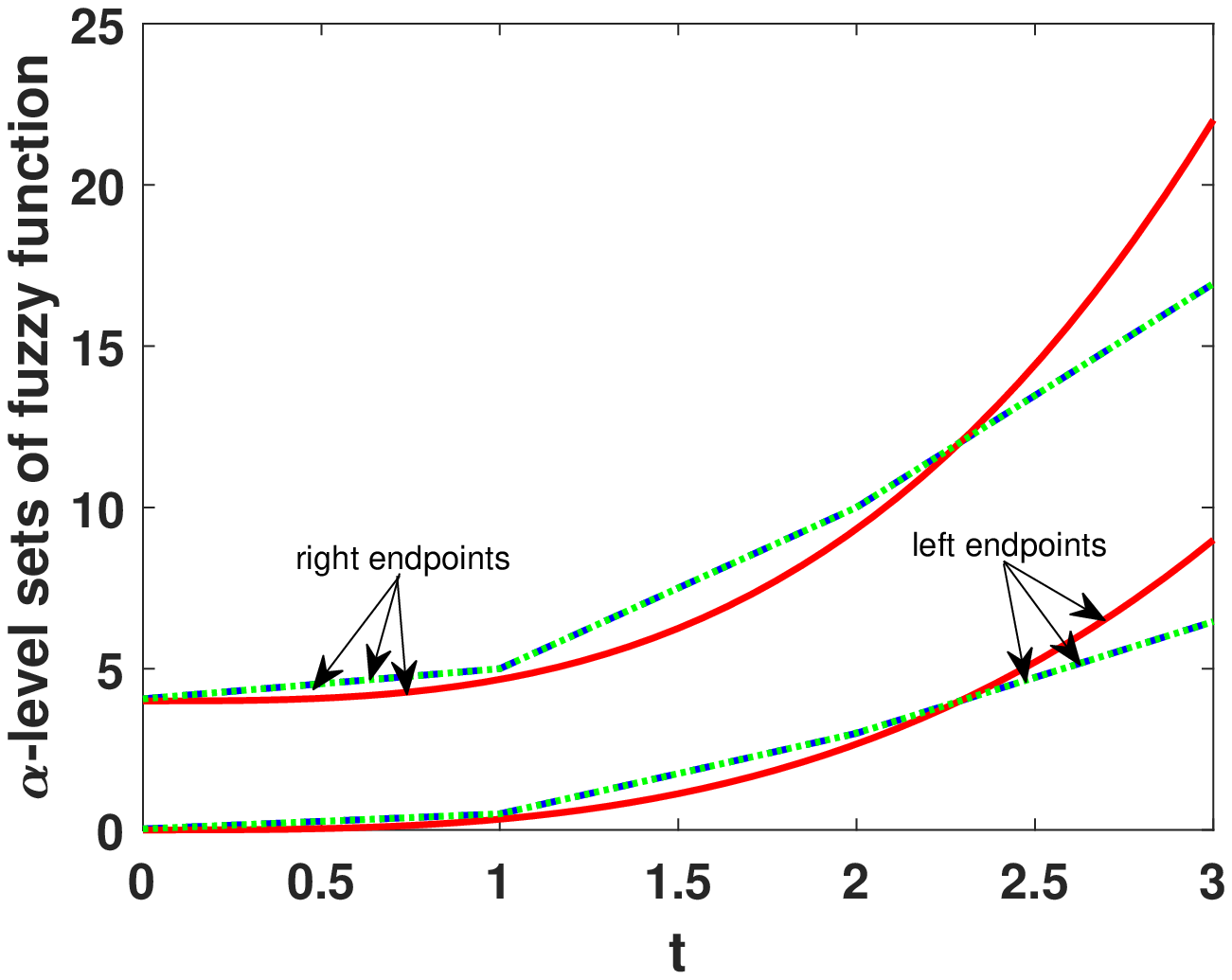} 
    \caption{Red, blue, green curves show $\alpha(=0)$-level sets of of fuzzy functions $\widehat{f}(t)$, granular $F$-transform and granular inverse $F$-transform presented in Example \ref{exa}, respectively. }
    \label{fig:00}
\end{figure}
\section{Fuzzy prey-predator model}
In this section, we formulate a fuzzy prey-predator model and study the dynamic behavior of the same. This section is divided into two subsections; the first is towards the formulation of the model, and the latter is towards its dynamical behavior.
\subsection{Formulation of fuzzy prey-predator model}
In this subsection, we present a fuzzy prey-predator model in which one predator team interacts with two teams of prey. In the presence of a predator, prey groups assist one another, but in the absence of a predator, they compete. We consider two teams of prey with densities $\widehat{p}(u)$ and $\widehat{q}(u)$, interacting with one team of predators with density $\widehat{r}(u)$ in a fuzzy environment, respectively, which is chiefly motivated from a crisp model given in \cite{elet}. The proposed fuzzy prey-predator model is as follows:
\begin{eqnarray}\label{gr1}
\begin{array}{ll}
{\mathcal{D}^{gr}\widehat{p}(u)}=\widehat{a_1}\widehat{p}(u)(1-\widehat{p}(u))-\widehat{p}(u)\widehat{r}(u)+\widehat{p}(u)\widehat{q}(u)\widehat{r}(u)=\widehat{g_1}(u,\widehat{p},\widehat{q},\widehat{r}),\\ 
{\mathcal{D}^{gr}\widehat{q}(u)}=\widehat{a_2}\widehat{q}(u)(1-\widehat{q}(u))-\widehat{q}(u)\widehat{r}(u)+\widehat{p}(u)\widehat{q}(u)\widehat{r}(u)=\widehat{g_2}(u,\widehat{p},\widehat{q},\widehat{r}),\\ 
{\mathcal{D}^{gr}\widehat{r}(u)}= -\widehat{a_3}\widehat{r}^2(u)+\widehat{a_4}\widehat{p}(u)\widehat{r}(u)+\widehat{a_5}\widehat{q}(u)\widehat{r}(u)=\widehat{g_3}(u,\widehat{p},\widehat{q},\widehat{r}),\\
\widehat{p}(0)=\widehat{p}_0,\,\widehat{q}(0)=\widehat{q}_0 ,\,\widehat{r}(0)=\widehat{r}_0,
\end{array}
\end{eqnarray}
where $\widehat{p},\widehat{q},\widehat{r}:[0,1]\subseteq \mathbb{R}\rightarrow E^1$ are gr-differentiable fuzzy functions and $\widehat{a_1},\widehat{a_2},\widehat{a_3},\widehat{a_4},\widehat{a_5}$, $\widehat{p}_0,\widehat{q}_0,\widehat{r}_0$ are positive fuzzy numbers. Based on
Remark \ref{grr1} and Proposition \ref{grgrdif1}, the system (\ref{gr1}) can be written as
\begin{eqnarray}\label{gr2}
\nonumber\mathbb{K}(\mathcal{D}^{gr}\widehat{p}(u))&=&{\frac{\partial}{\partial t}\mathbb{K}(\widehat{p}(u))}=\mathbb{K}(\widehat{a_1})\mathbb{K}(\widehat{p}(u))({1}-\mathbb{K}(\widehat{p}))-\mathbb{K}(\widehat{p}(u))\mathbb{K}(\widehat{r}(u))+\\
\nonumber&&\mathbb{K}(\widehat{p}(u))\mathbb{K}(\widehat{q}(u))\mathbb{K}(\widehat{r}(u))=\mathbb{K}(\widehat{g_1}(u,\mathbb{K}(\widehat{p}),\mathbb{K}(\widehat{q}),\mathbb{K}(\widehat{r}))),\\ 
\nonumber\mathbb{K}(\mathcal{D}^{gr}\widehat{q}(u))&=&{\frac{\partial}{\partial t}\mathbb{K}(\widehat{q}(u))}=\mathbb{K}(\widehat{a_2})\mathbb{K}(\widehat{q}(u))({1}-\mathbb{K}(\widehat{q}(u)))-\mathbb{K}(\widehat{q}(u))\mathbb{K}(\widehat{r}(u))+\\
\nonumber&&\mathbb{K}(\widehat{p}(u))\mathbb{K}(\widehat{q}(u))\mathbb{K}(\widehat{r}(u))=\mathbb{K}(\widehat{g_2}(\mathbb{K}(u,\widehat{p}),\mathbb{K}(\widehat{q}),\mathbb{K}(\widehat{r}))),\\ 
\nonumber\mathbb{K}(\mathcal{D}^{gr}\widehat{r}(u))&=&{\frac{\partial}{\partial t}\mathbb{K}(\widehat{r}(u))}= -\mathbb{K}(\widehat{a_3})\mathbb{K}(\widehat{r}(u))^2+\mathbb{K}(\widehat{a_4})\mathbb{K}(\widehat{p}(u))\mathbb{K}(\widehat{r}(u))+\\
\nonumber&&\mathbb{K}(\widehat{a_5})\mathbb{K}(\widehat{q}(u))\mathbb{K}(\widehat{r}(u))=\mathbb{K}(\widehat{g_3}(u,\mathbb{K}(\widehat{p}),\mathbb{K}(\widehat{q}),\mathbb{K}(\widehat{r}))),\\
\mathbb{K}(\widehat{p}(0))&=&\mathbb{K}(\widehat{p}_0),\,\mathbb{K}(\widehat{q}(0))=\mathbb{K}(\widehat{q}_0)
 ,\,\mathbb{K}(\widehat{r}(0))=\mathbb{K}(\widehat{r}_0).
\end{eqnarray}
The system (\ref{gr2}) can be given as
\begin{eqnarray}\label{gr2a}
\nonumber\frac{\partial}{\partial t}p^{gr}(u,\alpha,\mu_p)&=&a_1^{gr}(\alpha,\mu_{a_1})p^{gr}(u,\alpha,\mu_p)({1}-p^{gr}(u,\alpha,\mu_p))-p^{gr}(u,\alpha,\mu_p)r^{gr}(u,\alpha,\mu_r)\\
\nonumber&&+p^{gr}(u,\alpha,\mu_p)q^{gr}(u,\alpha,\mu_q)r^{gr}(u,\alpha,\mu_r),\\
\nonumber\frac{\partial}{\partial t}q^{gr}(u,\alpha,\mu_q)&=&a_2^{gr}(\alpha,\mu_{a_2})q^{gr}(u,\alpha,\mu_q)({1}-q^{gr}(u,\alpha,\mu_q))-q^{gr}(u,\alpha,\mu_q)r^{gr}(u,\alpha,\mu_r)\\
\nonumber&&+p^{gr}(u,\alpha,\mu_p)q^{gr}(u,\alpha,\mu_q)r^{gr}(u,\alpha,\mu_r),\\
\nonumber\frac{\partial}{\partial t}r^{gr}(u,\alpha,\mu_r)&=& -a_2^{gr}(\alpha,\mu_{a_3})r^{2,gr}(u,\alpha,\mu_r)+d^{gr}(\alpha,\mu_{a_4})p^{gr}(u,\alpha,\mu_p)r^{gr}(u,\alpha,\mu_r)+\\
&&a_5^{gr}(\alpha,\mu_{a_5})q^{gr}(u,\alpha,\mu_q)r^{gr}(u,\alpha,\mu_r),
\end{eqnarray}
$p^{gr}(0,\alpha,\mu_p)=p_0^{gr}(\alpha,\mu_{p_0}),q^{gr}(0,\alpha,\mu_q)=q_0^{gr}(\alpha,\mu_{q_0}),r^{gr}(0,\alpha,\mu_r)=r_0^{gr}(\alpha,\mu_{r_0}).$
\subsection{Dynamical behavior of the fuzzy prey-predator model}
\begin{figure}
  \centering
\includegraphics[scale=0.5]{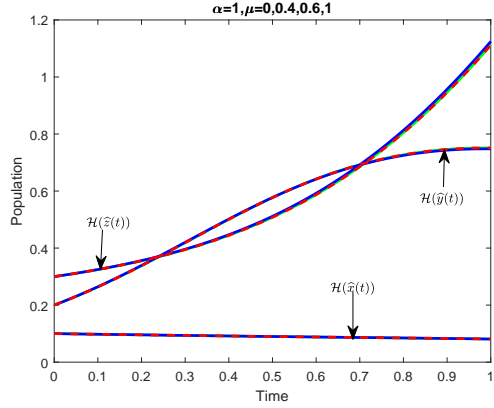}
  \caption{Variation of preys and predator populations against the time for the system (\ref{gr2a}) for different $\alpha,\mu$. Blue, green, red curves show population densities corresponding to exact, FT-Euler mid-point and Euler method, respectively.}
  \label{grfig:1}
\end{figure}

In this subsection, we investigate the equilibrium points and their stability of the proposed fuzzy prey-predator model. Now, we initiate with the following.
\begin{def1}\label{greq}
 A point $(\widehat{p_e},\widehat{q}_e,\widehat{r}_e)$ is called the {\bf fuzzy equilibrium point} of the system (\ref{gr1}) if
 \begin{eqnarray*}
 \mathcal{D}^{gr}\widehat{p}(u)&=&\widehat{g_1}(\widehat{p_e},\widehat{q}_e,\widehat{r}_e)=0,\\
 \mathcal{D}^{gr}\widehat{q}(u)&=&\widehat{g_2}(\widehat{p_e},\widehat{q}_e,\widehat{r}_e)=0,\\
 \mathcal{D}^{gr}\widehat{r}(u)&=&\widehat{g_3}(\widehat{p_e},\widehat{q}_e,\widehat{r}_e)=0.
 \end{eqnarray*}
\end{def1}
\begin{rem}\label{grremeq}
It can be easily seen that $(\widehat{p_e},\widehat{q}_e,\widehat{r}_e)$ is a fuzzy equilibrium point of the system (\ref{gr1}) iff $(\mathbb{K}(\widehat{p_e}),\mathbb{K}(\widehat{q}_e),\mathbb{K}(\widehat{r}_e))$ is an equilibrium point of the system (\ref{gr2}).
\end{rem}
 After algebraic calculation for the system (\ref{gr1}), we get several fuzzy equilibrium points $(\widehat{p},\widehat{q},\widehat{r})$, e.g.
\begin{itemize}
\item[(i)] $\widehat{E}_0({0},{0},{0})$,
\item[(ii)] $\widehat{E}_1({1},{0},{0})$,
\item[(iii)] $\widehat{E}_2({0},{1},{0})$,
\item[(iv)] $\widehat{E}_3({1},{1},{0})$,
\item[(v)] $\widehat{E}_4({0},\widehat{q_{e_4}},\widehat{r_{e_4}})$, where
\begin{eqnarray*}\widehat{q_{e_4}}&=&\frac{\widehat{a_2}\widehat{a_3}}{\widehat{a_2}\widehat{a_3}+\widehat{a_5}},\\
\widehat{r_{e_4}}&=&\frac{\widehat{a_2}\widehat{a_5}}{\widehat{a_2}\widehat{a_3}+\widehat{a_5}},
\end{eqnarray*}
\item[(vi)] $\widehat{E}_5(\widehat{p_{e_5}},{0},\widehat{r_{e_5}})$, where 
\begin{eqnarray*}
\widehat{p_{e_5}}&=&\frac{\widehat{a_1}\widehat{a_3}}{\widehat{a_1}\widehat{a_3}+\widehat{a_4}},\\
\widehat{r_{e_5}}&=&\frac{\widehat{a_1}\widehat{a_4}}{\widehat{a_1}\widehat{a_3}+\widehat{a_4}},
\end{eqnarray*}
\item[(vii)] $\widehat{E}_6({1},{1},\widehat{r_{e_6}})$, where 
\[\widehat{r_{e_6}}=\frac{\widehat{a_4}+\widehat{a_5}}{\widehat{a_3}},\,\mbox{and}\]
\item[(viii)] $\widehat{E}_7(\widehat{p_{e_7}},\widehat{q_{e_7}},\widehat{r_{e_7}})$, where
\begin{eqnarray*}
\widehat{p_{e_7}}&=&\frac{\widehat{a_2}\widehat{a_3}+\widehat{a_5}\left({1}-\sqrt{\frac{\widehat{a_2}}{\widehat{a_1}}}\right)}{\widehat{a_5}+\widehat{a_4}\sqrt{\frac{\widehat{a_2}}{\widehat{a_1}}}},\\
\widehat{q_{e_7}}&=&\frac{\widehat{a_3}\sqrt{\widehat{a_1}\widehat{a_2}}-\widehat{a_4}\left({1}-\sqrt{\frac{\widehat{a_2}}{\widehat{a_1}}}\right)}{\widehat{a_5}+\widehat{a_4}\sqrt{\frac{\widehat{a_2}}{\widehat{a_1}}}},\\
\widehat{r_{e_7}}&=&\sqrt{\widehat{a_1}\widehat{a_2}}.
\end{eqnarray*}
\end{itemize}
The existence of ${E}_0,{E}_1, {E}_2,{E}_3,{E}_4,{E}_5$ and ${E}_6$ are obvious. Now, we are concentrating on the existence of fuzzy equilibrium point $\widehat{E}_7$. Therefore the fuzzy equilibrium point ${E}_7$ exists if 
\begin{eqnarray*}
&&\widehat{a_3}\sqrt{\widehat{a_1}\widehat{a_2}}\leq\widehat{a_4}+\widehat{a_5},\\
&&\widehat{a_1}\widehat{a_3}+\widehat{a_4}>\widehat{a_4}\sqrt{\frac{\widehat{a_2}}{\widehat{a_1}}},\\
&&\widehat{a_2}\widehat{a_3}+\widehat{a_5}>\widehat{a_5}\sqrt{\frac{\widehat{a_1}}{\widehat{a_2}}}.
\end{eqnarray*}
From Remark \ref{grremeq}, $(\mathbb{K}(\widehat{p}),\mathbb{K}(\widehat{q}),\mathbb{K}(\widehat{r}))$ is an equilibrium point of the system (\ref{gr2}), e.g. 
\begin{itemize}
\item[(i)] $E_0({0},{0},{0})$,
\item[(ii)] $E_1({1},{0},{0})$,
\item[(iii)] $E_2({0},{1},{0})$,
\item[(iv)] $E_3({1},{1},{0})$,
\item[(v)] $E_4({0},\mathbb{K}(\widehat{q_{e_4}}),\mathbb{K}(\widehat{r_{e_4}}))$, where
\begin{eqnarray*}
\mathbb{K}(\widehat{q_{e_4}})&=&\frac{\mathbb{K}(\widehat{a_2})\mathbb{K}(\widehat{a_3})}{\mathbb{K}(\widehat{a_2})\mathbb{K}(\widehat{a_3})+\mathbb{K}(\widehat{a_5})},\\
\mathbb{K}(\widehat{r_{e_4}})&=&\frac{\mathbb{K}(\widehat{a_2})\mathbb{K}(\widehat{a_5})}{\mathbb{K}(\widehat{a_2})\mathbb{K}(\widehat{a_3})+\mathbb{K}(\widehat{a_5})},
\end{eqnarray*}
\item[(vi)]$E_5(\mathbb{K}(\widehat{p_{e_5}}),{0},\mathbb{K}(\widehat{r_{e_5}}))$, where
\begin{eqnarray*}
\mathbb{K}(\widehat{p_{e_5}})&=&\frac{\mathbb{K}(\widehat{a_1})\mathbb{K}(\widehat{a_3})}{\mathbb{K}(\widehat{a_1})\mathbb{K}(\widehat{a_3})+\mathbb{K}(\widehat{a_4})},\\
\mathbb{K}(\widehat{r_{e_5}})&=&\frac{\mathbb{K}(\widehat{a_1})\mathbb{K}(\widehat{a_4})}{\mathbb{K}(\widehat{a_1})\mathbb{K}(\widehat{a_3})+\mathbb{K}(\widehat{a_4})}
\end{eqnarray*}
\item[(vii)] $E_6({1},{1},\mathbb{K}(\widehat{r_{e_6}}))$, where 
\[\mathbb{K}(\widehat{r_{e_6}}=\frac{\mathbb{K}(\widehat{a_4})+\mathbb{K}(\widehat{a_5})}{\mathbb{K}(\widehat{a_3})},\,\mbox{and}\]
\item[(viii)] $E_7(\mathbb{K}(\widehat{p_{e_7}}),\mathbb{K}(\widehat{q_{e_7}}),\mathbb{K}(\widehat{r_{e_7}}))$, where 
\begin{eqnarray*}
\mathbb{K}(\widehat{p_{e_7}})&=&\frac{\mathbb{K}(\widehat{a_2})\mathbb{K}(\widehat{a_3})+\mathbb{K}(\widehat{a_5})\left({1}-\sqrt{\frac{\mathbb{K}(\widehat{a_2})}{\mathbb{K}(\widehat{a_1})}}\right)}{\mathbb{K}(\widehat{a_5})+\mathbb{K}(\widehat{a_4})\sqrt{\frac{\mathbb{K}(\widehat{a_2})}{\mathbb{K}(\widehat{a_1})}}},\\
\mathbb{K}(\widehat{q_{e_7}})&=&\frac{\mathbb{K}(\widehat{a_3})\sqrt{\mathbb{K}(\widehat{a_1})\mathbb{K}(\widehat{a_2})}-\mathbb{K}(\widehat{a_4})\left({1}-\sqrt{\frac{\mathbb{K}(\widehat{a_2})}{\mathbb{K}(\widehat{a_1})}}\right)}{\mathbb{K}(\widehat{a_5})+\mathbb{K}(\widehat{a_4})\sqrt{\frac{\mathbb{K}(\widehat{a_2})}{\mathbb{K}(\widehat{a_1})}}},\\
\mathbb{K}(\widehat{r_{e_7}})&=&\sqrt{\mathbb{K}(\widehat{a_1})\mathbb{K}(\widehat{a_2})},
\end{eqnarray*}
\end{itemize}
are the equilibrium points of the system (\ref{gr2}).\\\\ 
Below, we discuss the local and global fuzzy stability of the fuzzy equilibrium points corresponding to the fuzzy eigenvalues, respectively. Before stating the next, we introduce the following.
\begin{def1}\label{grstable} Let $\widehat{\lambda}_i,\,i=1,2,...,m$, be a fuzzy eigenvalue of the fuzzy matrix $\widehat{M}=[\widehat{a_1}_{ij}]_{m\times m}$.
Then the fuzzy equilibrium points of the system of $\mathcal{D}^{gr}\widehat{p}(u)=\widehat{M}\widehat{p}(u)$ are called {\bf fuzzy stable} and {\bf fuzzy unstable} iff $Re(\widehat{\lambda}_i)<0$ and $Re(\widehat{\lambda}_i)>0$, respectively.
\end{def1}
\begin{def1}\label{grlocal}
A fuzzy equilibrium point $\widehat{p_e}$ of $\mathcal{D}^{gr}\widehat{p}(u)=\widehat{M}\widehat{p}(u)$ is called {\bf locally asymptotically fuzzy stable} if all eigen values of the corresponding variational fuzzy matrix have negative real parts.
\end{def1}
\begin{def1}\label{grglobal} A fuzzy equilibrium point $\widehat{p_e}$ of $\mathcal{D}^{gr}\widehat{p}(u)=\widehat{M}\widehat{p}(u)$ is called {\bf globally asymptotically fuzzy stable} if there exists a $gr$-continuous and differentiable function $\widehat{U}:D\subseteq E^m\rightarrow E^1$ such that 
\begin{itemize}
  \item[(i)] $\widehat{U}(\widehat{p_e})=0$;
    \item[(ii)] $\widehat{U}(\widehat{p})>0,\,\forall\,\widehat{p}\in D-\{\widehat{p_e}\}$;
      \item[(iii)] $\widehat{U}(\widehat{p})$ is radially unbounded; and
        \item[(iv)] $\mathcal{D}^{gr}\widehat{U}(\widehat{p})<0,\,\forall\,\widehat{p}\in D-\{\widehat{p_e}\}$.
\end{itemize}
Also, the fuzzy function $\widehat{U}$ is called {\bf fuzzy Lyapunov function}.
\end{def1}
\begin{thm}\label{grlogo}
A fuzzy equilibrium point $\widehat{p_e}$ of $\mathcal{D}^{gr}\widehat{p}(u)=\widehat{M}\widehat{p}(u)$ is locally or globally asymptotically fuzzy stable iff the equilibrium point $\mathbb{K}(\widehat{p_e})$ of $\mathbb{K}(\mathcal{D}^{gr}\widehat{p}(u))=\mathbb{K}(\widehat{M})\mathbb{K}(\widehat{p}(u))$ is locally or globally asymptotically stable, respectively.
\end{thm}
\textbf{Proof:} Let $\widehat{p_e}$ be a fuzzy equilibrium point of $\mathcal{D}^{gr}\widehat{p}(u)=\widehat{M}\widehat{p}(u)$. Then from Remark \ref{grremeq}, $\mathbb{K}(\widehat{p_e})$ is an equilibrium point of $\mathbb{K}(\mathcal{D}^{gr}\widehat{p}(u))=\mathbb{K}(\widehat{M})\mathbb{K}(\widehat{p}(u))$. Also, let $\mathbb{K}(\widehat{\lambda}_i),\,i=1,2,...,m$ be an eigen value of the corresponding variational matrix. Then the equilibrium point $\mathbb{K}(\widehat{p_e})$ is locally asymptotically stable iff $Re(\mathbb{K}(\widehat{\lambda}_i))<0$, i.e., $\mathbb{K}(Re(\widehat{\lambda}_i))<0$. Therefore from Remark \ref{grr1}, we have $Re(\widehat{\lambda}_i)<0$. Thus the fuzzy equilibrium point $\widehat{p_e}$ is locally asymptotically
fuzzy stable. Also, the equilibrium point $\mathbb{K}(\widehat{p_e})$ is globally asymptotically
stable if there exists a suitable Lyapunov function $\mathbb{K}({\widehat{U}})$ such that $\mathbb{K}({\widehat{U}}(\mathbb{K}(\widehat{p_e})))=0$, $\mathbb{K}({\widehat{U}}(\mathbb{K}(\widehat{p})))>0,\frac{\partial\mathbb{K}({\widehat{U}})}{\partial t}<0,\forall \,\mathbb{K}(\widehat{p})$ (except for the equilibrium point $\mathbb{K}(\widehat{p_e})$ and $\mathbb{K}({\widehat{U}}(\mathbb{K}(\widehat{p})))$ is radially unbounded. Therefore from Remark \ref{grr1}, the fuzzy function $\widehat{U}$ is the fuzzy Lyapunov function, i.e., satisfies all conditions given in Definition \ref{grglobal}. Thus the fuzzy equilibrium point $\widehat{p_e}$ is globally asymptotically
fuzzy stable and conversely. \\\\
The variational fuzzy matrix of the system (\ref{gr1}) is defined as
\begin{equation}\label{grm}
\widehat{M_1} = 
\begin{bmatrix}
\widehat{a_1}-2\widehat{a_1}\widehat{p}-\widehat{r}+\widehat{q}\widehat{r} & \widehat{p}\widehat{r} & -\widehat{p}+\widehat{p}\widehat{q} \\
\widehat{q}\widehat{r} & \widehat{a_2}-2\widehat{a_2}\widehat{q}-\widehat{r}+\widehat{p}\widehat{r} & -\widehat{q} +\widehat{p}\widehat{q}\\
\widehat{a_4}\widehat{r} & \widehat{a_5}\widehat{r} & -2\widehat{a_3}\widehat{r}+\widehat{a_4}\widehat{p}+\widehat{a_5}\widehat{q} 
\end{bmatrix}.
\end{equation}
Now, the variational matrix of the system (\ref{gr2}) is given by
\begin{equation}\label{grm1}
\mathbb{K}(\widehat{M_1}) = 
\begin{bmatrix}
\mathbb{K}(\widehat{a}_{11})& \mathbb{K}(\widehat{p})\mathbb{K}(\widehat{r}) & -\mathbb{K}(\widehat{p})+\mathbb{K}(\widehat{p})\mathbb{K}(\widehat{q}) \\
\mathbb{K}(\widehat{q})\mathbb{K}(\widehat{r} )& \mathbb{K}(\widehat{a}_{22}) & -\mathbb{K}(\widehat{q})+\mathbb{K}(\widehat{p})\mathbb{K}(\widehat{q}) \\
\mathbb{K}(\widehat{a_4})\mathbb{K}(\widehat{r}) & \mathbb{K}(\widehat{a_5})\mathbb{K}(\widehat{r}) & \mathbb{K}(\widehat{a}_{33})
\end{bmatrix},
\end{equation}
where
\[\mathbb{K}(\widehat{a}_{11})=\mathbb{K}(\widehat{a_1})-2\mathbb{K}(\widehat{a_1})\mathbb{K}(\widehat{p})-\mathbb{K}(\widehat{r})+\mathbb{K}(\widehat{q})\mathbb{K}(\widehat{r}), \]
\[\mathbb{K}(\widehat{a}_{22})=\mathbb{K}(\widehat{a_2})-2\mathbb{K}(\widehat{a_2})\mathbb{K}(\widehat{q})-\mathbb{K}(\widehat{r})+\mathbb{K}(\widehat{p})\mathbb{K}(\widehat{r} ),\]
\[\mathbb{K}(\widehat{a}_{33})=-2\mathbb{K}(\widehat{a_3})\mathbb{K}(\widehat{r})+\mathbb{K}(\widehat{a_4})\mathbb{K}(\widehat{p})+\mathbb{K}(\widehat{a_5})\mathbb{K}(\widehat{q}).
\] 
To check the stability of the fuzzy equilibrium points of the systems (\ref{gr1}), we define the characteristic equation of the matrix (\ref{grm1}) by $det(\mathbb{K}(\widehat{\lambda})I-\mathbb{K}(\widehat{M_1}))=0$.
\begin{itemize}
\item[(i)] The variational matrix (\ref{grm1}) corresponding to the equilibrium point $E_0(0,0,0)$ has eigenvalues
$\mathbb{K}(\widehat{\lambda})=0,\mathbb{K}(\widehat{a_1}),\mathbb{K}(\widehat{a_2})$. In which two eigenvalues are positive. Thus $E_0(0,0,0)$ is an unstable equilibrium point. From Remark \ref{grr1}, the variational fuzzy matrix (\ref{grm}) corresponding to the fuzzy equilibrium point $\widehat{E}_0(0,0,0)$ has fuzzy eigenvalues $\widehat{\lambda}=\widehat{0},\widehat{a_1},\widehat{a_2}$ and $\widehat{E}_0(0,0,0)$ is an unstable fuzzy equilibrium point.
\item[(ii)] The variational matrix (\ref{grm1}) corresponding to the equilibrium point $E_1(1,0,0)$ has
$\mathbb{K}(\widehat{\lambda})=-\mathbb{K}(\widehat{a_1}),\mathbb{K}(\widehat{a_2}),\mathbb{K}(\widehat{a_4})$. In which two eigenvalues are positive. Thus $E_1(1,0,0)$ is an unstable equilibrium point. From Remark \ref{grr1}, the variational fuzzy matrix (\ref{grm}) corresponding to the fuzzy equilibrium point $\widehat{E}_1(1,0,0)$ has fuzzy eigenvalues $\widehat{\lambda}=-\widehat{a_1},\widehat{a_2},\widehat{a_4}$ and $\widehat{E}_1(1,0,0)$ is an unstable fuzzy equilibrium point.
\item[(iii)] The variational matrix (\ref{grm1}) corresponding to the equilibrium point $E_2(0,1,0)$ has eigenvalues
$\mathbb{K}(\widehat{\lambda})=\mathbb{K}(\widehat{a_1}),-\mathbb{K}(\widehat{a_2}),\mathbb{K}(\widehat{a_5})$. In which two eigenvalues are positive. Thus $E_2(0,1,0)$ is an unstable equilibrium point. From Remark \ref{grr1}, the variational fuzzy matrix (\ref{grm}) corresponding to the fuzzy equilibrium point $\widehat{E}_2({0},{1},{0})$ has fuzzy eigenvalues $\widehat{\lambda}=\widehat{a_1},-\widehat{a_2},\widehat{a_5}$ and $\widehat{E}_2({0},{1},{0})$ is an unstable fuzzy equilibrium point.
\item[(iv)] The variational matrix (\ref{grm1}) corresponding to the equilibrium point $E_3(1,1,0)$ has eigenvalues
$\mathbb{K}(\widehat{\lambda})=-\mathbb{K}(\widehat{a_1}),-\mathbb{K}(\widehat{a_2}),\mathbb{K}(\widehat{a_4})+\mathbb{K}(\widehat{a_5})$. In which one eigenvalues are positive. Thus $E_3(1,1,0)$ is an unstable equilibrium point. From Remark \ref{grr1}, the variational fuzzy matrix (\ref{grm}) corresponding to the fuzzy equilibrium point $\widehat{E}_3(1,1,0)$ has fuzzy eigenvalues $\widehat{\lambda}=-\widehat{a_1},-\widehat{a_2},\widehat{a_4}+\widehat{a_5}$ and $\widehat{E}_3(1,1,0)$ is an unstable fuzzy equilibrium point.
\item[(v)] The equilibrium point $E_4(0,\mathbb{K}(\widehat{q_{e_4}}),\mathbb{K}(\widehat{r_{e_4}}))$ is locally asymptotically stable if 
\[\mathbb{K}(\widehat{a_1})<\frac{\mathbb{K}(\widehat{a_2})\mathbb{K}(\widehat{a_5})^2}{(\mathbb{K}(\widehat{a_2})\mathbb{K}(\widehat{a_3})+\mathbb{K}(\widehat{a_5}))^2}.\]
Therefore $\widehat{E}_4({0},\widehat{q_{e_4}},\widehat{r_{e_4}})$ is locally asymptotically fuzzy stable if \[\widehat{a_1}<\frac{\widehat{a_2}\widehat{a_5}^2}{(\widehat{a_2}\widehat{a_3}+\widehat{a_5})^2}.\]
\begin{exa}\label{grex1}
\vspace{-5mm}
\end{exa} Let $\widehat{a_1}=(0.01,0.02,0.03),\widehat{a_2}=(2,4,6),\widehat{a_3}=(0.1,0.2,0.3),\widehat{a_4}=(3,4,5),\widehat{a_5}=(1,2,3),\widehat{p}_0=0.1,\widehat{q}_0=0.2,\widehat{r}_0=0.3$, where $\widehat{a_1},\widehat{a_2},\widehat{a_3},\widehat{a_4},\widehat{a_5}$ are triangular fuzzy numbers. Now, the horizontal membership functions of the given triangular fuzzy numbers are given by
\begin{eqnarray*}
\mathbb{K}(\widehat{a_1})&=&a_1^{gr}(\alpha,\mu_{a_1})=0.01+0.01\alpha+\mu_{a_1}(0.02-0.02\alpha),\\
\mathbb{K}(\widehat{a_2})&=&a_2^{gr}(\alpha,\mu_{a_2})=2+2\alpha+\mu_{a_2}(4-4\alpha),\\
\mathbb{K}(\widehat{a_3})&=&a_3^{gr}(\alpha,\mu_{a_3})=0.1+0.1\alpha+\mu_{a_3}(0.2-0.2\alpha),\\
\mathbb{K}(\widehat{a_4})&=&a_4^{gr}(\alpha,\mu_{a_4})=3+\alpha+\mu_{a_4}(2-2\alpha),\\
\mathbb{K}(\widehat{a_5})&=&a_5^{gr}(\alpha,\mu_{a_5})=1+\alpha+\mu_{a_5}(2-2\alpha),\\
\mathbb{K}(\widehat{p}_0)&=&0.1,\mathbb{K}(\widehat{q}_0)=0.2,\mathbb{K}(\widehat{r}_0)=0.3.
\end{eqnarray*}
Further, we assume $\mu_p,\mu_q,\mu_r,\mu_{a_1},\mu_{a_2},\mu_{a_3},\mu_{a_4},\mu_{a_5}=\mu\in\{0,0.4,0.6,1\}$ and $\alpha\in\{0,0.5,1\}$.
\begin{table}[ht!]
 \begin{center}
    \begin{adjustbox}{max width=0.9\linewidth}
  \begin{tabular}{|p{.4cm}|p{3.1cm}| p{3.1cm}|p{3.1cm}|} 
  \hline
   $\mu$ & $\alpha=0$&  $\alpha=0.5$& $\alpha=1$\\
   \hline
   0&(0,0.1667,1.6667)&(0,0.2308,2.3077)&(0,0.2857,2.8571) \\
   \hline
    0.4&(0,0.2647,2.6471)&(0,0.2754,2.7536)&(0,0.2857,2.8571) \\
   \hline
    0.6&(0,0.3056,3.0556)&(0,0.2958,2.9578)&(0,0.2857,2.8571)\\
   \hline
    1&(0,0.375,3.75)&(0,0.3333,3.3333)&(0,0.2857,2.8571)\\
    \hline
  \end{tabular}
  \end{adjustbox}
   \caption{Equilibrium point $E_4(0,\mathbb{K}(\widehat{q_{e_4}}),\mathbb{K}(\widehat{r_{e_4}}))$ for different values of $\alpha,\mu$}
  \label{grtab:table1}
 \end{center}
\end{table}
For the given data set, we realize that the stability condition of $E_4$ is well satisfied and the system (\ref{gr2a}) have different equilibrium points corresponding to different values of $\alpha,\mu$, as shown in Table \ref{grtab:table1}. From Table \ref{grtab:table1}, we realize that for a fixed value of $\alpha$ the equilibrium point increases with increasing value of $\mu$ and at $\alpha=1$ remain same because the left and right interval of triangular fuzzy number coincides with each other. For $\mu=0,0.4$ the equilibrium point increases while for $\mu=0.6,1$ the equilibrium point decreases with increasing value of $\alpha$. Figure \ref{grfig:1} shows that for the given data in Example \ref{grex1}, the population density $\mathbb{K}(\widehat{p}(u))$ of prey eventually extinct while the population densities $\mathbb{K}(\widehat{q}(u))$ and $\mathbb{K}(\widehat{r}(u))$ of prey and predator persist, respectively and eventually get their steady states given in Table \ref{grtab:table1} corresponding to different values of $\alpha,\mu$ and become asymptotically stable. At $\alpha=1$, it shows crisp behavior.
\begin{figure}
 \centering
   \includegraphics[scale=0.75]{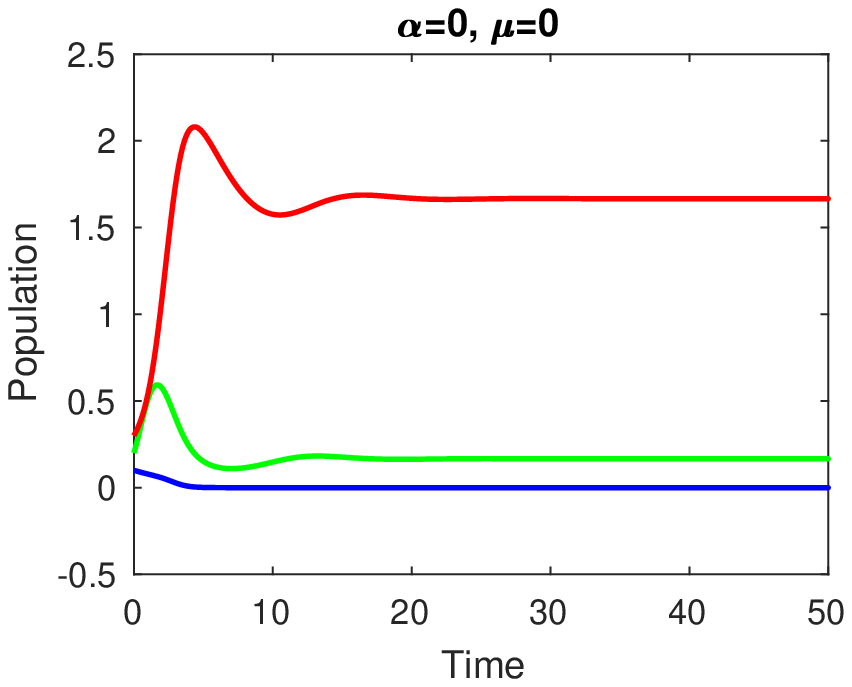} 
   \includegraphics[scale=0.75]{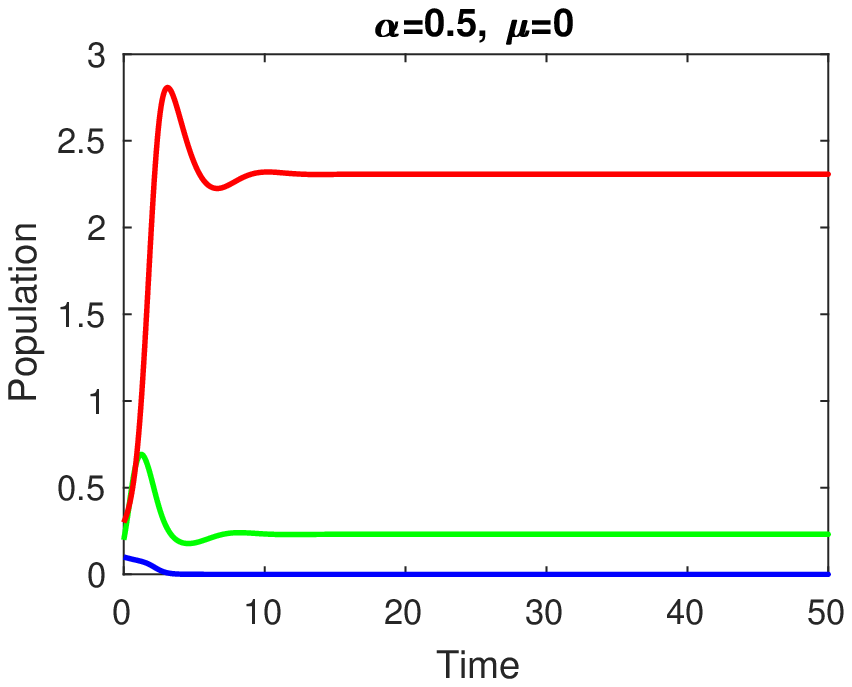}
  \includegraphics[scale= 0.75]{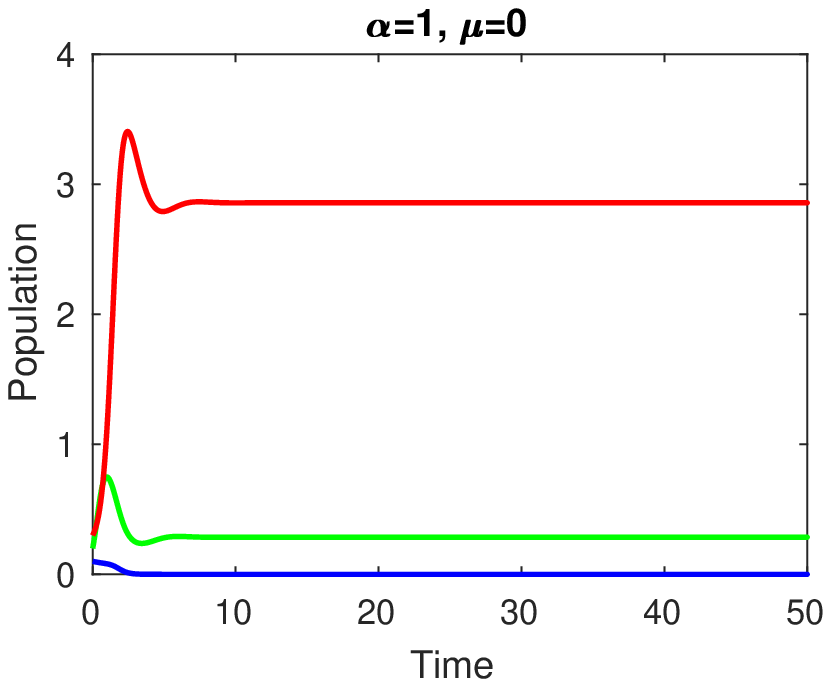}
 \includegraphics[scale=0.75]{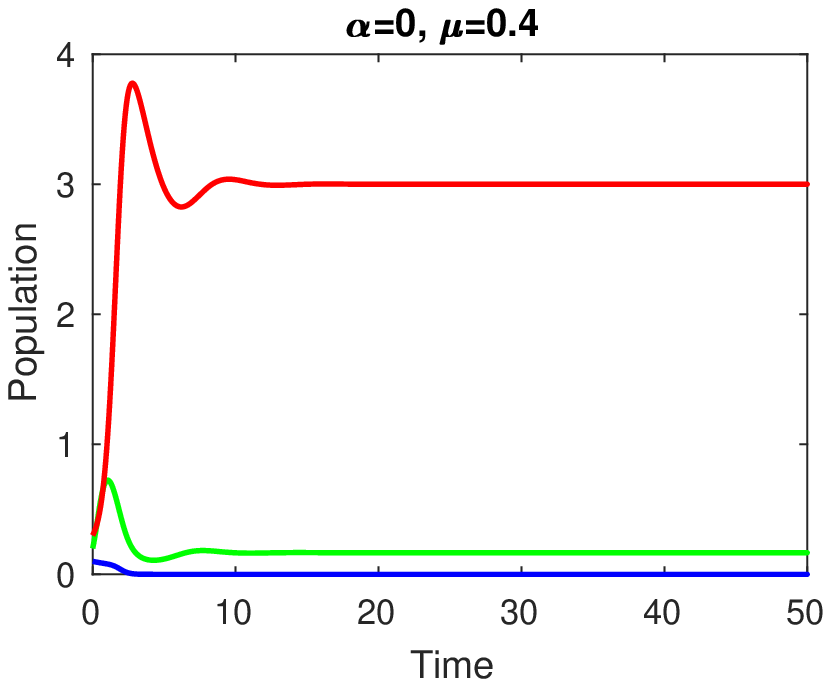}
   \includegraphics[scale=0.75]{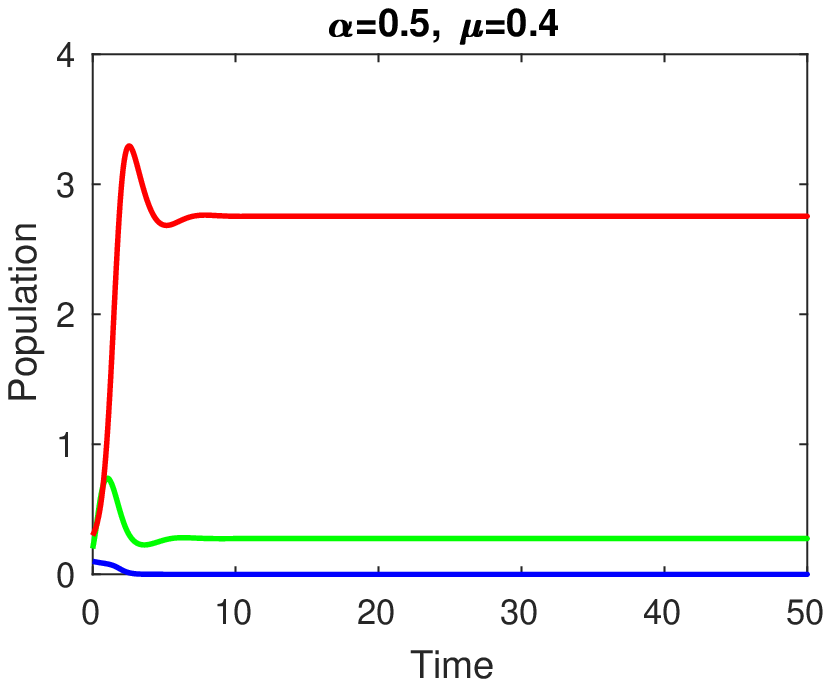}
   \includegraphics[scale=0.75]{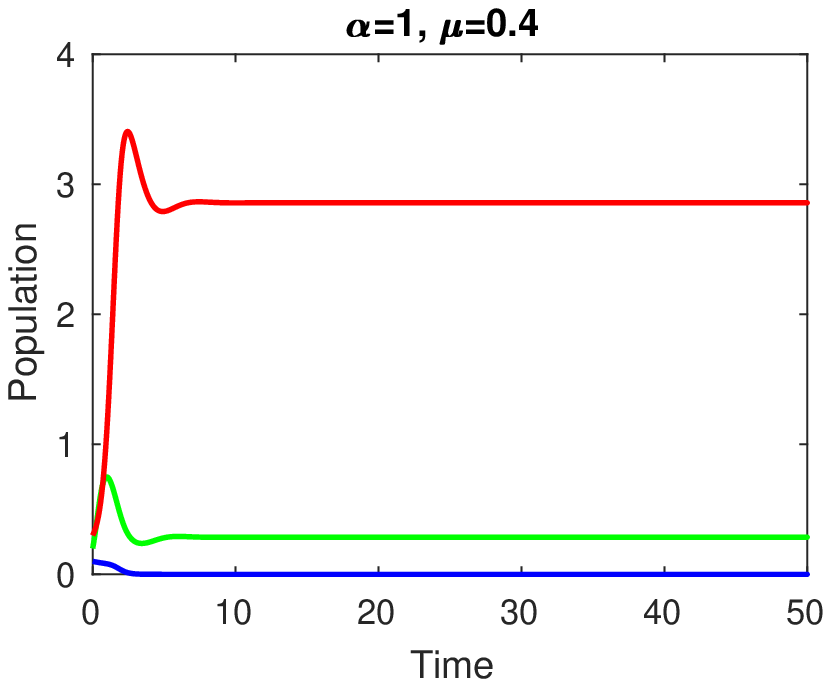}
  \includegraphics[scale=0.75]{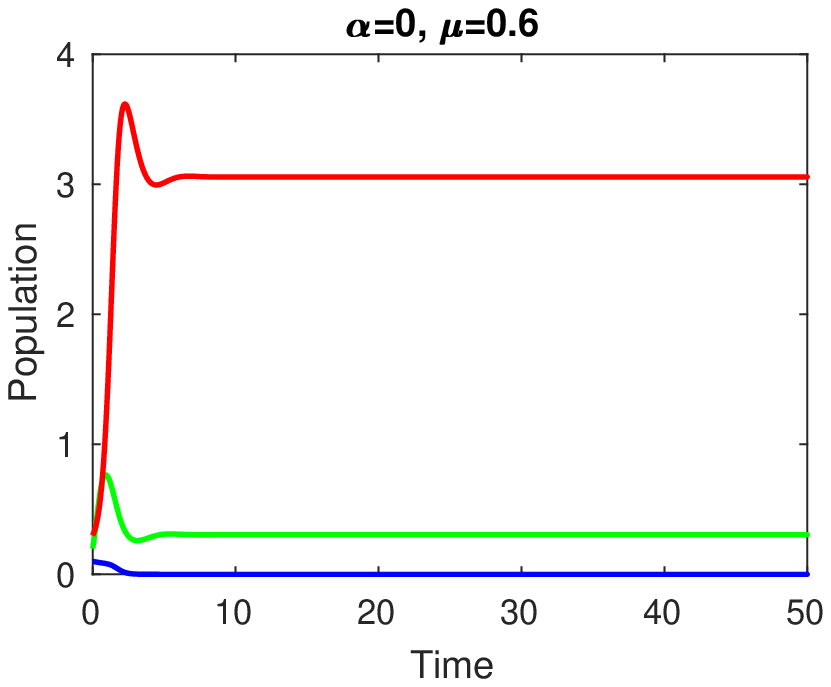}
    \includegraphics[scale=0.75]{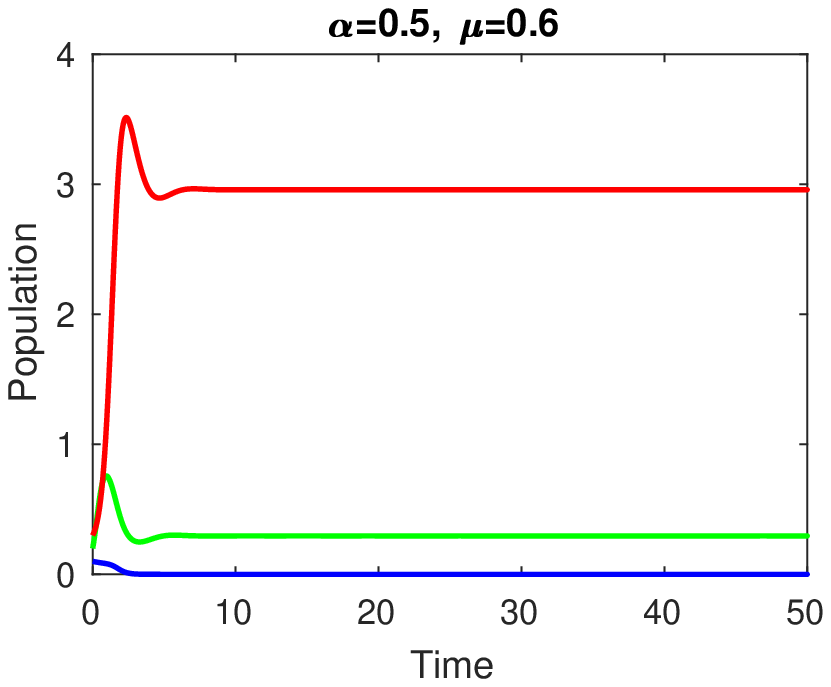}
      \end{figure}
    \begin{figure}
  \centering
  \includegraphics[scale=0.75]{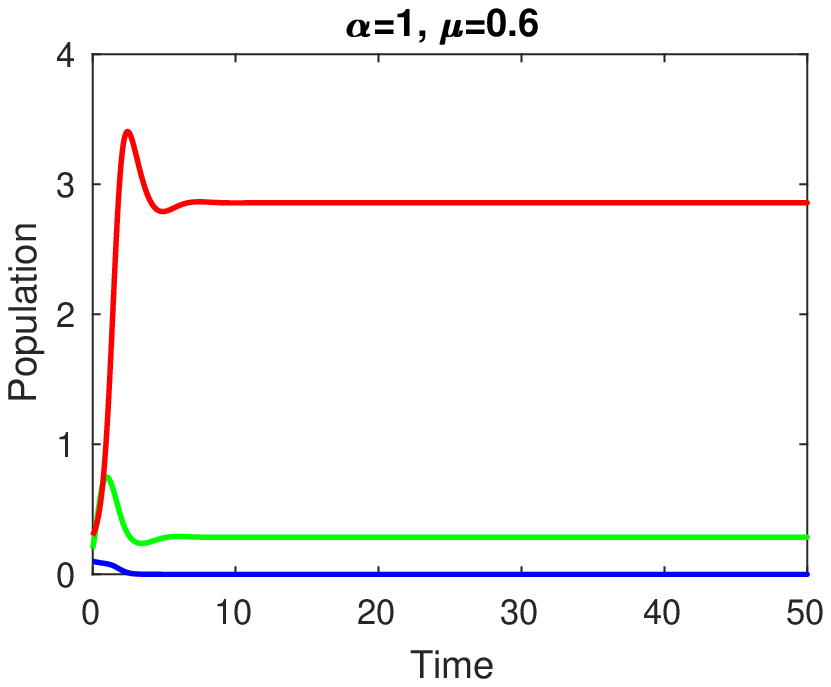} 
     \includegraphics[scale=0.75]{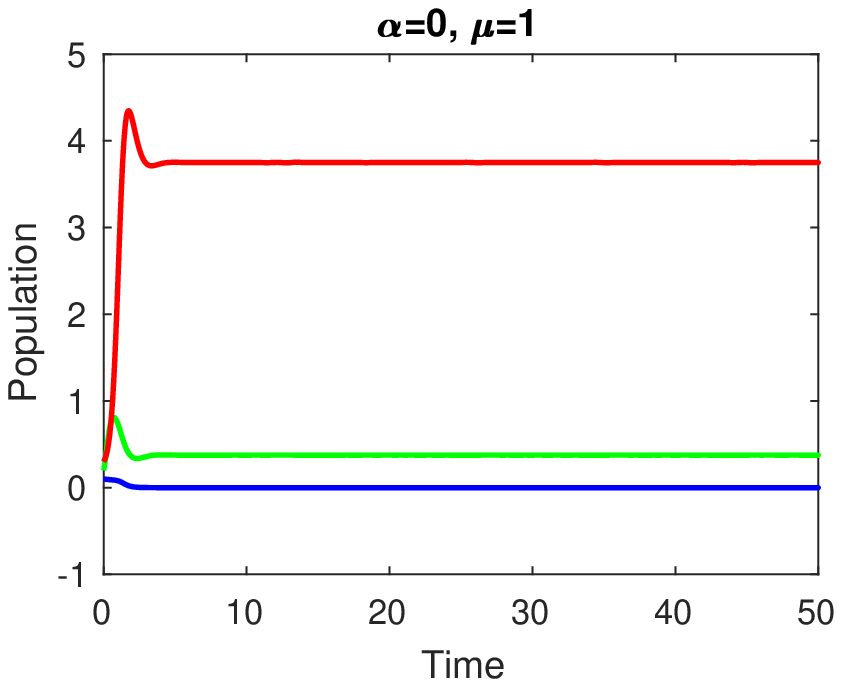}
     \includegraphics[scale=0.75]{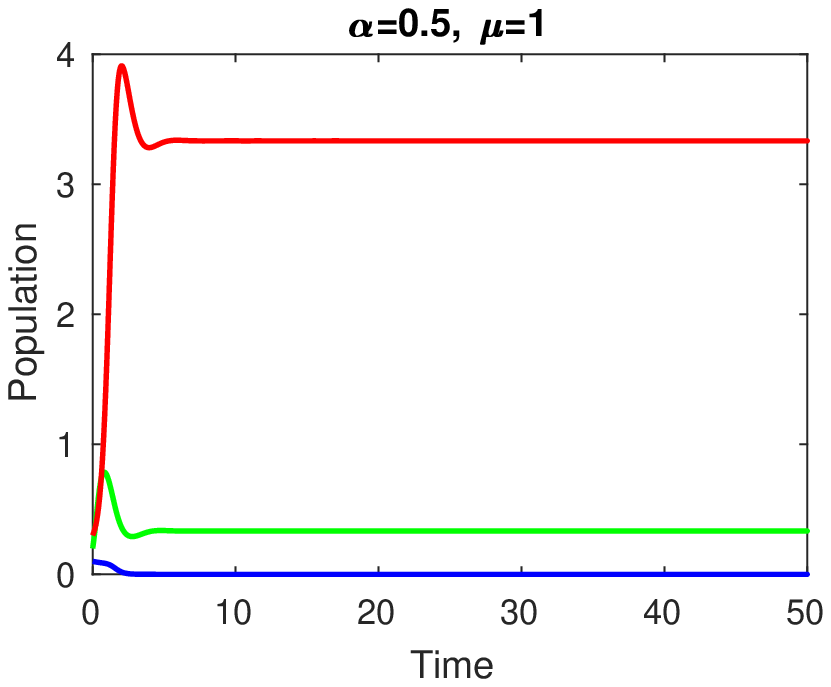}
 \includegraphics[scale=0.75]{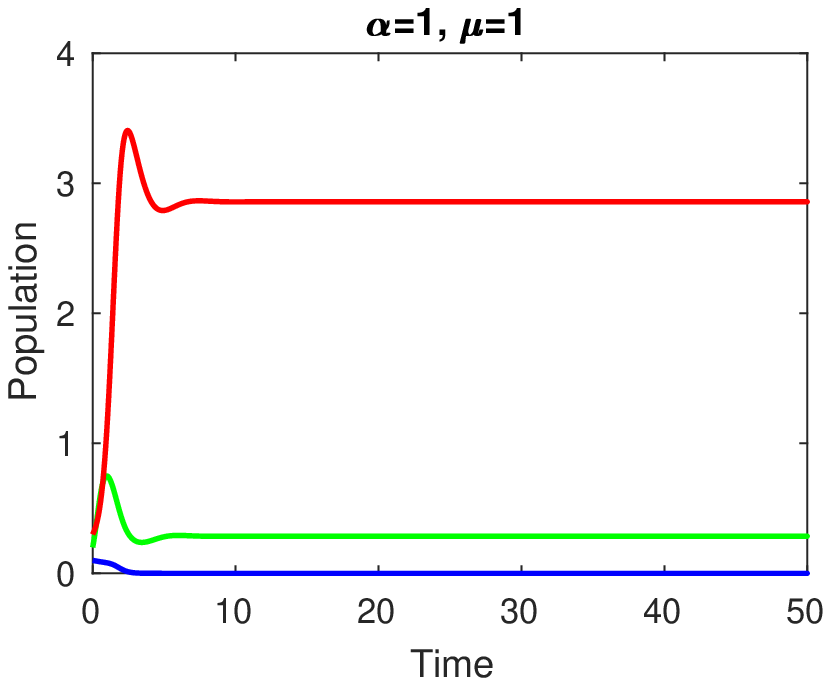}        \caption{Variation of preys and predator populations against the time for the system (\ref{gr2a}) for different values of $\alpha,\mu$. Blue, green and red curves show the horizontal membership function of the population densities $\mathbb{K}(\widehat{p}(u)), \mathbb{K}(\widehat{q}(u))$ and $\mathbb{K}(\widehat{r}(u))$, respectively.}
  \label{grfig:F}
\end{figure}
\item[(vi)] The equilibrium point $E_5(\mathbb{K}(\widehat{p_{e_5}}),0,\mathbb{K}(\widehat{r_{e_5}}))$ is locally asymptotically stable if 
\[\mathbb{K}(\widehat{a_2})<\frac{\mathbb{K}(\widehat{a_1})\mathbb{K}(\widehat{a_4})^2}{(\mathbb{K}(\widehat{a_1})\mathbb{K}(\widehat{a_3})+\mathbb{K}(\widehat{a_4}))^2}.\]
Therefore the fuzzy equilibrium point $\widehat{E}_5(\widehat{p_{e_5}},{0},\widehat{r_{e_5}})$ is locally asymptotically fuzzy stable if 
\[\widehat{a_2}<\frac{\widehat{a_1}\widehat{a_4}^2}{(\widehat{a_1}\widehat{a_3}+\widehat{a_4})^2}.\]
\begin{exa}\label{grex2}
\end{exa}
 Let $\widehat{a_1}=(2,4,6),\widehat{a_2}=(0.01,0.02,0.03),\widehat{a_3}=(0.1,0.2,0.3),\widehat{a_4}=(1,2,3),\widehat{a_5}=(3,4,5),\widehat{p}_0=0.1,\widehat{q}_0=0.2,\widehat{r}_0=0.3$, where $\widehat{a_1},\widehat{a_2},\widehat{a_3},\widehat{a_4},\widehat{a_5}$ are triangular fuzzy numbers. Now, the horizontal membership functions of the given triangular fuzzy numbers are given by
\begin{eqnarray*}
\mathbb{K}(\widehat{a_1})&=&a_1^{gr}(\alpha,\mu_{a_1})=2+2\alpha+\mu_{a_1}(4-4\alpha),\\
\mathbb{K}(\widehat{a_2})&=&a_2^{gr}(\alpha,\mu_{a_2})=0.01+0.01\alpha+\mu_{a_2}(0.02-0.02\alpha),\\
\mathbb{K}(\widehat{a_3})&=&a_3^{gr}(\alpha,\mu_{a_3})=0.1+0.1\alpha+\mu_{a_3}(0.2-0.2\alpha),\\
\mathbb{K}(\widehat{a_4})&=&a_4^{gr}(\alpha,\mu_{a_4})=1+\alpha+\mu_{a_4}(2-2\alpha),\\
\mathbb{K}(\widehat{a_5})&=&a_5^{gr}(\alpha,\mu_{a_5})=3+\alpha+\mu_{a_5}(2-2\alpha),\\
\mathbb{K}(\widehat{p}_0)&=&0.1,\mathbb{K}(\widehat{q}_0)=0.2,\mathbb{K}(\widehat{r}_0)=0.3.
\end{eqnarray*}
Further, we assume $\mu_p,\mu_q,\mu_r,\mu_{a_1},\mu_{a_2},\mu_{a_3},\mu_{a_4},\mu_{a_5}=\mu\in\{0,0.4,0.6,1\}$ and $\alpha\in\{0,0.5,1\}$.
\begin{table}[ht!]
 \begin{center}
   \begin{adjustbox}{max width=0.9\linewidth}
  \begin{tabular}{|p{.4cm}|p{3.1cm}| p{3.1cm}|p{3.1cm}|} 
  \hline
   $\mu$ & $\alpha=0$&  $\alpha=0.5$& $\alpha=1$\\
   \hline
   0&(0.1667,0,1.6667)&(0.2308,0,2.3077)&(0.2857,0,2.8571) \\
   \hline
    0.4&(0.2647,0,2.6471)&(0.2754,0,2.7536)&(0.2857,0,2.8571) \\
   \hline
    0.6&(0.3056,0,3.0556)&(0.2958,0,2.9577)&(0.2857,0,2.8571)\\
   \hline
    1&(0.375,0,3.75)&(0.3333,0,3.3333)&(0.28571,0,2.8571)\\
    \hline
  \end{tabular}
  \end{adjustbox}
   \caption{Equilibrium point $E_5(\mathbb{K}(\widehat{p_{e_5}}),0,\mathbb{K}(\widehat{r_{e_5}}))$ for different values of $\alpha,\mu$}
  \label{grtab:table2}
 \end{center}
\end{table}
For the given data set, we observe that the stability condition of $E_5$ is well satisfied and the system (\ref{gr2a}) have different equilibrium points corresponding to different values of $\alpha,\mu$, as shown in Table \ref{grtab:table2}. From Table \ref{grtab:table2}, we realize that for a fixed value of $\alpha$ the equilibrium point increases with increasing value of $\mu$ and at $\alpha=1$ remain same because the left and right interval of triangular fuzzy number coincides with each other. For $\mu=0,0.4$ the equilibrium point increases while for $\mu=0.6,1$ the equilibrium point decreases with increasing value of $\alpha$. Figure \ref{grfig:2} shows that for the given data in Example \ref{grex2}, the population density $\mathbb{K}(\widehat{q}(u))$ of prey eventually extinct while the population densities $\mathbb{K}(\widehat{p}(u))$ and $\mathbb{K}(\widehat{r}(u))$ of prey and predator persist, respectively and eventually get their steady states given in Table \ref{grtab:table2} corresponding to different values of $\alpha,\mu$ and become asymptotically stable. At $\alpha=1$, it shows crisp behavior.
\begin{figure}
  \centering
  \includegraphics[scale=0.75]{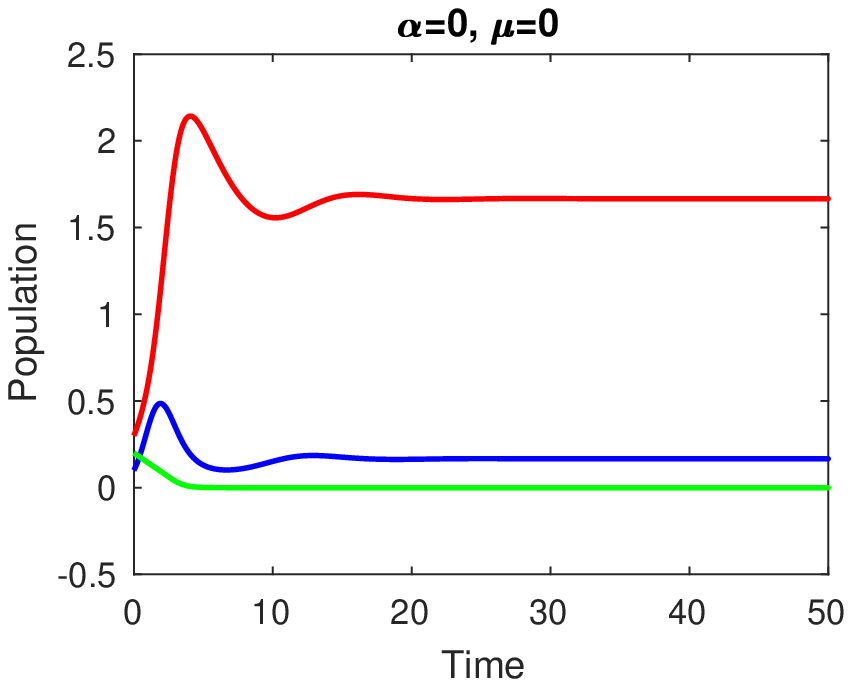}  
  \includegraphics[scale=0.75]{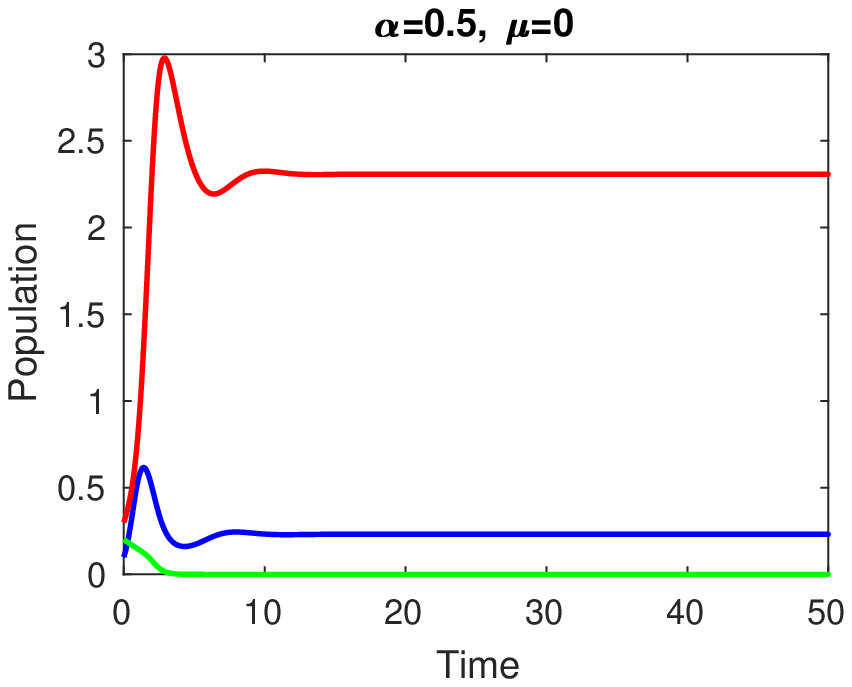}
\includegraphics[scale= 0.75]{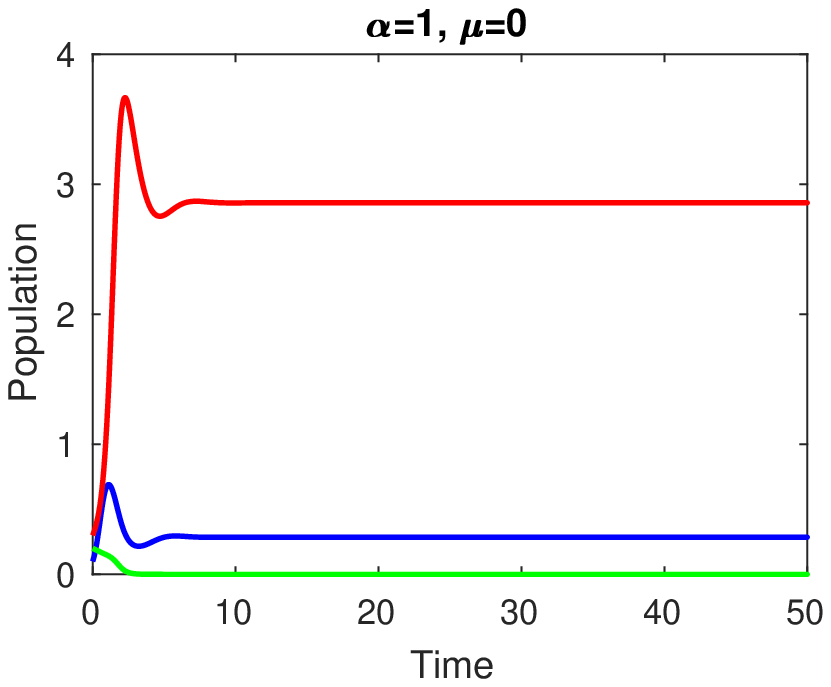}
  \includegraphics[scale=0.75]{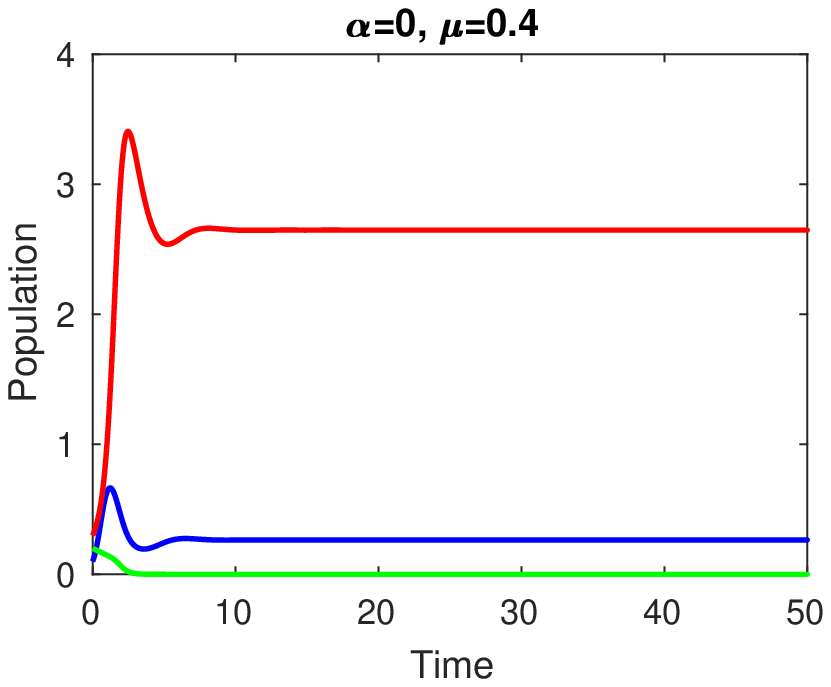}  
  \includegraphics[scale=0.75]{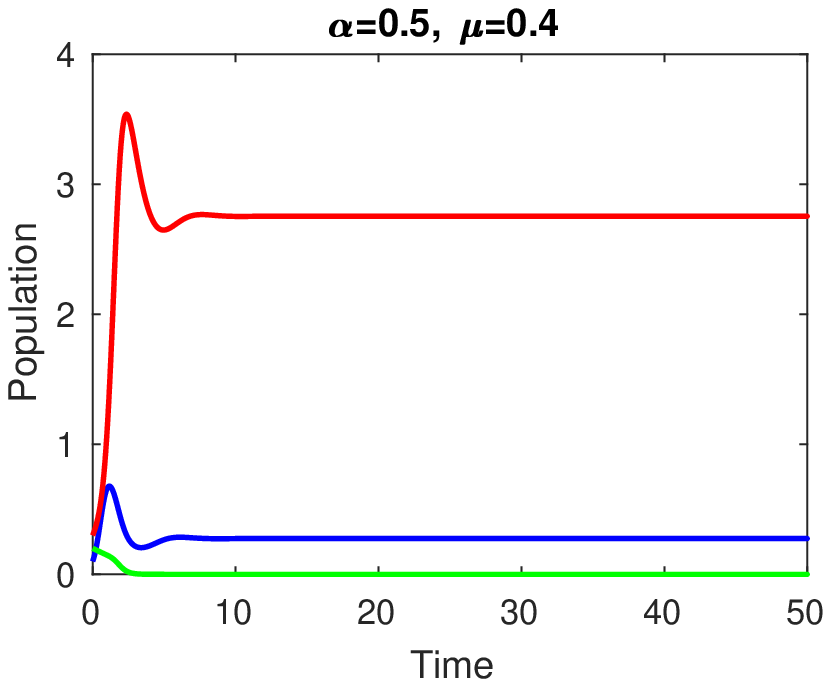}
   \includegraphics[scale=0.75]{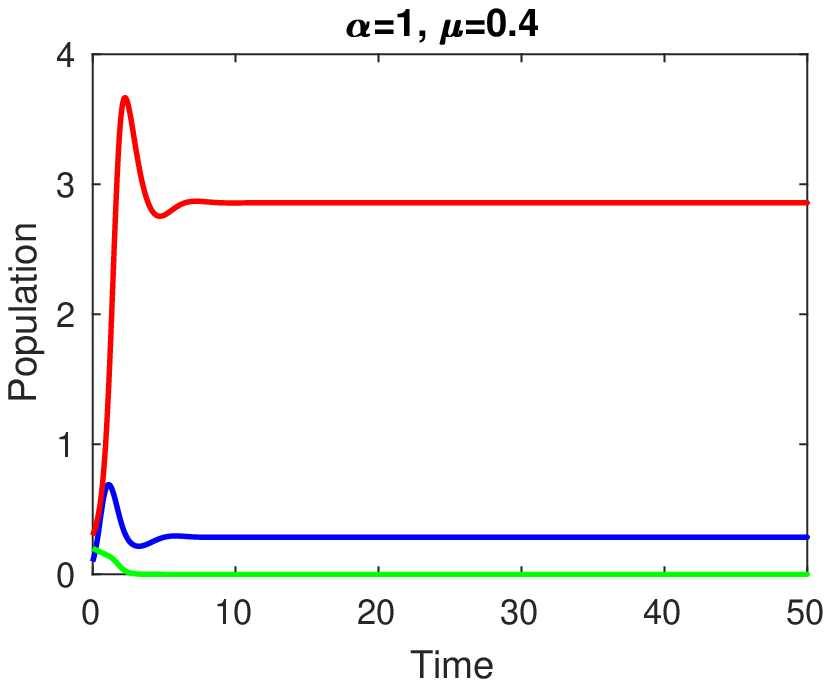}
    \includegraphics[scale=0.75]{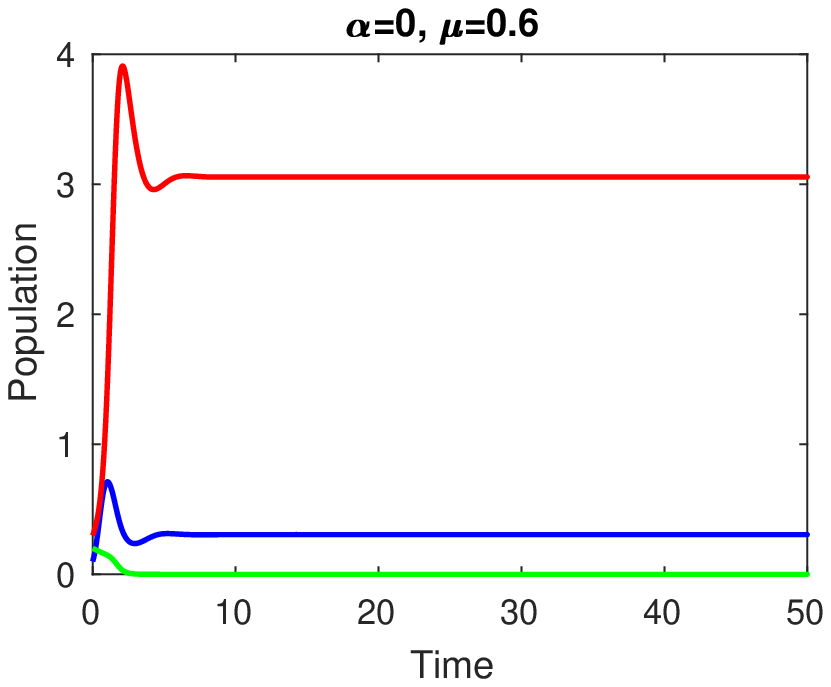}
    \includegraphics[scale=0.75]{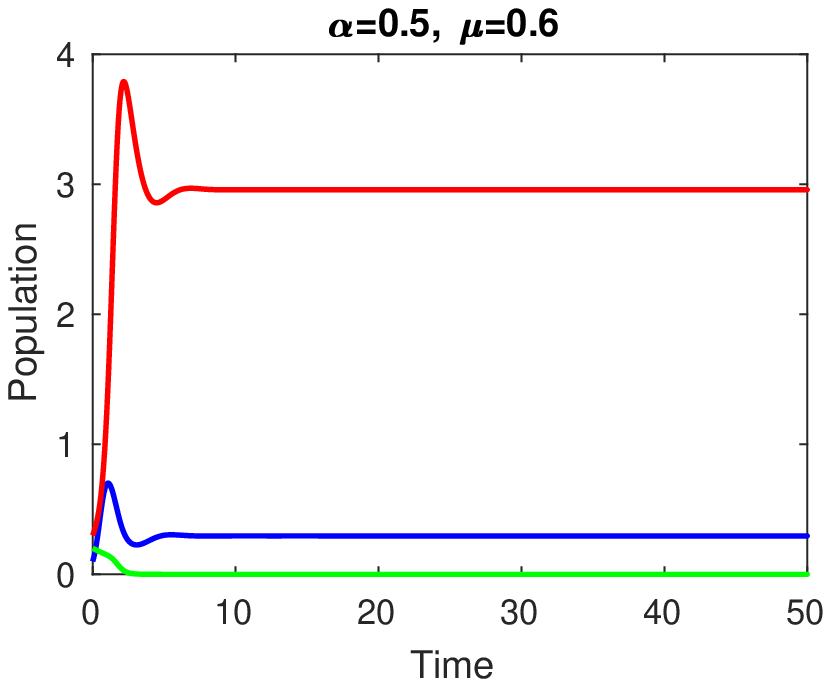}
     \end{figure}
     \begin{figure}
  \centering
    \includegraphics[scale=0.75]{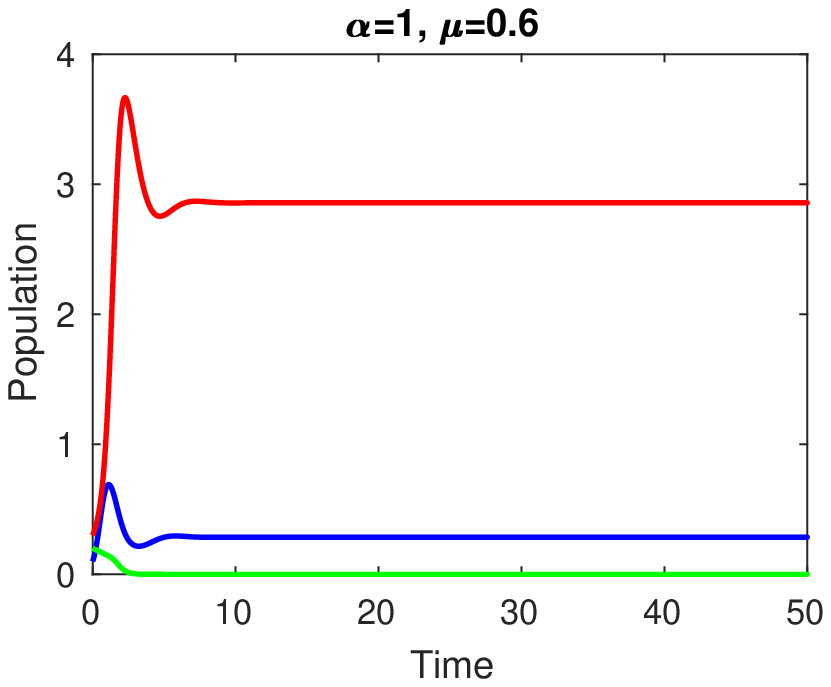}
     \includegraphics[scale=0.75]{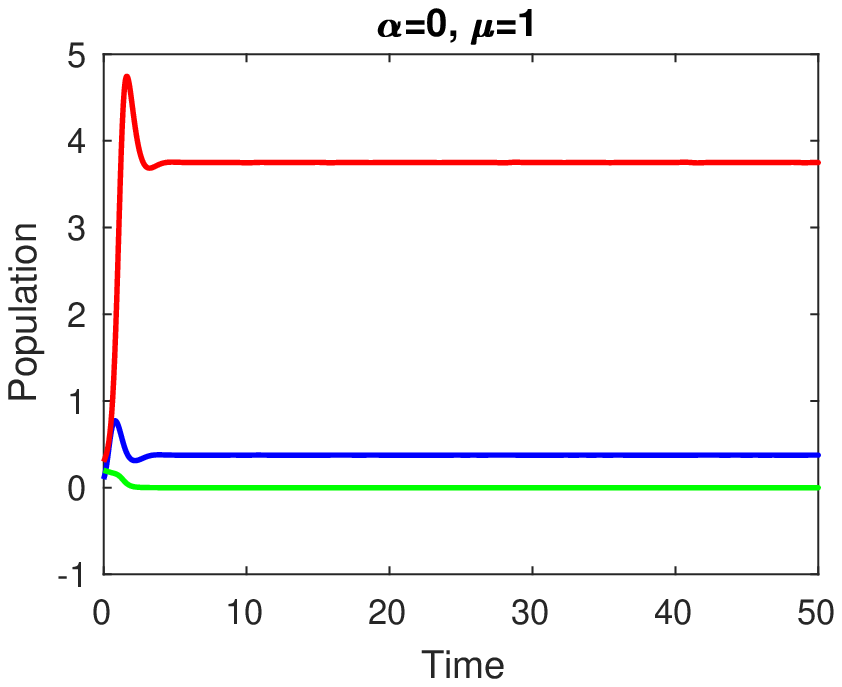}
      \includegraphics[scale=0.75]{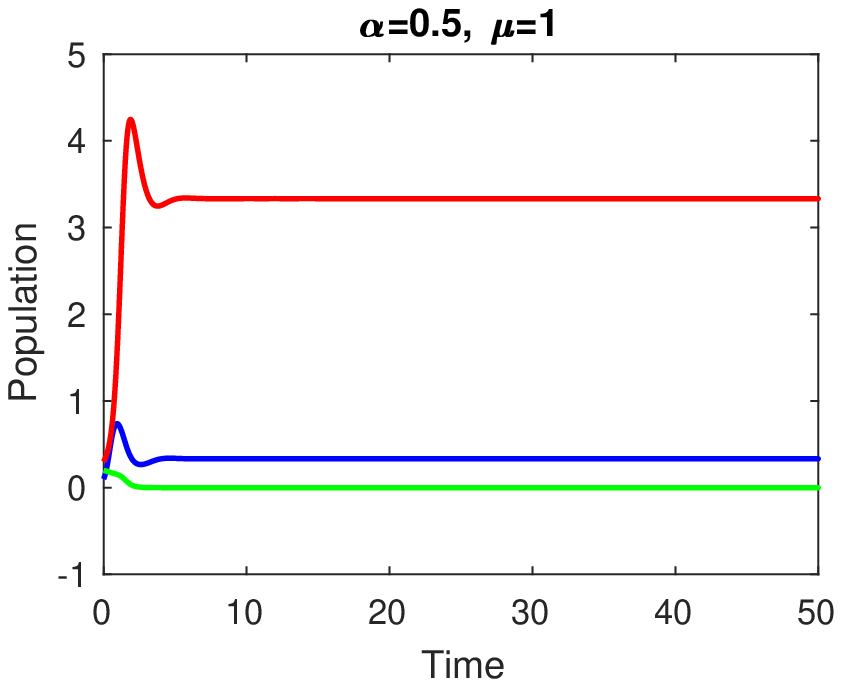}
     \includegraphics[scale=.75]{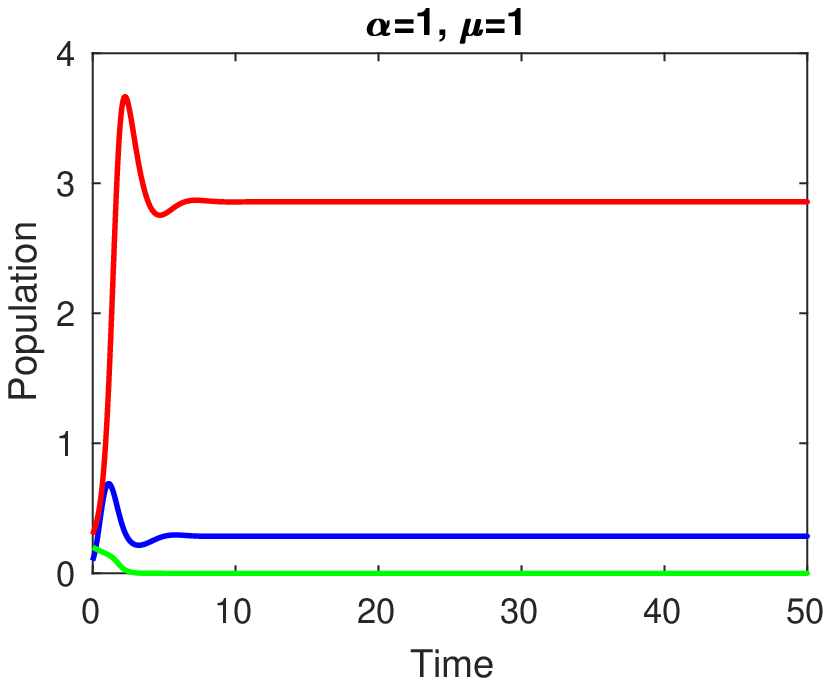}
     \caption{Variation of preys and predator populations against the time for the system (\ref{gr2a}) for different values of $\alpha,\mu$. Blue, green and red curves show the horizontal membership function of the population densities $\mathbb{K}(\widehat{p}(u)), \mathbb{K}(\widehat{q}(u))$ and $\mathbb{K}(\widehat{r}(u))$, respectively.}
  \label{grfig:2}
\end{figure}
\item[(vii)] The equilibrium point $E_6(1,1,\mathbb{K}(\widehat{r_{e_6}}))$ is locally asymptotically stable if 
\[\mathbb{K}(\widehat{a_1})\mathbb{K}(\widehat{a_2})>\frac{(\mathbb{K}(\widehat{a_4})+\mathbb{K}(\widehat{a_5}))^2}{\mathbb{K}(\widehat{a_3})^2}.\]
Therefore the fuzzy equilibrium point $\widehat{E}_6({1},{1},\widehat{r_{e_6}})$ is locally asymptotically fuzzy stable if 
\[\widehat{a_1}\widehat{a_2}>\frac{(\widehat{a_4}+\widehat{a_5})^2}{\widehat{a_3}^2}.\]
\begin{exa}\label{grex3}
\end{exa} Let $\widehat{a_1}=(3,4,5),\widehat{a_2}=(2,4,6),\widehat{a_3}=(3,4,5),\widehat{a_4}=(1,2,3),\widehat{a_5}=(0.01,0.02,0.03),\widehat{p}_0=0.1,\widehat{q}_0=0.2,\widehat{r}_0=0.3$, where $\widehat{a_1},\widehat{a_2},\widehat{a_3},\widehat{a_4},\widehat{a_5}$ are triangular fuzzy numbers. Now, the horizontal membership functions of the given triangular fuzzy numbers are given by
\begin{eqnarray*}
\mathbb{K}(\widehat{a_1})&=&a_1^{gr}(\alpha,\mu_{a_1})=3+\alpha+\mu_{a_1}(2-2\alpha),\\
\mathbb{K}(\widehat{a_2})&=&a_2^{gr}(\alpha,\mu_{a_2})=2+2\alpha+\mu_{a_2}(4-4\alpha),\\
\mathbb{K}(\widehat{a_3})&=&a_3^{gr}(\alpha,\mu_{a_3})=3+\alpha+\mu_{a_3}(2-2\alpha),\\
\mathbb{K}(\widehat{a_4})&=&a_4^{gr}(\alpha,\mu_{a_4})=1+\alpha+\mu_{a_4}(2-2\alpha),\\
\mathbb{K}(\widehat{a_5})&=&a_5^{gr}(\alpha,\mu_{a_5})=0.01+0.01\alpha+\mu_{a_5}(0.02-0.02\alpha),\\
\mathbb{K}(\widehat{p}_0)&=&0.1,\mathbb{K}(\widehat{q}_0)=0.2,\mathbb{K}(\widehat{r}_0)=0.3.
\end{eqnarray*}
Further, we assume $\mu_p,\mu_q,\mu_r,\mu_{a_1},\mu_{a_2},\mu_{a_3},\mu_{a_4},\mu_{a_5}=\mu\in\{0,0.4,0.6,1\}$ and $\alpha\in\{0,0.5,1\}$.
\begin{table}[ht!]
 \begin{center}
  \begin{adjustbox}{max width=0.9\linewidth}
  \begin{tabular}{|p{.4cm}|p{2.0cm}| p{2.0cm}|p{2.0cm}|} 
  \hline
   $\mu$ & $\alpha=0$&  $\alpha=0.5$& $\alpha=1$\\
   \hline
   0&(1,1,0.3367)&(1,1,0.4329)&(1,1,0.5050) \\
   \hline
    0.4&(1,1,0.4784)&(1,1,0.4921)&(1,1,0.5050) \\
   \hline
    0.6&(1,1,0.5290)&(1,1,0.5173)&(1,1,0.5050)\\
   \hline
    1&(1,1,0.6060)&(1,1,0.5611)&(1,1,0.5050)\\
    \hline
  \end{tabular}
  \end{adjustbox}
   \caption{Equilibrium point $E_6({1},{1},\mathbb{K}(\widehat{r_{e_6}}))$ for different values of $\alpha,\mu$}
  \label{grtab:table3}
 \end{center}
 \end{table}
For the given data set, we observe that the stability condition of $E_6$ is well satisfied and the system (\ref{gr2a}) have different equilibrium points corresponding to different values of $\alpha,\mu$, as shown in Table \ref{grtab:table3}. From Table \ref{grtab:table3}, we realize that for a fixed value of $\alpha$ the equilibrium point increases with increasing value of $\mu$ and at $\alpha=1$ remain same because the left and right interval of triangular fuzzy number coincides with each other. For $\mu=0,0.4$ the equilibrium point increases while for $\mu=0.6,1$ the equilibrium point decreases with increasing value of $\alpha$. Figure \ref{grfig:3} shows that for the given data in Example \ref{grex3}, initially the population density $\mathbb{K}(\widehat{r}(u))$ of predator decreases while the population densities $\mathbb{K}(\widehat{p}(u))$ and $\mathbb{K}(\widehat{q}(u))$ of preys increases, respectively and eventually get their steady states given in Table \ref{grtab:table3} corresponding to different values of $\alpha,\mu$ and become asymptotically stable. At $\alpha=1$, it shows crisp behavior.
\begin{figure}
  \centering
  \includegraphics[scale=0.75]{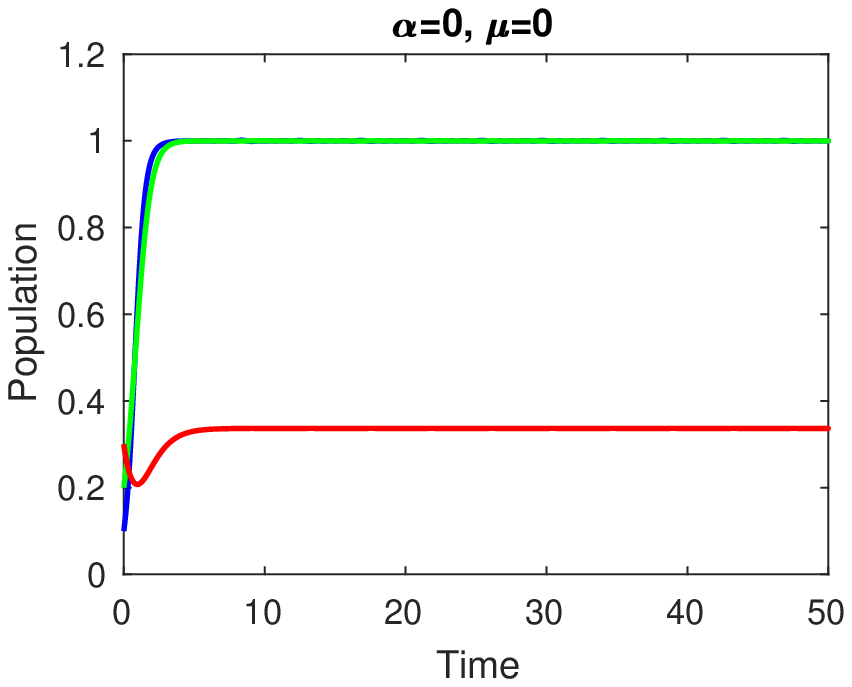}
  \includegraphics[scale=0.75]{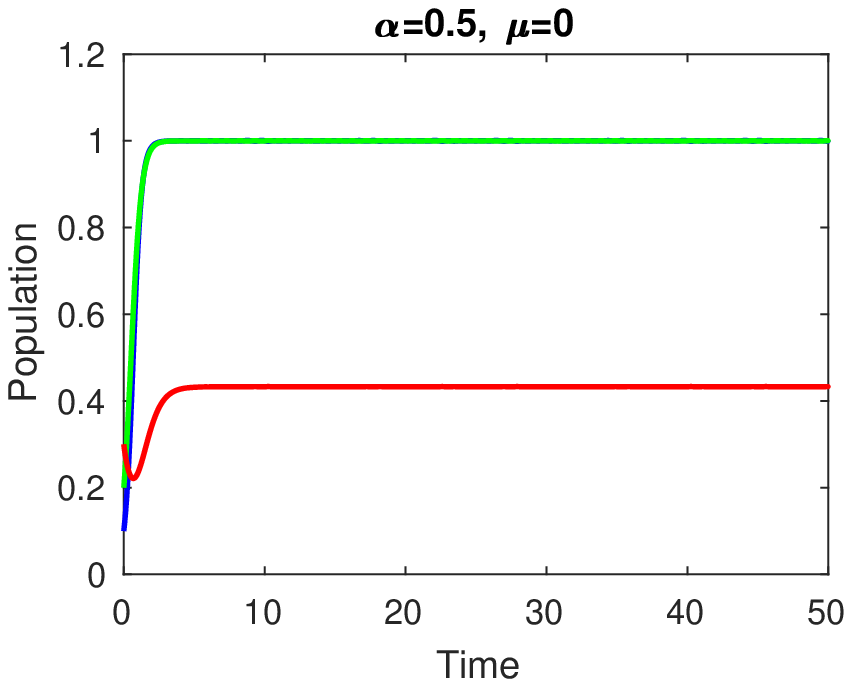}
  \includegraphics[scale= 0.75]{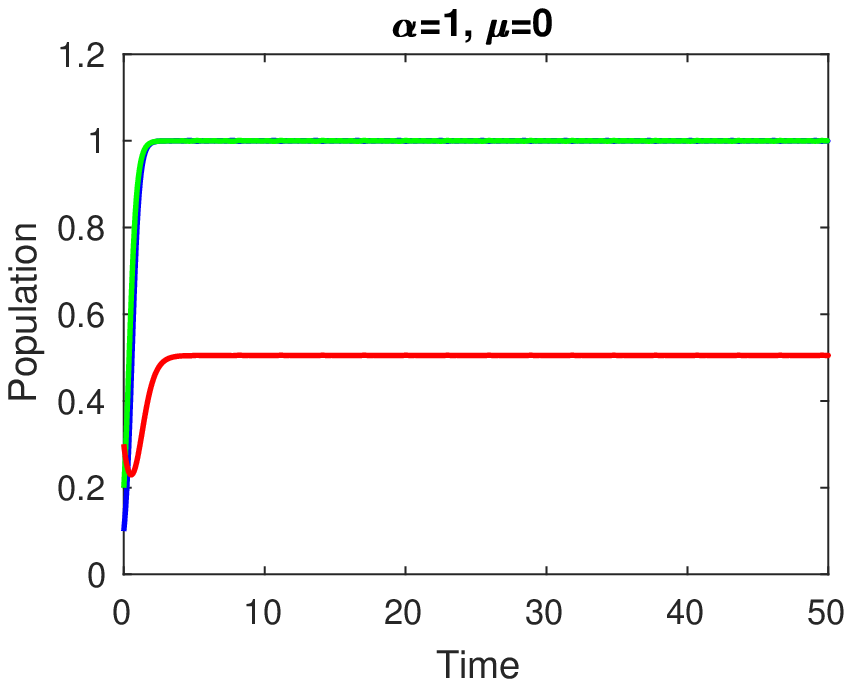}
   \includegraphics[scale=0.75]{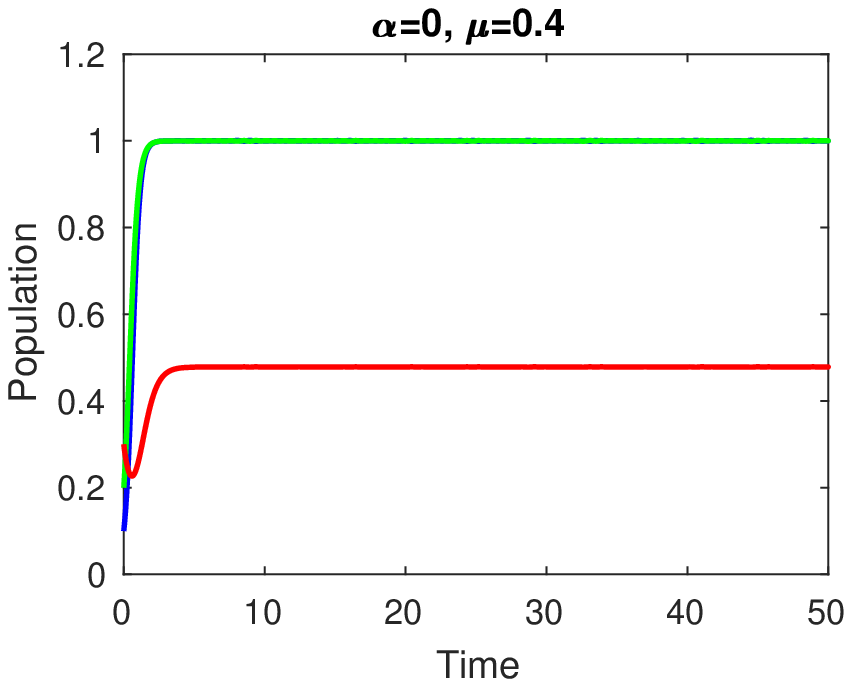}
   \includegraphics[scale=0.75]{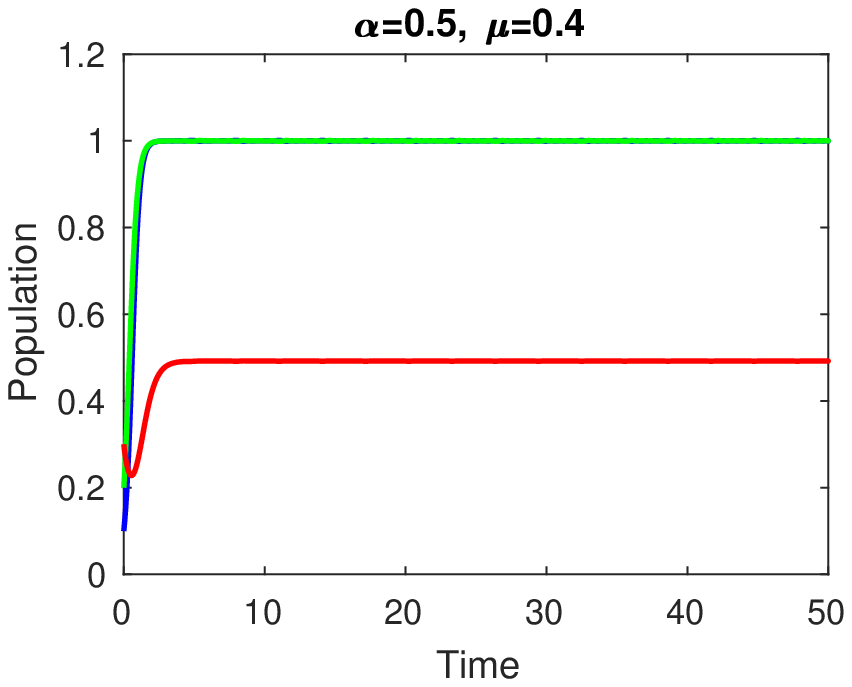}
   \includegraphics[scale= 0.65]{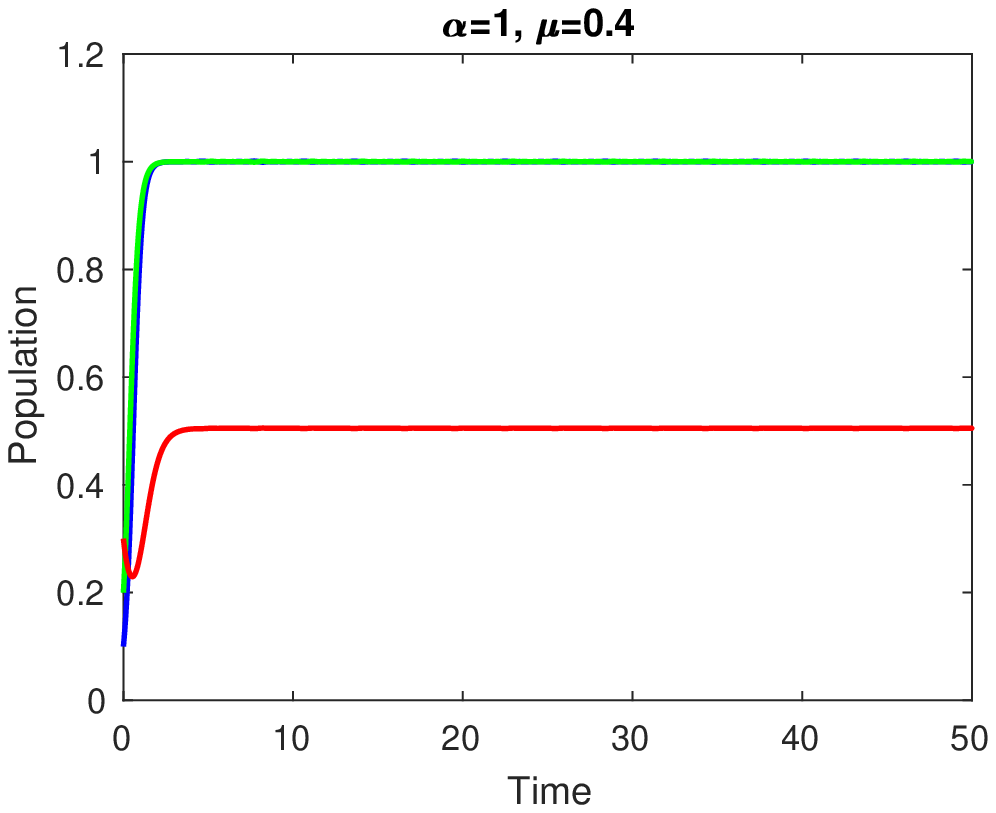}
\includegraphics[scale=0.75]{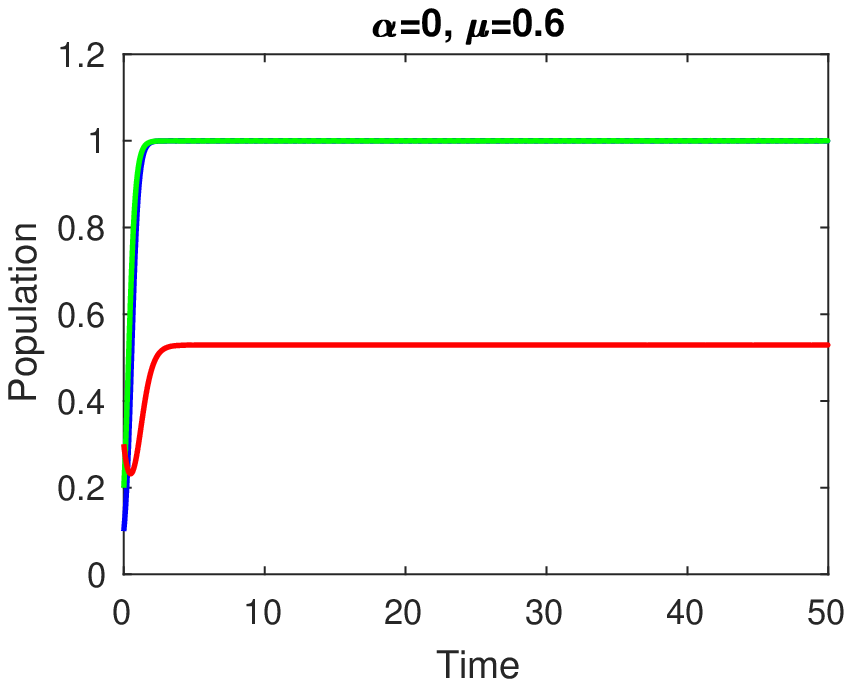}
\includegraphics[scale=0.75]{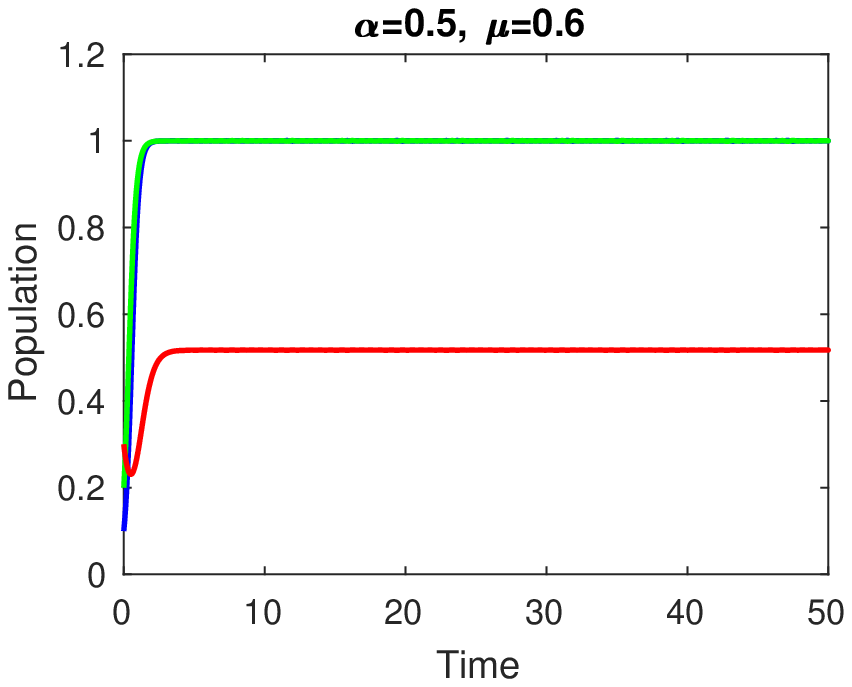}
\end{figure}
\begin{figure}
  \centering
\includegraphics[scale= 0.65]{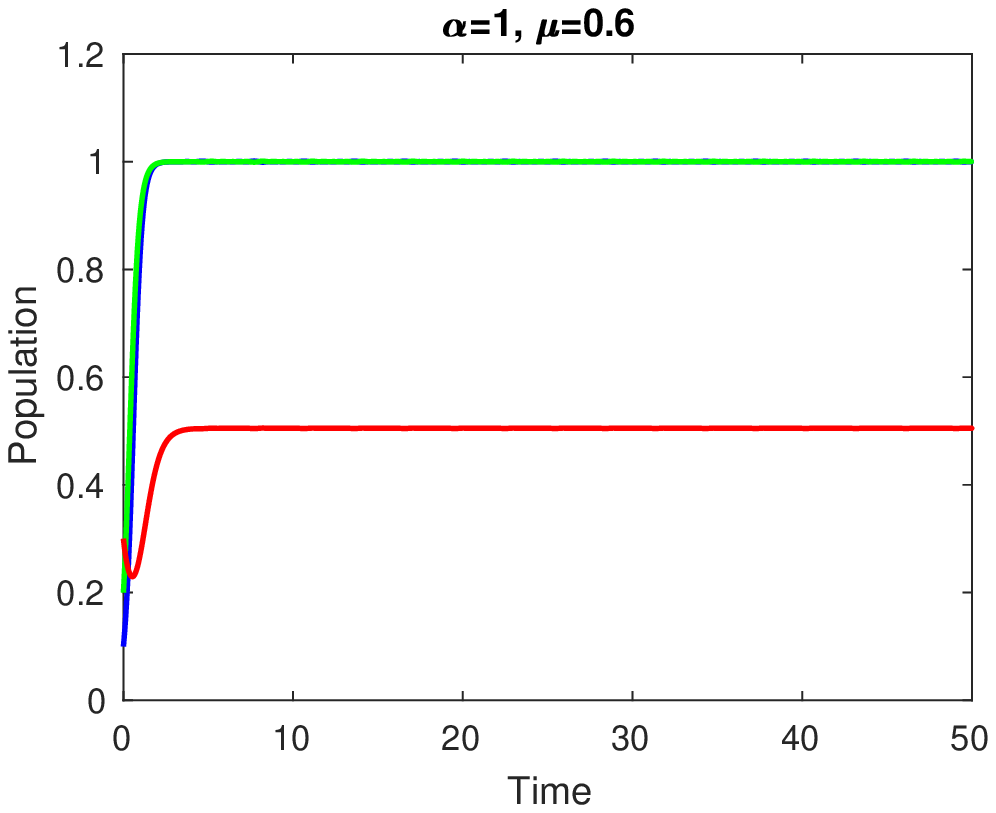}
\includegraphics[scale=0.8]{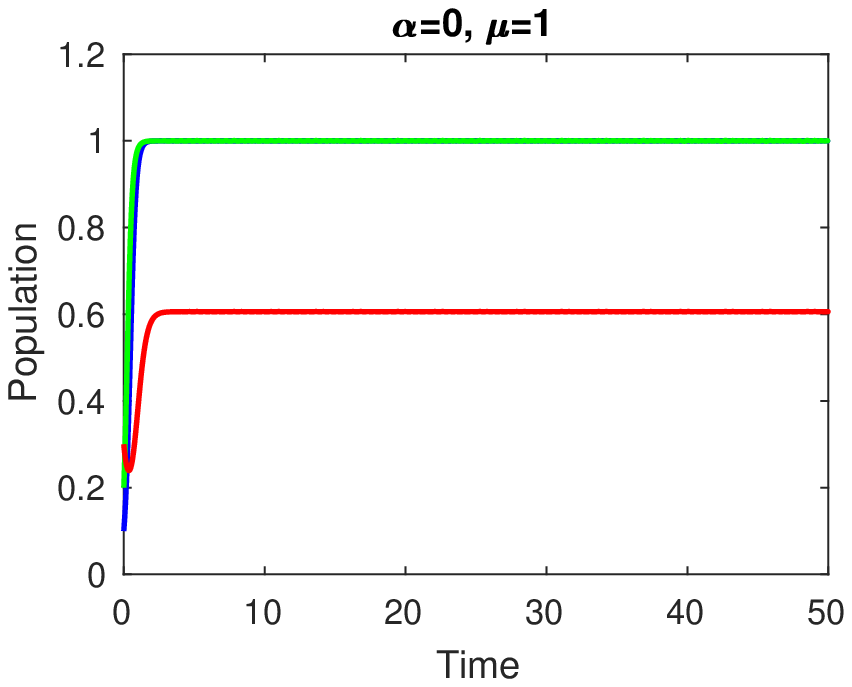}
\includegraphics[scale=0.8]{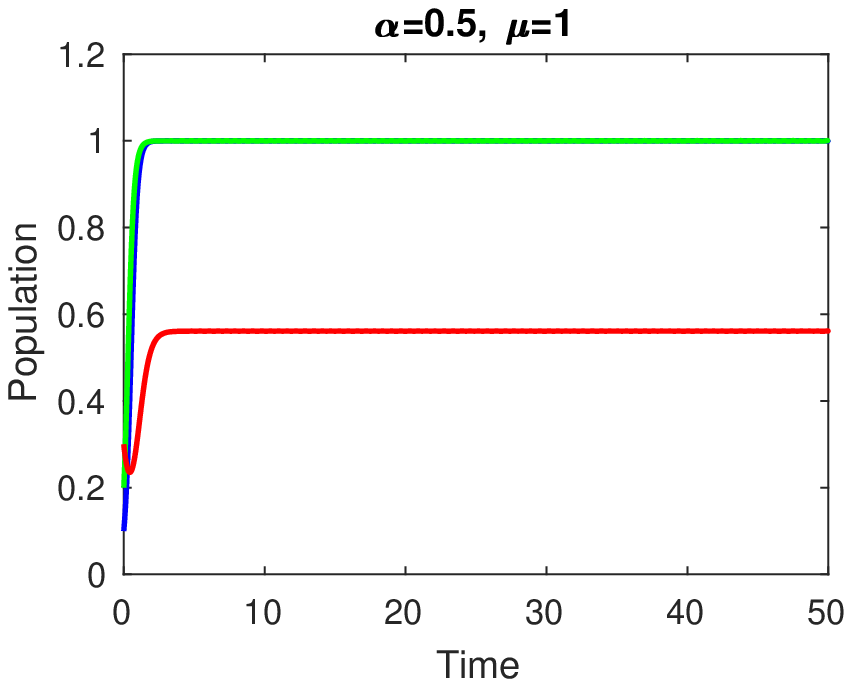}
\includegraphics[scale= 0.65]{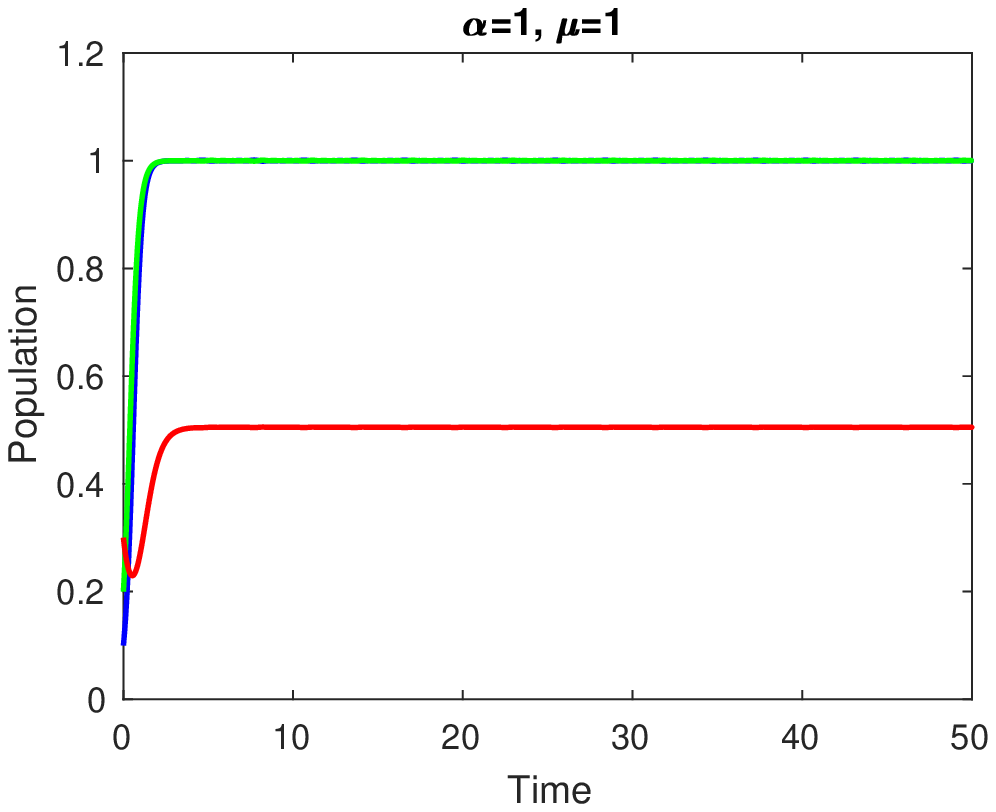} \caption{Variation of preys and predator populations against the time for the system (\ref{gr2a}) for different values of $\alpha,\mu$. Blue, green and red curves show the horizontal membership function of the population densities $\mathbb{K}(\widehat{p}(u)), \mathbb{K}(\widehat{q}(u))$ and $\mathbb{K}(\widehat{r}(u))$, respectively.}
 \label{grfig:3}
\end{figure}
\item[(viii)] The equilibrium point $E_7(\mathbb{K}(\widehat{p_{e_7}}),\mathbb{K}(\widehat{q_{e_7}}),\mathbb{K}(\widehat{r_{e_7}}))$ is locally asymptotically stable if 
\begin{eqnarray*}
\mathbb{K}(\widehat{a_4})+\mathbb{K}(\widehat{a_5})&>&\mathbb{K}(\widehat{a_3})\sqrt{\mathbb{K}(\widehat{a_1})\mathbb{K}(\widehat{a_2})},\\
\mathbb{K}(\widehat{a_2})&\geq& \mathbb{K}(\widehat{a_1}),\\
\mathbb{K}(\widehat{a_1})\mathbb{K}(\widehat{a_3})&\geq& 1.
\end{eqnarray*}
Therefore the fuzzy equilibrium point $\widehat{E}_7(\widehat{p_{e_7}},\widehat{q_{e_7}},\widehat{r_{e_7}})$ is locally asymptotically fuzzy stable if 
\begin{eqnarray*}
\widehat{a_4}+\widehat{a_5}&>&\widehat{a_3}\sqrt{\widehat{a_1}\widehat{a_2}},\\
\widehat{a_2}&\geq&\widehat{a_1},\\
\widehat{a_1}\widehat{a_3}&\geq& 1.
\end{eqnarray*}
\begin{exa}\label{grex4}
\end{exa} Let $\widehat{a_1}=(1,2,3),\widehat{a_2}=(2,4,6),\widehat{a_3}=(1,2,3),\widehat{a_4}=(3,4,5),\widehat{a_5}=(1,2,3),\widehat{p}_0=0.1,\widehat{q}_0=0.2,\widehat{r}_0=0.3$, where $\widehat{a_1},\widehat{a_2},\widehat{a_3},\widehat{a_4},\widehat{a_5}$ are triangular fuzzy numbers. Now, the horizontal membership functions of the given triangular fuzzy numbers are given by
\begin{eqnarray*}
\mathbb{K}(\widehat{a_1})&=&a_1^{gr}(\alpha,\mu_{a_1})=1+\alpha+\mu_{a_1}(2-2\alpha),\\
\mathbb{K}(\widehat{a_2})&=&a_2^{gr}(\alpha,\mu_{a_2})=2+2\alpha+\mu_{a_2}(4-4\alpha),\\
\mathbb{K}(\widehat{a_3})&=&a_3^{gr}(\alpha,\mu_{a_3})=1+\alpha+\mu_{a_3}(2-2\alpha),\\
\mathbb{K}(\widehat{a_4})&=&a_4^{gr}(\alpha,\mu_{a_4})=3+\alpha+\mu_{a_4}(2-2\alpha),\\
\mathbb{K}(\widehat{a_5})&=&a_5^{gr}(\alpha,\mu_{a_5})=1+\alpha+\mu_{a_5}(2-2\alpha),\\
\mathbb{K}(\widehat{p}_0)&=&0.1,\mathbb{K}(\widehat{q}_0)=0.2,\mathbb{K}(\widehat{r}_0)=0.3.
\end{eqnarray*}
Further, we assume $\mu_p,\mu_q,\mu_r,\mu_{a_1},\mu_{a_2},\mu_{a_3},\mu_{a_4},\mu_{a_5}=\mu\in\{0,0.4,0.6,1\}$ and $\alpha\in\{0,0.5,1\}$.
\begin{table}[ht!]
 \begin{center}
  \begin{adjustbox}{max width=0.9\linewidth}
  \begin{tabular}{|p{0.5cm}|p{4cm}| p{4cm}|p{4cm}|} 
  \hline
   $\mu$ & $\alpha=0$&  $\alpha=0.5$& $\alpha=1$\\
   \hline
   0&(0.3025,0.5068,1.414)&(0.6014,0.7181,2.1213)&(0.9366,0.9552,2.8284) \\
   \hline
    0.4&(0.7993,0.8581,2.5456)&(0.8675,0.9063,2.6870)&
    (0.9366,0.9552,2.8284) \\
   \hline
    0.6&(1.0773,1.0546,3.1113)&(1.0066,1.0046,2.9698)&
    (0.9366,0.9552,2.8284)\\
   \hline
   1&(1.6639,1.4695,4.2426)&
   (1.2934,1.2075,3.5355)&
    (0.9366,0.9552,2.8284)\\
    \hline
  \end{tabular}
  \end{adjustbox}
 \caption{Equilibrium point $E_7(\mathbb{K}(\widehat{p_{e_7}}),\mathbb{K}(\widehat{q_{e_7}}),\mathbb{K}(\widehat{r_{e_7}}))$ for different values of $\alpha,\mu$}
  \label{grtab:table4}
 \end{center}
 \end{table}
For the given data set, we observe that the stability condition of $E_7(\widehat{p_{e_7}},\widehat{q_{e_7}},\widehat{r_{e_7}})$ is well satisfied and the system (\ref{gr2a}) have different equilibrium points corresponding to different values of $\alpha,\mu$, as shown in Table \ref{grtab:table4}. From Table \ref{grtab:table4}, we realize that for a fixed value of $\alpha$ the equilibrium point increases with increasing value of $\mu$ and at $\alpha=1$ remain same because the left and right interval of triangular fuzzy number coincides with each other. For $\mu=0,0.4$ the equilibrium point increases while for $\mu=0.6,1$ the equilibrium point decreases with increasing value of $\alpha$. Figure \ref{grfig:4} shows that for the given data in Example \ref{grex4}, initially the population densities $\mathbb{K}(\widehat{p}(u)),\mathbb{K}(\widehat{q}(u))$ and $\mathbb{K}(\widehat{r}(u))$ of preys and predator, respectively increases and eventually get their steady states given in Table \ref{grtab:table4} corresponding to different values of $\alpha,\mu$ and become asymptotically stable. At $\alpha=1$, it shows crisp behavior.
\begin{figure}
  \centering
  \includegraphics[scale=0.75]{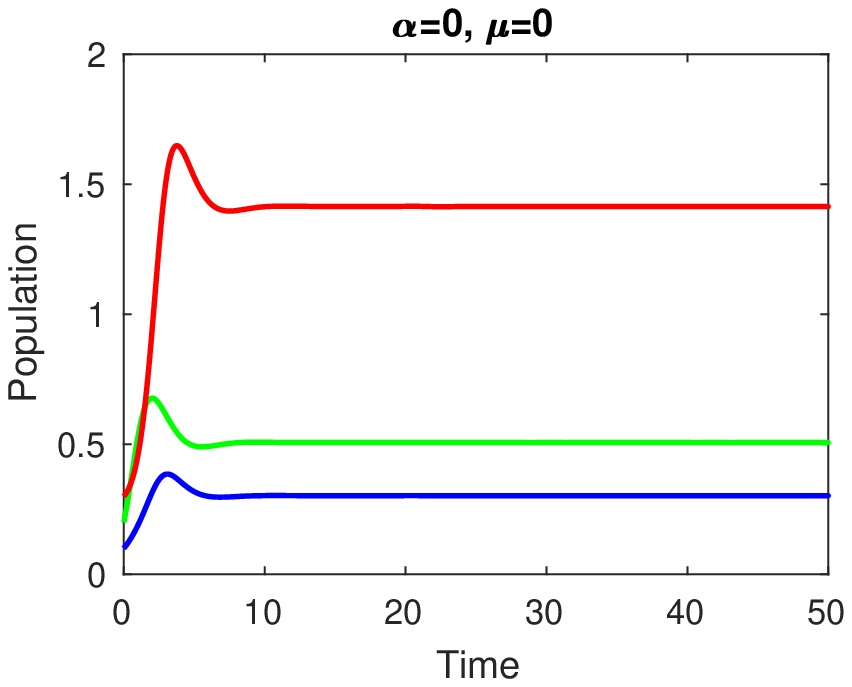}
   \includegraphics[scale=0.75]{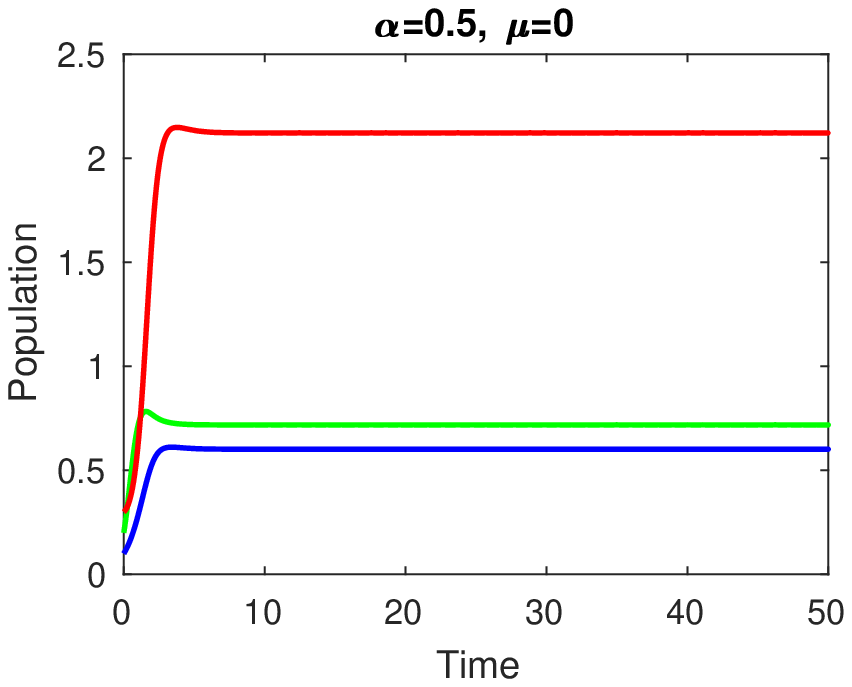}
     \includegraphics[scale= 0.75]{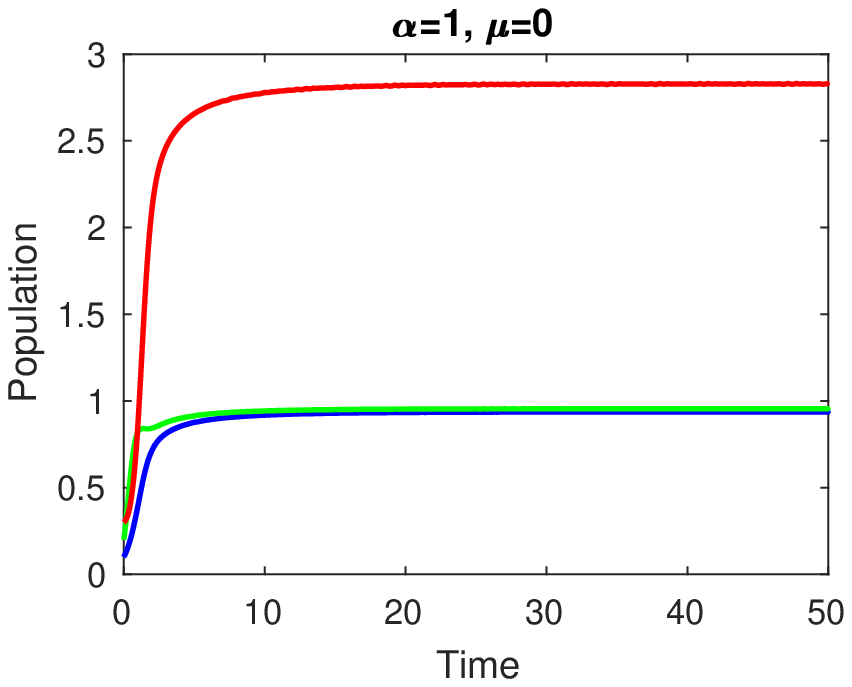}
   \includegraphics[scale=0.75]{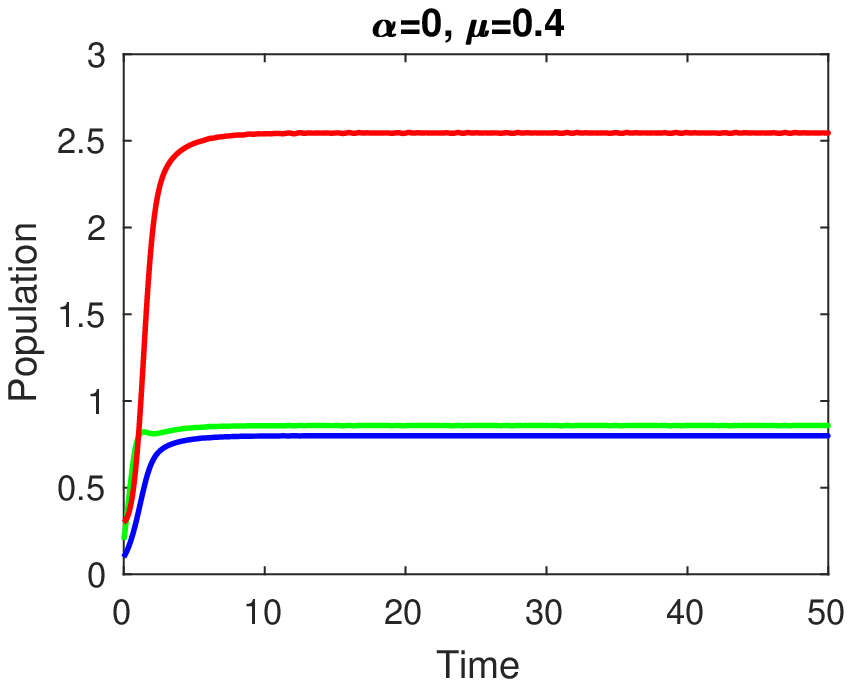}
    \includegraphics[scale=0.75]{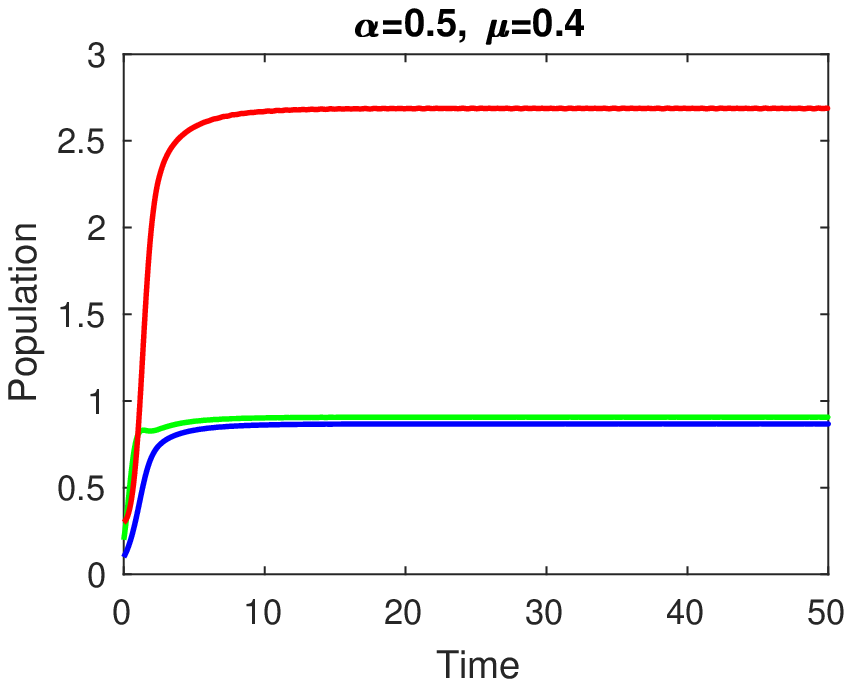}
  \includegraphics[scale= 0.65]{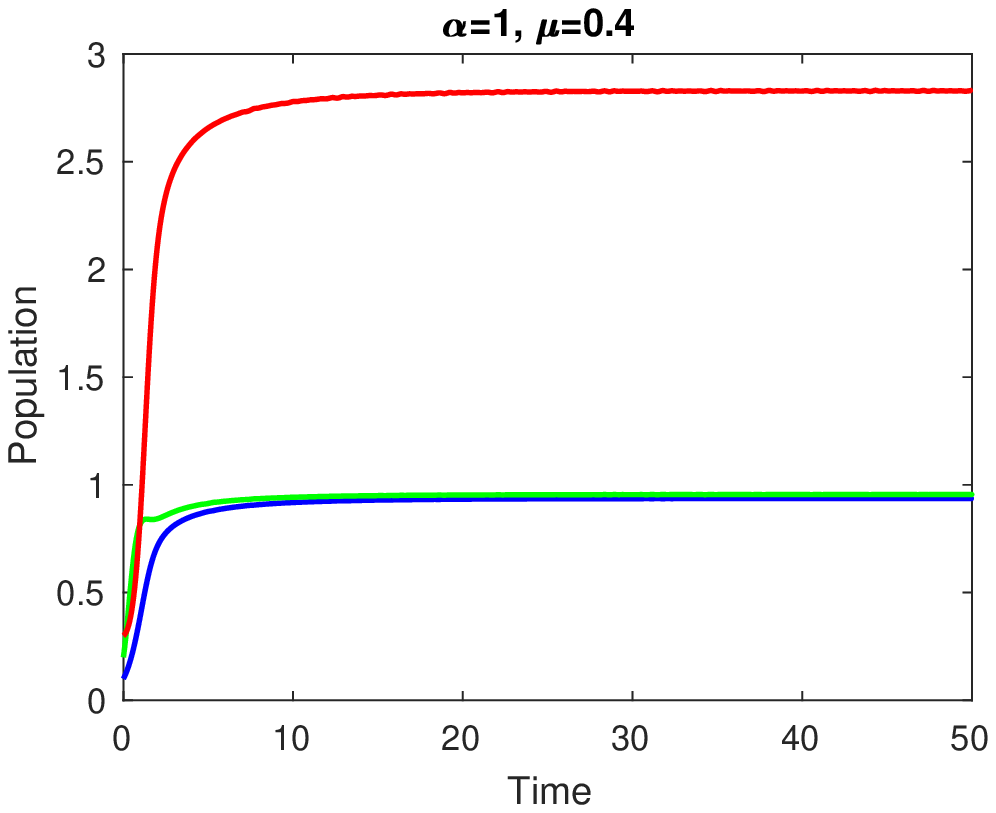}
    \includegraphics[scale=0.75]{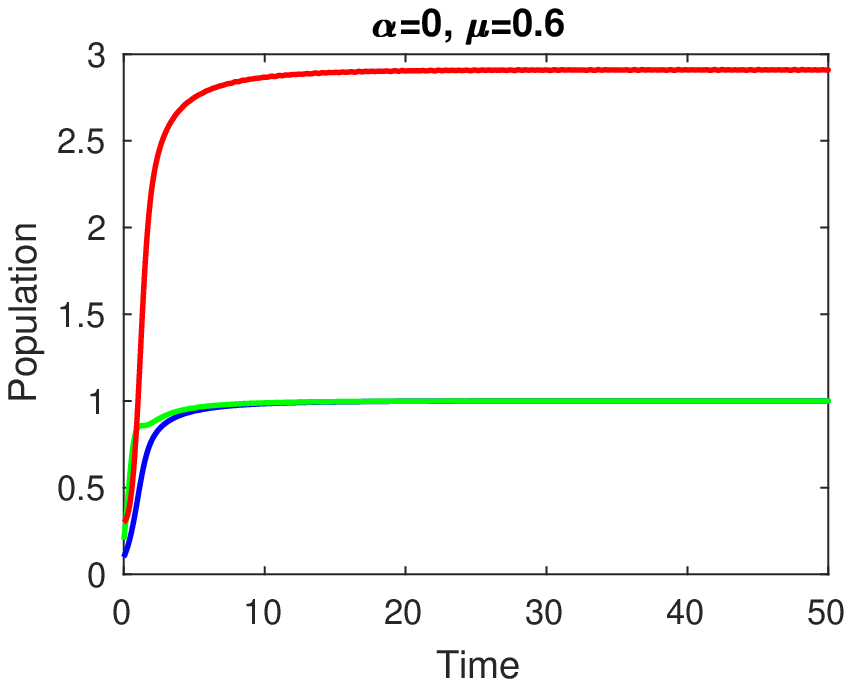} 
    \includegraphics[scale=0.75]{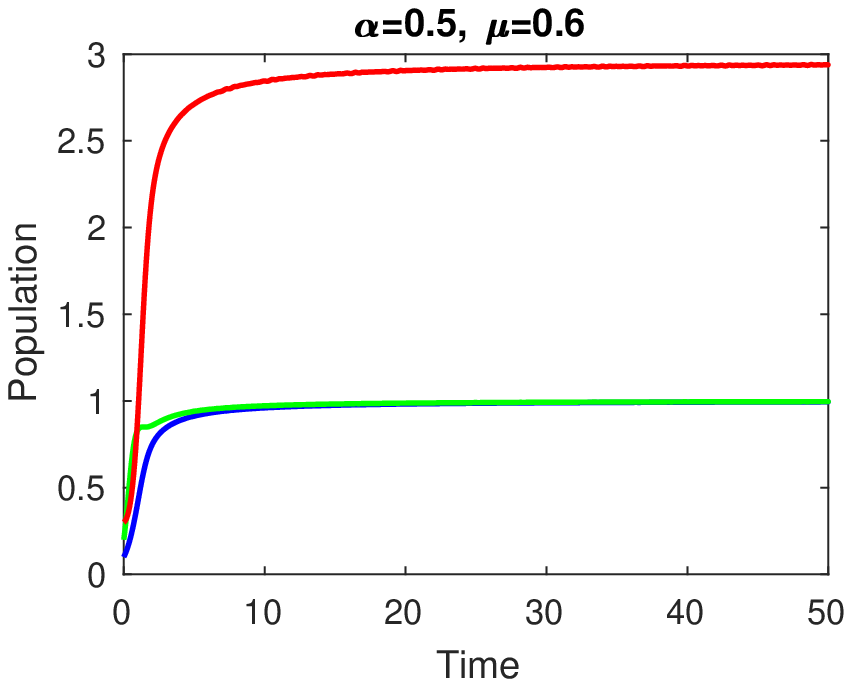}
    \end{figure}
  \begin{figure}
  \centering
    \includegraphics[scale= 0.65]{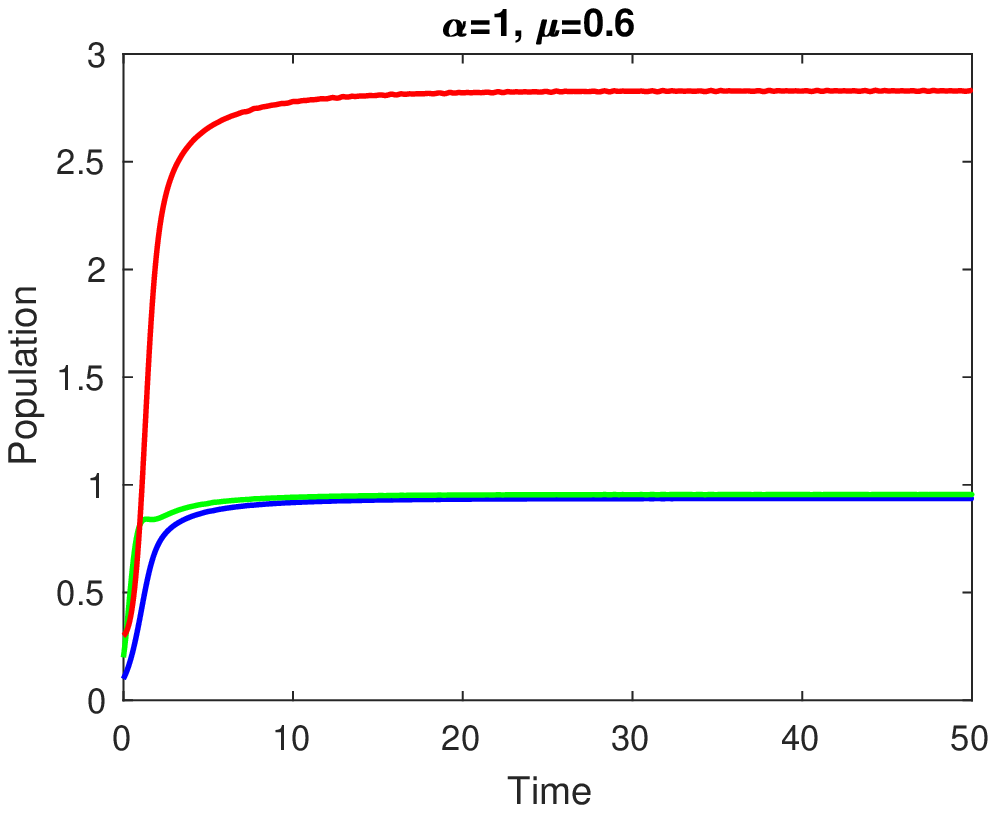}
     \includegraphics[scale=0.75]{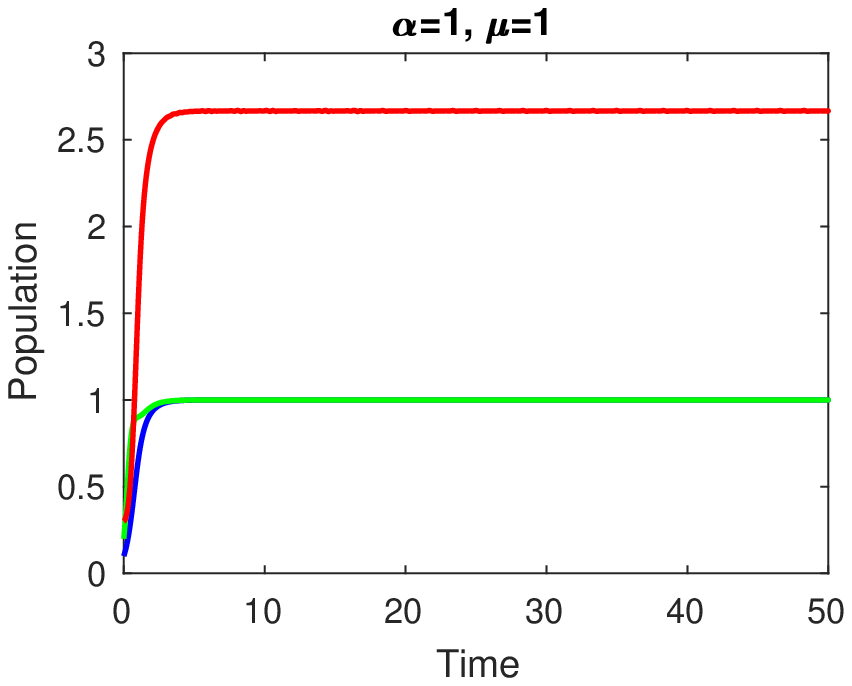} 
     \includegraphics[scale=0.75]{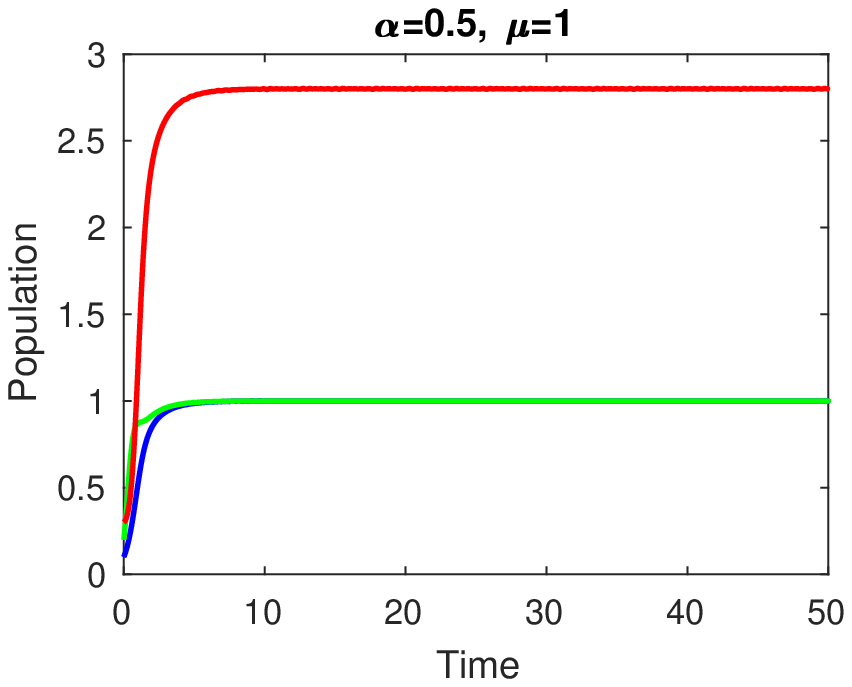}  \includegraphics[scale=0.65]{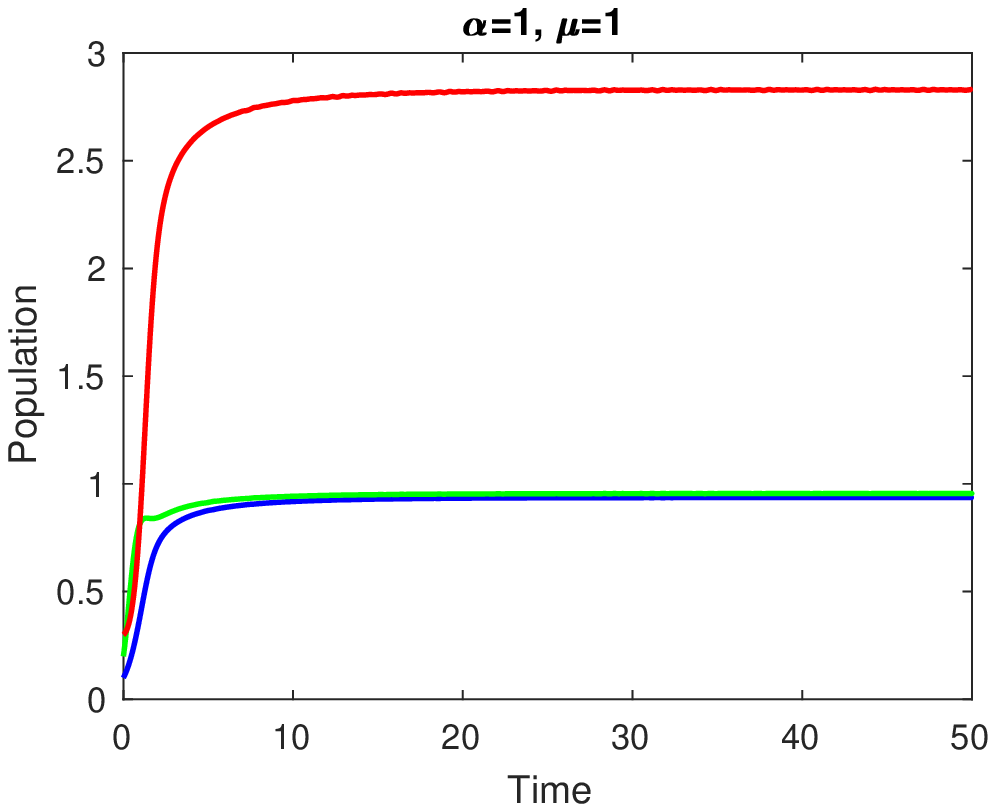}
  \caption{Variation of preys and predator populations against the time for the system (\ref{gr2a}) for different values of $\alpha,\mu$. Blue, green and red curves show the horizontal membership function of the population densities $\widehat{p}(u), \widehat{q}(u)$ and $\widehat{r}(u)$, respectively.}
  \label{grfig:4}
\end{figure}
\end{itemize}
\begin{thm}
The nontrivial fuzzy equilibrium point $\widehat{E}_7({\widehat{p_{e_7}}},{\widehat{q_{e_7}}},{\widehat{r_{e_7}}})$ is globally asymptotically fuzzy stable.
\end{thm}
\textbf{Proof:} To study the global fuzzy stability of the system (\ref{gr1}), it is enough to check that the global stability of the system (\ref{gr2}). For which, we construct a suitable Lyapunov function $\mathbb{K}(\widehat{U})=\left(\mathbb{K}(\widehat{p})-\mathbb{K}(\widehat{p_{e_7}})-\mathbb{K}(\widehat{p_{e_7}})\log\left(\frac{\mathbb{K}(\widehat{p})}{\mathbb{K}(\widehat{p_{e_7}})}\right)\right)+\left(\mathbb{K}(\widehat{q})-\mathbb{K}(\widehat{q_{e_7}})-\mathbb{K}(\widehat{q_{e_7}})\log\left(\frac{\mathbb{K}(\widehat{q})}{\mathbb{K}(\widehat{q_{e_7}})}\right)\right)+\linebreak \left(\mathbb{K}(\widehat{r})-\mathbb{K}(\widehat{r_{e_7}})-\mathbb{K}(\widehat{r_{e_7}})\log\left(\frac{\mathbb{K}(\widehat{r})}{\mathbb{K}(\widehat{r_{e_7}})}\right)\right)$. It is simple to verify that the function $\mathbb{K}(\widehat{U})$ is zero at the equilibrium point $E_7(\mathbb{K}(\widehat{p_{e_7}}),\mathbb{K}(\widehat{q_{e_7}}),\mathbb{K}(\widehat{r_{e_7}}))$ and positive for other values of\linebreak $(\mathbb{K}(\widehat{p}),\mathbb{K}(\widehat{q}),\mathbb{K}(\widehat{r}))$. Then
\begin{eqnarray}
    \begin{array}{ll}\label{grlya}
       \frac{\partial \mathbb{K}(\widehat{U})}{\partial t}=\frac{(\mathbb{K}(\widehat{p})-\mathbb{K}(\widehat{p_{e_7}}))}{\mathbb{K}(\widehat{p})}\frac{\partial \mathbb{K}(\widehat{p})}{\partial t}+\frac{(\mathbb{K}(\widehat{q})-\mathbb{K}(\widehat{q_{e_7}}))}{\mathbb{K}(\widehat{q})}\frac{\partial \mathbb{K}(\widehat{q})}{\partial t}+\frac{(\mathbb{K}(\widehat{r})-\mathbb{K}(\widehat{r_{e_7}}))}{\mathbb{K}(\widehat{r})}\frac{\partial \mathbb{K}(\widehat{r})}{\partial t}\\
    \hspace{1cm}= (\mathbb{K}(\widehat{p})-\mathbb{K}(\widehat{p_{e_7}}))(\mathbb{K}(\widehat{a_1})(1-\mathbb{K}(\widehat{p}))-
    \mathbb{K}(\widehat{r})+    \mathbb{K}(\widehat{q})    \mathbb{K}(\widehat{r}))\\
    \hspace{1cm}+ (\mathbb{K}(\widehat{q})-\mathbb{K}(\widehat{q_{e_7}}))(\mathbb{K}(\widehat{a_2})(1-\mathbb{K}(\widehat{q}))-
    \mathbb{K}(\widehat{r})+    \mathbb{K}(\widehat{p})    \mathbb{K}(\widehat{r}))\\
    \hspace{1cm}+ (\mathbb{K}(\widehat{r})-\mathbb{K}(\widehat{r_{e_7}}))(-\mathbb{K}(\widehat{a_3})\mathbb{K}(\widehat{r})+
    \mathbb{K}(\widehat{a_4})\mathbb{K}(\widehat{p})+    \mathbb{K}(\widehat{a_5})\mathbb{K}(\widehat{q})).
    \end{array}
\end{eqnarray}
For the equilibrium point $E_7(\mathbb{K}(\widehat{p_{e_7}}),\mathbb{K}(\widehat{q_{e_7}}),\mathbb{K}(\widehat{r_{e_7}}))$, we have the set of equilibrium equations
\begin{eqnarray}\label{grseteq}
    \begin{array}{ll}
      \mathbb{K}(\widehat{a_1})(1-\mathbb{K}(\widehat{p_{e_7}}))-
    \mathbb{K}(\widehat{r_{e_7}})+    \mathbb{K}(\widehat{q_{e_7}})    \mathbb{K}(\widehat{r_{e_7}})=0\\
    \mathbb{K}(\widehat{a_2})(1-\mathbb{K}(\widehat{q_{e_7}}))-
    \mathbb{K}(\widehat{r_{e_7}})+    \mathbb{K}(\widehat{p_{e_7}})    \mathbb{K}(\widehat{r_{e_7}})=0\\
-\mathbb{K}(\widehat{a_3})\mathbb{K}(\widehat{r_{e_7}})+
    \mathbb{K}(\widehat{a_4})\mathbb{K}(\widehat{p_{e_7}})+    \mathbb{K}(\widehat{a_5})\mathbb{K}(\widehat{q_{e_7}})=0.
    \end{array}
\end{eqnarray} 
Now, we simplify the Equations (\ref{grlya}) and (\ref{grseteq}) in the form
\begin{eqnarray}
    \begin{array}{ll}
       \frac{\partial \mathbb{K}(\widehat{U})}{\partial t}= -\mathbb{K}(\widehat{a_1})(\mathbb{K}(\widehat{p})-\mathbb{K}(\widehat{p_{e_7}}))^2-\mathbb{K}(\widehat{a_2})(\mathbb{K}(\widehat{q})-\mathbb{K}(\widehat{q_{e_7}}))^2
       \\
      \hspace{1cm} -\mathbb{K}(\widehat{a_3})(\mathbb{K}(\widehat{r})-\mathbb{K}(\widehat{r_{e_7}}))^2+(\mathbb{K}(\widehat{a_4})-1)(\mathbb{K}(\widehat{p})-\mathbb{K}(\widehat{p_{e_7}}))(\mathbb{K}(\widehat{r})-\mathbb{K}(\widehat{r_{e_7}}))\\
      \hspace{1cm} +\left(\frac{\mathbb{K}(\widehat{a_3})\mathbb{K}(\widehat{r_{e_7}})}{\mathbb{K}(\widehat{q_{e_7}})}-\frac{\mathbb{K}(\widehat{a_4})\mathbb{K}(\widehat{p_{e_7}})}{\mathbb{K}(\widehat{q_{e_7}})}-1\right)(\mathbb{K}(\widehat{q})-\mathbb{K}(\widehat{q_{e_7}}))(\mathbb{K}(\widehat{r})-\mathbb{K}(\widehat{r_{e_7}})).
      \end{array}
\end{eqnarray}
Next, we assume that $\mathbb{K}(\widehat{a_4})<1,\,\frac{\mathbb{K}(\widehat{a_3})\mathbb{K}(\widehat{r_{e_7}})}{\mathbb{K}(\widehat{q_{e_7}})}<\frac{\mathbb{K}(\widehat{a_4})\mathbb{K}(\widehat{p_{e_7}})}{\mathbb{K}(\widehat{q_{e_7}})}+1$ and $(\mathbb{K}(\widehat{p})-\mathbb{K}(\widehat{p_{e_7}})),\,(\mathbb{K}(\widehat{q})-\mathbb{K}(\widehat{q_{e_7}})),\,(\mathbb{K}(\widehat{r})-\mathbb{K}(\widehat{r_{e_7}}))$ are of same sign. Then $ \frac{\partial \mathbb{K}(\widehat{U})}{\partial t}<0$ except for the equilibrium point $E_7(\mathbb{K}(\widehat{p_{e_7}}),\mathbb{K}(\widehat{q_{e_7}}),\mathbb{K}(\widehat{r_{e_7}}))$. Thus the equilibrium point $E_7(\mathbb{K}(\widehat{p_{e_7}}),\mathbb{K}(\widehat{q_{e_7}}),\mathbb{K}(\widehat{r_{e_7}}))$ is globally asymptotically stable. Now, from Remark \ref{grr1} and Theorem \ref{grlogo}, we have $\mathcal{D}^{gr}\widehat{U}<0$, except for the fuzzy equilibrium point $\widehat{E}_7(\widehat{p_{e_7}},\widehat{q_{e_7}},\widehat{r_{e_7}})$. Thus the equilibrium point $\widehat{E}_7(\widehat{p_{e_7}},\widehat{q_{e_7}},\widehat{r_{e_7}})$ is globally asymptotically fuzzy stable.

\section{Numerical solution of the fuzzy prey-predator model by using granular $F$-transform}
This section presents a numerical method by using the concept of granular $F$-transform to solve the system (\ref{gr1}). To solve the system (\ref{gr1}), we first need to solve the system (\ref{gr2}) by using the concept of $F$-transform. The method of $F$-transform consists in applying the $F$-transform to both sides of the system (\ref{gr2}). It leads to a system of algebraic equations, whose solution is a discrete representation of an analytical solution of the system (\ref{gr2}). Accordingly, we apply granular $F$-transform to both sides of the system (\ref{gr1}), whereby we get the following equations:
\begin{eqnarray}\label{gr5}
\nonumber{F^{gr}[\mathcal{D}^{gr}\widehat{p}}]&=&F^{gr}[\widehat{g_1}(u,\widehat{p},\widehat{q},\widehat{r})],\\ 
\nonumber{G^{gr}[\mathcal{D}^{gr}\widehat{q}}]&=&G^{gr}[\widehat{g_2}(u,\widehat{p},\widehat{q},\widehat{r})],\\ 
{H^{gr}[\mathcal{D}^{gr}\widehat{r}}]&=& H^{gr}[\widehat{g_3}(u,\widehat{p},\widehat{q},\widehat{r})].
\end{eqnarray}
From Remark \ref{grr1}, the system (\ref{gr5}) can be written as
\begin{eqnarray}\label{gr6}
\nonumber{\mathbb{K}(F^{gr}[\mathcal{D}^{gr}\widehat{p}}])&=&\mathbb{K}(F^{gr}[\widehat{g_1}(u,\widehat{p},\widehat{q},\widehat{r})]),\\ 
\nonumber{\mathbb{K}(G^{gr}[\mathcal{D}^{gr}\widehat{q}])}&=&\mathbb{K}(G^{gr}[\widehat{g_2}(u,\widehat{p},\widehat{q},\widehat{r})]),\\ 
\nonumber{\mathbb{K}(H^{gr}[\mathcal{D}^{gr}\widehat{r}]})&=& \mathbb{K}(H^{gr}[\widehat{g_3}(u,\widehat{p},\widehat{q},\widehat{r})]).
\end{eqnarray}
From Proposition \ref{grgrdif1}, the above system can be rewritten as
\begin{eqnarray}\label{gr7}
\nonumber{F[\frac{\partial }{\partial t}\mathbb{K}(\widehat{p})]}&=&F[\mathbb{K}(\widehat{g_1}(u,\mathbb{K}(\widehat{p}),\mathbb{K}(\widehat{q}),\mathbb{K}(\widehat{r})))]=F[g_1],\\
\nonumber{G[\frac{\partial }{\partial t}\mathbb{K}(\widehat{p})]}&=&G[\mathbb{K}(\widehat{g_2}(u,\mathbb{K}(\widehat{p}),\mathbb{K}(\widehat{q}),\mathbb{K}(\widehat{r})))]=G[g_2],\\ 
{H[\frac{\partial }{\partial t}\mathbb{K}(\widehat{p})]}&=&H[\mathbb{K}(\widehat{g_3}(u,\mathbb{K}(\widehat{p}),\mathbb{K}(\widehat{q}),\mathbb{K}(\widehat{r})))]=H[g_3],
\end{eqnarray}
where 
\[g_1=\mathbb{K}(\widehat{g_1}(u,\mathbb{K}(\widehat{p}),\mathbb{K}(\widehat{q}),\mathbb{K}(\widehat{r}))),\]
\[g_2=\mathbb{K}(\widehat{g_2}(u,\mathbb{K}(\widehat{p}),\mathbb{K}(\widehat{q}),\mathbb{K}(\widehat{r}))),\]
\[g_3=\mathbb{K}(\widehat{g_3}(u,\mathbb{K}(\widehat{p}),\mathbb{K}(\widehat{q}),\mathbb{K}(\widehat{r}))).\]
 Let us choose $n > 2$ and create an $h$-uniform fuzzy partition $P_1,... P_m$ of $[a_1,a_2]$, where $h =\dfrac{a_2-a_1}{m-1}$.
 We denote the $F$-transforms of the functions on left and right sides of the system (\ref{gr7}) by 
 \[F[\mathbb{K}(\widehat{p})] = (X_1,...,X_m)^T,\]
 \[G[\mathbb{K}(\widehat{q})] = (Y_1,...,Y_m)^T ,\]
 \[H[\mathbb{K}(\widehat{r})] = (Z_1,...,Z_m)^T,\]
 \[F[\frac{\partial }{\partial t}\mathbb{K}(\widehat{p})] = (X'_1,...,X'_m)^T,\]
 \[G[\frac{\partial }{\partial t}\mathbb{K}(\widehat{q})] = (Y'_1,...,Y'_m)^T ,\]
 \[H[\frac{\partial }{\partial t}\mathbb{K}(\widehat{r})] = (Z'_1,...,Z'_m)^T ,\]
 \[F[g_1] = (F_1,...,F_m)^T,\]
 \[G[g_2] = (G_1,...,G_m)^T,\]
 \[H[g_3] = (H_1,...,H_m)^T.\]
 Also, we set \[X_1=\mathbb{K}(\widehat{p}_a),\,Y_1=\mathbb{K}(\widehat{q}_a),\,Z_1=\mathbb{K}(\widehat{r}_a),\] respectively. Thus from the system (\ref{gr7}), we obtain the following system of algebraic equations:
\begin{eqnarray}\label{gr8}
\nonumber&&X'_i=F_i,\,Y'_i=G_i,\,Z'_i=H_i,\,i=2,....m-1,\\
&&X_1=\mathbb{K}(\widehat{p}_a),\,Y_1=\mathbb{K}(\widehat{q}_a),Z_1=\mathbb{K}(\widehat{r}_a),
\end{eqnarray}
where \[X'_i=F_i[\frac{\partial}{\partial t}\mathbb{K}(\widehat{p})],\,Y'_i=G_i[\frac{\partial }{\partial t}\mathbb{K}(\widehat{q})],\,Z'_i=H_i[\frac{\partial }{\partial t}\mathbb{K}(\widehat{r})],\]
\[X_i=F_i[\mathbb{K}(\widehat{p})],\,Y_i=G_i[\mathbb{K}(\widehat{q})],\,Z_i=H_i[\mathbb{K}(\widehat{r})],\] \[F_i=F_i[g_1],\,G_i=G_i[g_2],\,H_i=H_i[g_3]\]
are $F$-transform components. Also, by using schemes of the method of central differences, we replace 
\begin{eqnarray*}
X'_i&=&\dfrac{X_{i+1}-X_{i-1}}{2h},\\
Y'_i&=&\dfrac{Y_{i+1}-Y_{i-1}}{2h},\\
Z'_i&=&\dfrac{Z_{i+1}-Z_{i-1}}{2h},\,i=2,...,m-1
\end{eqnarray*}
and assume $X_2=X_1+hF_1,Y_2=Y_1+hG_1,Z_2=Z_1+hH_1$. Thus we can introduce $(m-2)\times m$ matrix
\begin{equation}\label{grm2}
{D} = \frac{1}{2h}
\begin{bmatrix}
-1 & 0 &1&0&0&0&...&0\\
0 & -1 &0&1&0&0&...&0\\
0 & 0&-1&0&1&0&...&0\\
 .& &&&.&& &.\\[-.2cm]
 . & &&&&.& &.\\[-.2cm]
 . & &&&&&. &.\\
  0 & 0&0&0&...&-1&0 &1\\
\end{bmatrix}.
\end{equation}
Now, we can rewrite the system (\ref{gr8}) to the following system of linear equations:
\[D F[\mathbb{K}(\widehat{p})]=F[g_1],\,D G[\mathbb{K}(\widehat{q})]=G[g_2],\,D H[\mathbb{K}(\widehat{r})]=H[g_3],\]
where \[
F[\mathbb{K}(\widehat{p})] = (X_1,...,X_{m-1})^T,\]
\[G[\mathbb{K}(\widehat{q})] = (Y_1,...,Y_{m-1})^T,\] \[H[\mathbb{K}(\widehat{r})] = (Z_1,...,Z_{m-1})^T,\]
\[F[\frac{\partial }{\partial t}\mathbb{K}(\widehat{p})] = (X'_2,...,X'_{m-1})^T,\]\[G[\frac{\partial }{\partial t}\mathbb{K}(\widehat{q})] = (Y'_2,...,Y'_{m-1})^T,\]\[H[\frac{\partial }{\partial t}\mathbb{K}(\widehat{r})] = (Z'_2,...,Z'_{m-1})^T,\]
\[F[g_1] = (F_2,...,F_{m-1})^T,\]\[ G[g_2] = (G_2,...,G_{m-1})^T,\]\[H[g_3] = (H_2,...,H_{m-1})^T.\] Next, the matrix $D$ is completed by adding the first and second row as initial values, as seen below:
\begin{equation*}\label{grm3}
{D'} = \frac{1}{2h}
\begin{bmatrix}
1 & 0 &0&0&0&0&...&0\\
0& 1 &0&0&0&0&...&0\\
-1 & 0 &1&0&0&0&...&0\\
0 & -1 &0&1&0&0&...&0\\
0 & 0&-1&0&1&0&...&0\\
 .& &&&.&& &.\\[-.2cm]
 . & &&&&.& &.\\[-.2cm]
 . & &&&&&. &.\\
  0 & 0&0&0&...&-1&0 &1\\
\end{bmatrix},
\end{equation*}
so that the matrix $D'$ is nom-singular. Now, we can expand $F[ g_1 ],G[g_2]$ and $G[g_3]$ by adding the first and second elements using initial conditions and the matrix $D'$, as follows:
\[F[g_1] = \left(\frac{\mathbb{K}(\widehat{p}_a)}{2h},\frac{\mathbb{K}(\widehat{p_{2}})}{2h},g_2,...,F_{m-1}\right)^T,\]
\[G[g_2] = \left(\frac{\mathbb{K}(\widehat{q}_a)}{2h},\frac{\mathbb{K}(\widehat{q}_2)}{2h},G_2,...,G_{m-1}\right)^T, \] 
\[H[g_3] = \left(\frac{\mathbb{K}(\widehat{r}_a)}{2h},\frac{\mathbb{K}(\widehat{r}_2)}{2h},H_2,...,H_{m-1}\right)^T,\]
where $\mathbb{K}(\widehat{p_{2}})=X_2$, $\mathbb{K}(\widehat{q}_2)=Y_2$ and $\mathbb{K}(\widehat{r}_2)=Z_2$. Thus we have 
\begin{eqnarray}\label{grmat}
D' F[\mathbb{K}(\widehat{p})]=F[g_1],\,D' G[\mathbb{K}(\widehat{q})]=G[g_2],\,D' H[\mathbb{K}(\widehat{r})]=H[g_3].
\end{eqnarray}
The system (\ref{grmat}) can be written as
\begin{eqnarray}\label{grmat1}
 F[\mathbb{K}(\widehat{p})]=(D')^{-1}F[g_1],\, G[\mathbb{K}(\widehat{q})]=(D')^{-1}G[g_2],\, H[\mathbb{K}(\widehat{r})]=(D')^{-1}H[g_3].
\end{eqnarray}
Now, to solve the system (\ref{grmat1}), we compute the inverse matrix $(D')^{-1}$. If $m$ is even then inverse matrix is given as
\begin{equation*}\label{grm4}
({D'})^{-1} = {2h}
\begin{bmatrix}
1 & 0 &0&0&0&0&...&0\\
0& 1 &0&0&0&0&...&0\\
1 & 0 &1&0&0&0&...&0\\
0 & 1 &0&1&0&0&...&0\\
1 & 0&1&0&1&0&...&0\\
 .& &&&.&& &.\\[-.2cm]
 . & &&&&.& &.\\[-.2cm]
 . & &&&&&. &.\\
  0 & ..&0&1&0&1&0 &1\\
\end{bmatrix},
\end{equation*}
and when $m$ is odd, we obtain
\begin{equation*}\label{grm5}
({D'})^{-1} = {2h}
\begin{bmatrix}
1 & 0 &0&0&0&0&...&0\\
0& 1 &0&0&0&0&...&0\\
1 & 0 &1&0&0&0&...&0\\
0 & 1 &0&1&0&0&...&0\\
 .& &&&.&& &.\\[-.2cm]
 . & &&&&.& &.\\[-.2cm]
 . & &&&&&. &.\\
  0 & ..&0&1&0&1&0 &1\\
   1 & 0&1&0&...&1&0 &1\\
\end{bmatrix}.
\end{equation*}
Thus from the system (\ref{grmat1}), we have the following \begin{eqnarray}\label{gr9}
\nonumber X_{i+1}&=&X_{i-1}+2hF_i,\\
\nonumber Y_{i+1}&=&Y_{i-1}+2hG_i,\\
\nonumber Z_{i+1}&=&Z_{i-1}+2hH_i,\\
\nonumber X_2&=&\mathbb{K}(\widehat{p_{2}}),\,Y_2=\mathbb{K}(\widehat{q}_2),\,
Z_2=\mathbb{K}(\widehat{r}_2),\\
X_1&=&\mathbb{K}(\widehat{p}_a),\,Y_1=\mathbb{K}(\widehat{q}_a),Z_1=\mathbb{K}(\widehat{r}_a),\,i=2,...,m-1.
\end{eqnarray}
The system (\ref{gr9}) can be used to the computation of $X_3, . . . , X_m,Y_3, . . . , Y_m$ and $Z_3, . . . , Z_m$. However, it can not be used directly by using the functions $\mathbb{K}(\widehat{g_1}(u,\mathbb{K}(\widehat{p}),\mathbb{K}(\widehat{q}),\mathbb{K}(\widehat{r}))),$ $\mathbb{K}(\widehat{g_2}(u,\mathbb{K}(\widehat{p}),\mathbb{K}(\widehat{q}),\mathbb{K}(\widehat{r})))$ and $\mathbb{K}(\widehat{g_3}(u,\mathbb{K}(\widehat{p}),\mathbb{K}(\widehat{q}),\mathbb{K}(\widehat{r})))$, because these functions use unknown functions $\mathbb{K}(\widehat{p}), \mathbb{K}(\widehat{q})$ and $\mathbb{K}(\widehat{r})$. Therefore we can use the same technique as in \cite{p} and substitute functions by their $F$-transform components:
\begin{eqnarray}
\nonumber\hat{F}_i[g_1]&=&\dfrac{\int_{a}^{b} \mathbb{K}(\widehat{g_1}(u,X_i,Y_i,Z_i))P_i(u)du}{\int_{a}^{b}P_i(u)du },\\
\nonumber\hat{G}_i[g_2]&=&\dfrac{\int_{a}^{b} \mathbb{K}(\widehat{g_2}(u,X_i,Y_i,Z_i))P_i(u)du}{\int_{a}^{b}P_i(u)du },\\
\hat{H}_i[g_3]&=&\dfrac{\int_{a}^{b} \mathbb{K}(\widehat{g_3}(u,X_i,Y_i,Z_i))P_i(u)du}{\int_{a}^{b}P_i(u)du },i=1,...,m-1.
\end{eqnarray}
Thus the system (\ref{gr9}) can be written as
\begin{eqnarray}\label{gr10}
\nonumber X_{i+1}&=&X_{i-1}+2h\hat{F}_i[g_1],\\
\nonumber Y_{i+1}&=&Y_{i-1}+2h\hat{G}_i[g_2],\\
\nonumber Z_{i+1}&=&Z_{i-1}+2h\hat{H}_i[g_3],\\
\nonumber X_{2}&=&X_{1}+h\hat{F}_1[g_1],\\
\nonumber Y_{2}&=&Y_{1}+h\hat{G}_1[g_2],\\
\nonumber Z_{2}&=&Z_{1}+h\hat{H}_1[g_3],\,i=2,...,m-1,\\
X_1&=&\mathbb{K}(\widehat{p}_a),\,Y_1=\mathbb{K}(\widehat{q}_a),Z_1=\mathbb{K}(\widehat{r}_a).
\end{eqnarray}
The proposed method is based on the same assumptions as the well-known Euler mid-point method. It has the same degree of accuracy. The solution vector $(F_m[\mathbb{K}(\widehat{p})],G_m[\mathbb{K}(\widehat{q})],H_m[\mathbb{K}(\widehat{r})])$ approximates the (unique) solution $(\mathbb{K}(\widehat{p}),\mathbb{K}(\widehat{q}),\mathbb{K}(\widehat{r}))$ of the system (\ref{gr2}) in the sense that $(\mathbb{K}(\widehat{p}_i),\mathbb{K}(\widehat{q}_i),\mathbb{K}(\widehat{r}_i))\approx (X_i,Y_i,Z_i), i= 1,...,m$, where nodes
$u_1,...,u_m$ are determined by fuzzy partition $P_1,...,P_m$.\\\\
The system (\ref{gr10}) can be written as in terms of granular $F$-transform, i.e.,
\begin{eqnarray}\label{gr11}
\nonumber\mathbb{K}(F^{gr}_{i+1}[\widehat{p}])&=&\mathbb{K}(F^{gr}_{i-1}[\widehat{p}])+2h\mathbb{K}(\hat{F}^{gr}_i[\widehat{g_1}]),\\
\nonumber\mathbb{K}(G^{gr}_{i+1}[\widehat{q}])&=&\mathbb{K}(G^{gr}_{i-1}[\widehat{q}])+2h\mathbb{K}(\hat{G}^{gr}_i[\widehat{g_2}]),\\
\nonumber\mathbb{K}(H^{gr}_{i+1}[\widehat{r}])&=&\mathbb{K}(H^{gr}_{i-1}[\widehat{r}])+2h\mathbb{K}(\hat{H}^{gr}_i[\widehat{g_3}]),\\
\nonumber\mathbb{K}(F^{gr}_{2}[\widehat{p}])&=&\mathbb{K}(F^{gr}_{1}[\widehat{p}])+h\mathbb{K}(\hat{F}^{gr}_{1}[{\widehat{g_1}}]),\\
\nonumber\mathbb{K}(G^{gr}_{2}[\widehat{q}])&=&\mathbb{K}(G^{gr}_{1}[\widehat{q}])+h\mathbb{K}(\hat{G}^{gr}_{1}[{\widehat{g_2}}]),\\
\nonumber\mathbb{K}(H^{gr}_{2}[\widehat{r}])&=&\mathbb{K}(H^{gr}_{1}[\widehat{r}])+h\mathbb{K}(\hat{H}^{gr}_{1}[{\widehat{g_3}}]),\\
\mathbb{K}(F^{gr}_{1}[\widehat{p}])&=&\mathbb{K}(\widehat{p}_a),\mathbb{K}(G^{gr}_{1}[\widehat{q}])=\mathbb{K}(\widehat{q}_a),\mathbb{K}(H^{gr}_{1}[\widehat{r}])=\mathbb{K}(\widehat{r}_a),
\end{eqnarray}
where 
\begin{eqnarray*}
\widehat{g_1}&=&\widehat{g_1}(F^{gr}_i[\widehat{p}],G^{gr}_i[\widehat{q}],H^{gr}_i[\widehat{r}]),\\ \widehat{g_2}&=&\widehat{g_2}(F^{gr}_i[\widehat{p}],G^{gr}_i[\widehat{q}],H^{gr}_i[\widehat{r}]),\\ \widehat{g_3}&=&\widehat{g_3}(F^{gr}_i[\widehat{p}],G^{gr}_i[\widehat{q}],H^{gr}_i[\widehat{r}]),\,i=2,...,m-1. 
\end{eqnarray*}
From Reamrk \ref{grr1}, the system (\ref{gr11}) becomes
\begin{eqnarray}\label{gr12}
\nonumber F^{gr}_{i+1}[\widehat{p}]&=&F^{gr}_{i-1}[\widehat{p}] +2h \hat{F}^{gr}_i[\widehat{g_1}],\\
\nonumber G^{gr}_{i+1}[\widehat{q}]&=&G^{gr}_{i-1}[\widehat{q}] +2h \hat{G}^{gr}_i[\widehat{g_2}],\\
\nonumber H^{gr}_{i+1}[\widehat{r}]&=&H^{gr}_{i-1}[\widehat{r}] +2h \hat{H}^{gr}_i[\widehat{g_3}],\\
\nonumber F^{gr}_{2}[\widehat{p}]&=&F^{gr}_{1}[\widehat{p}])+h\hat{F}^{gr}_{1}[{\widehat{g_1}}],\\
\nonumber G^{gr}_{2}[\widehat{q}]&=&G^{gr}_{1}[\widehat{q}])+h\hat{G}^{gr}_{1}[{\widehat{g_2}}],\\
 H^{gr}_{2}[\widehat{r}]&=&H^{gr}_{1}[\widehat{r}]+h\hat{H}^{gr}_{1}[{\widehat{g_3}}],\,i=2,...,m-1
\end{eqnarray}
where $F^{gr}_{1}[\widehat{p}]=\widehat{p}_a,G^{gr}_{1}[\widehat{q}]=\widehat{q}_a,H^{gr}_{1}[\widehat{r}]=\widehat{r}_a$. The solution vector $(F^{gr}[\widehat{p}],G^{gr}[\widehat{q}],H^{gr}[\widehat{r}])$ approximates the (unique) solution $(\widehat{p},\widehat{q},\widehat{r})$ of the system (\ref{gr1}) in the sense that $(\widehat{p}_i,\widehat{q}_i,\widehat{r}_i)\approx (F^{gr}_i,G^{gr}_i,H^{gr}_i), i= 1,...,m$, where nodes
$u_1,...,u_m$ are determined by fuzzy partition $P_1,...,P_m$. Applying the granular inverse $F$-transform to the extended vector $(F^{gr}[\widehat{p}],G^{gr}[\widehat{q}],H^{gr}[\widehat{r}])$,
we obtain an approximate solution $(\widehat{p},\widehat{q},\widehat{r})$ of the system (\ref{gr1}) in the form of fuzzy function function:
\begin{eqnarray*}
\widehat{{p}}^{gr}_m(u)&=&\sum_{n=1}^{m} F^{gr}_{i}[\widehat{p}]P_i(u),\\
\widehat{{q}}^{gr}_m(u)&=&\sum_{n=1}^{m} G^{gr}_{i}[\widehat{q}]P_i(u),\\ 
\widehat{{r}}^{gr}_m(u)&=&\sum_{n=1}^{m} H^{gr}_{i}[\widehat{r}]P_i(u),
\end{eqnarray*}
where $P_i(u), i= 1, . . . ,m$ are given basic functions. 
\begin{exa}\label{grex5}
\end{exa} Let $\widehat{a_1}=(0.01,0.02,0.03),\widehat{a_2}=(2,4,6),\widehat{a_3}=(0.1,0.2,0.3),\widehat{a_4}=(3,4,5),\widehat{a_5}=(1,2,3),\widehat{p}_0=0.1,\widehat{q}_0=0.2,\widehat{r}_0=0.3$, where $\widehat{a_1},\widehat{a_2},\widehat{a_3},\widehat{a_4},\widehat{a_5}$ are triangular fuzzy numbers and assume that fuzzy partition is a triangular fuzzy partition. Now, the horizontal membership functions of the given triangular fuzzy numbers are given by
\begin{eqnarray*}
\mathbb{K}(\widehat{a_1})&=&a_1^{gr}(\alpha,\mu_{a_1})=0.01+0.01\alpha+\mu_{a_1}(0.02-0.02\alpha),\\
\mathbb{K}(\widehat{a_2})&=&a_2^{gr}(\alpha,\mu_{a_2})=2+2\alpha+\mu_{a_2}(4-4\alpha),\\
\mathbb{K}(\widehat{a_3})&=&a_3^{gr}(\alpha,\mu_{a_3})=0.1+0.1\alpha+\mu_{a_3}(0.2-0.2\alpha),\\
\mathbb{K}(\widehat{a_4})&=&a_4^{gr}(\alpha,\mu_{a_4})=3+\alpha+\mu_{a_4}(2-2\alpha),\\
\mathbb{K}(\widehat{a_5})&=&a_5^{gr}(\alpha,\mu_{a_5})=1+\alpha+\mu_{a_5}(2-2\alpha),\\
\mathbb{K}(\widehat{p}_0)&=&0.1,\mathbb{K}(\widehat{q}_0)=0.2,\mathbb{K}(\widehat{r}_0)=0.3.
\end{eqnarray*}
Further, we assume $\mu_p,\mu_q,\mu_r,\mu_{a_1},\mu_{a_2},\mu_{a_3},\mu_{a_4},\mu_{a_5}=\mu\in\{0,0.4,0.6,1\}$, and $\alpha\in\{0,0.5,1\}$. A comparison of the numerical solutions obtained by F-transform (FT-Euler mid-point method) and Euler method with exact solution of the system (\ref{gr2}) corresponding to different values of $\alpha,\mu$ and $h=0.01$ . For the accuracy of numerical solutions over exact solution, root mean square error (RMS error) is used, which is defined by formula
\[\mbox{RMS error}=\sqrt{\frac{1}{m}\sum_{i=1}^{m}({\Delta}y_k)^2},\] where ${\Delta}y_i=y_k^{exact}-y_k^{numerical}$ and $y_k^{exact},y_k^{numerical}$ are exact and numerical solution, respectively. Tables \ref{grtab:5}, \ref{grtab:6}, \ref{grtab:7} and \ref{grtab:8} show comparison among the exact solution, Euler method and the FT-Euler mid-point method for the system (\ref{gr2}). The graphical representation of this system shows that exact solution and the FT-Euler mid-point method are very close to each other in comparison to Euler method. Therefore from above comparison, we can say that FT-Euler mid-point method is a trustworthy numerical method and the above figures show that the population density $\mathbb{K}(\widehat{p}(u))$ of prey decreases whereas the population densities $\mathbb{K}(\widehat{q}(u)),\mathbb{K}(\widehat{r}(u))$ of prey and predator, respectively increase corresponding to different values of $\alpha,\mu$.
\begin{figure}
  \centering
  \includegraphics[scale=0.5]{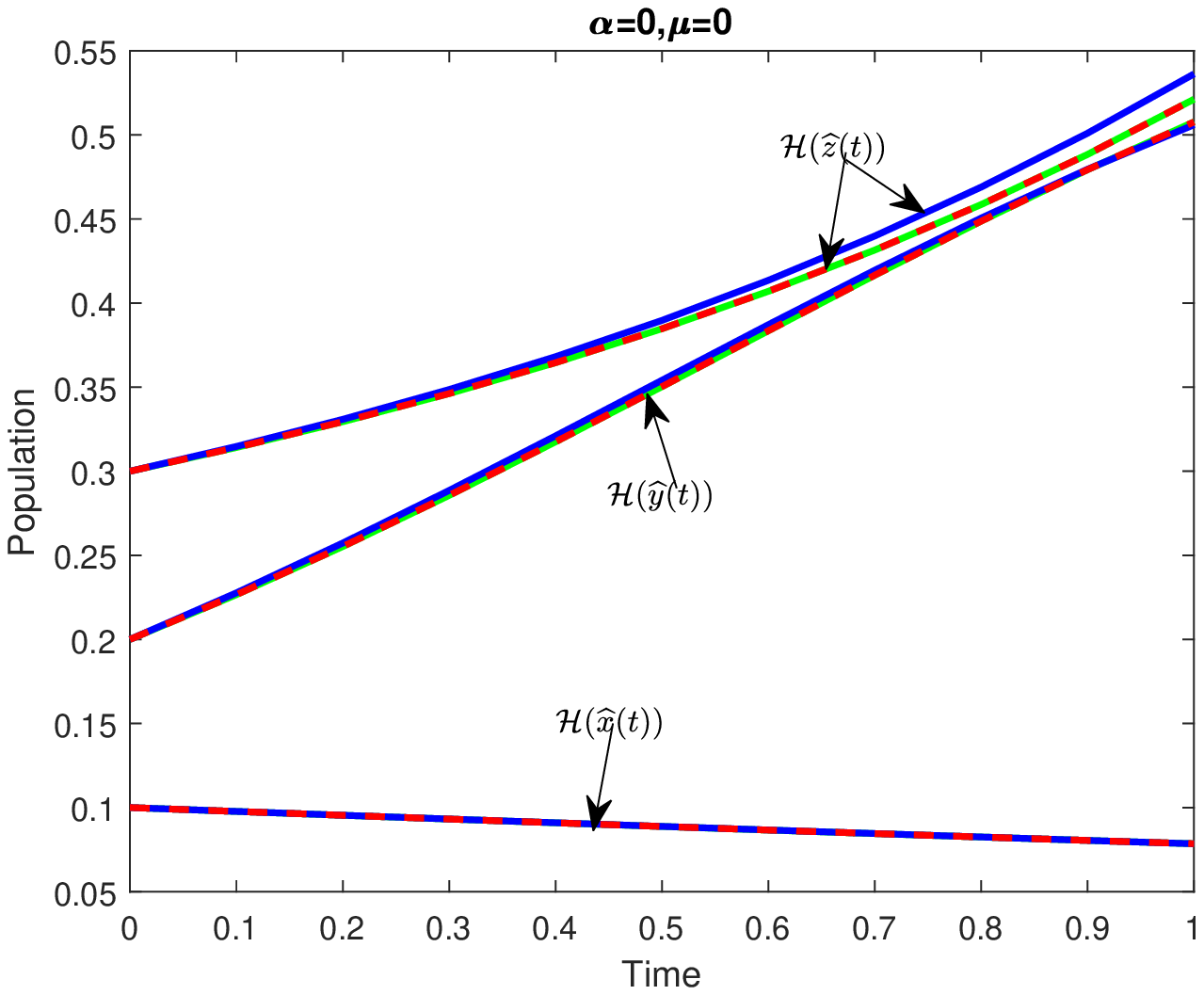} \includegraphics[scale=0.5]{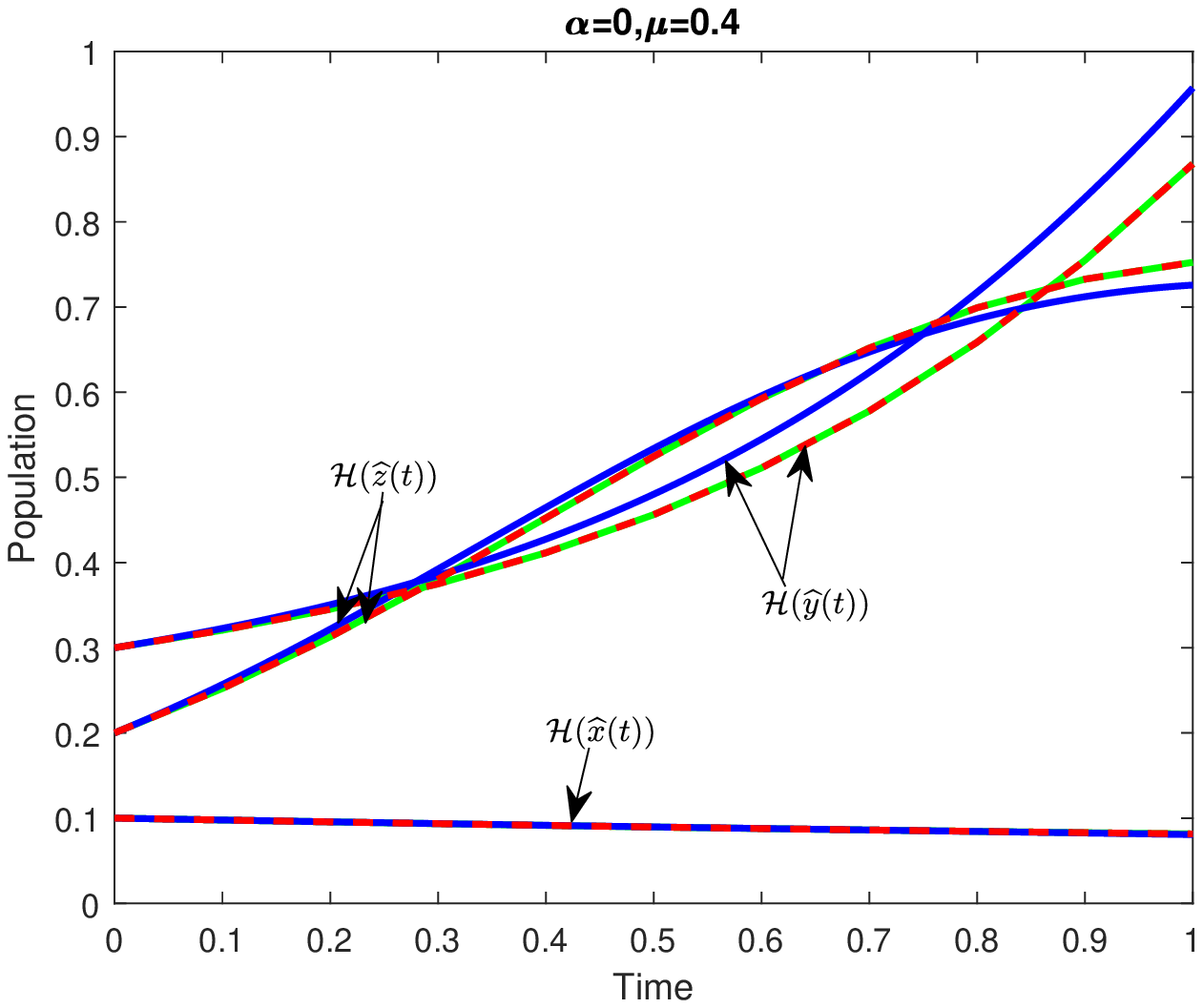}
   \includegraphics[scale=0.5]{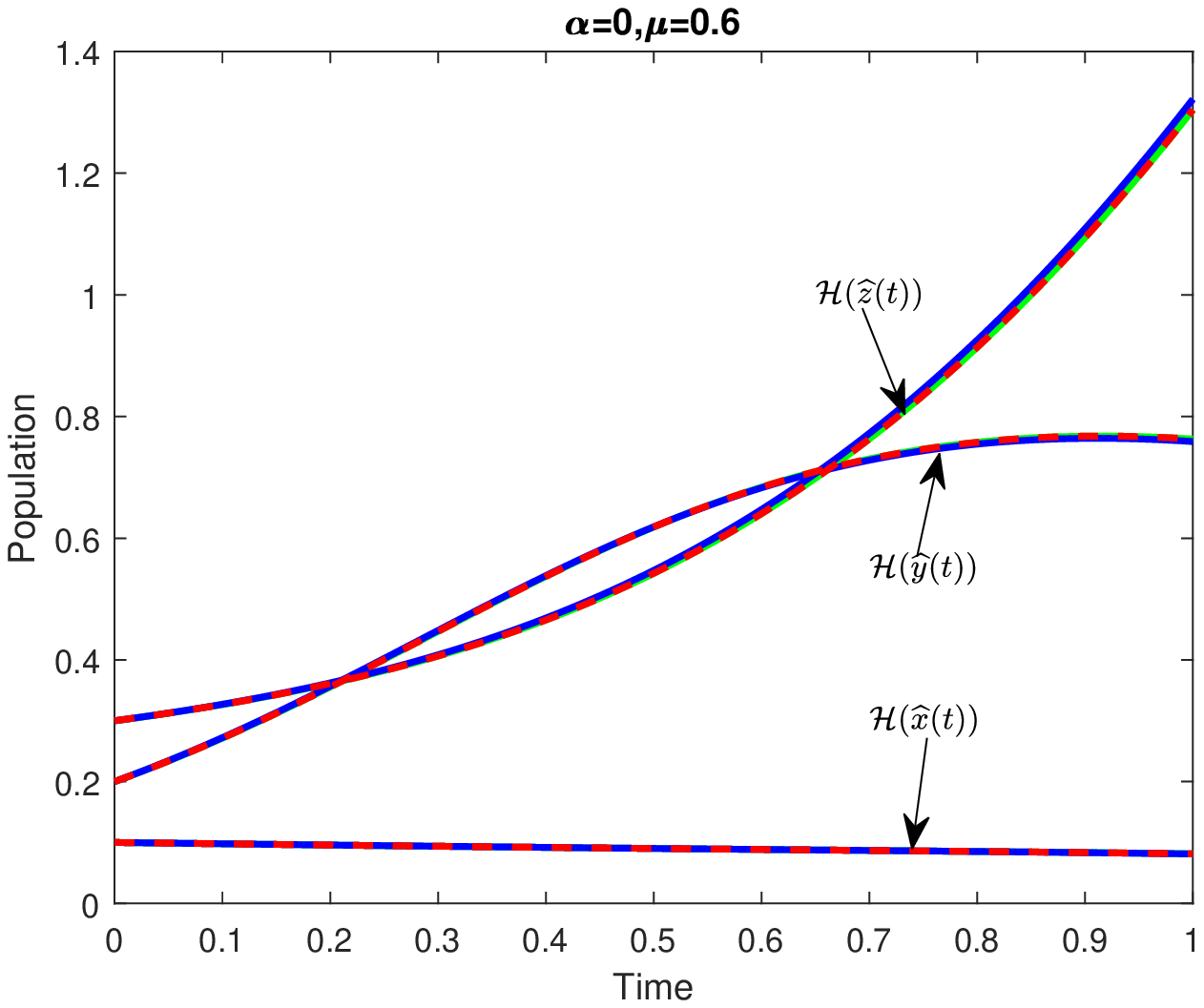}
   \includegraphics[scale=0.5]{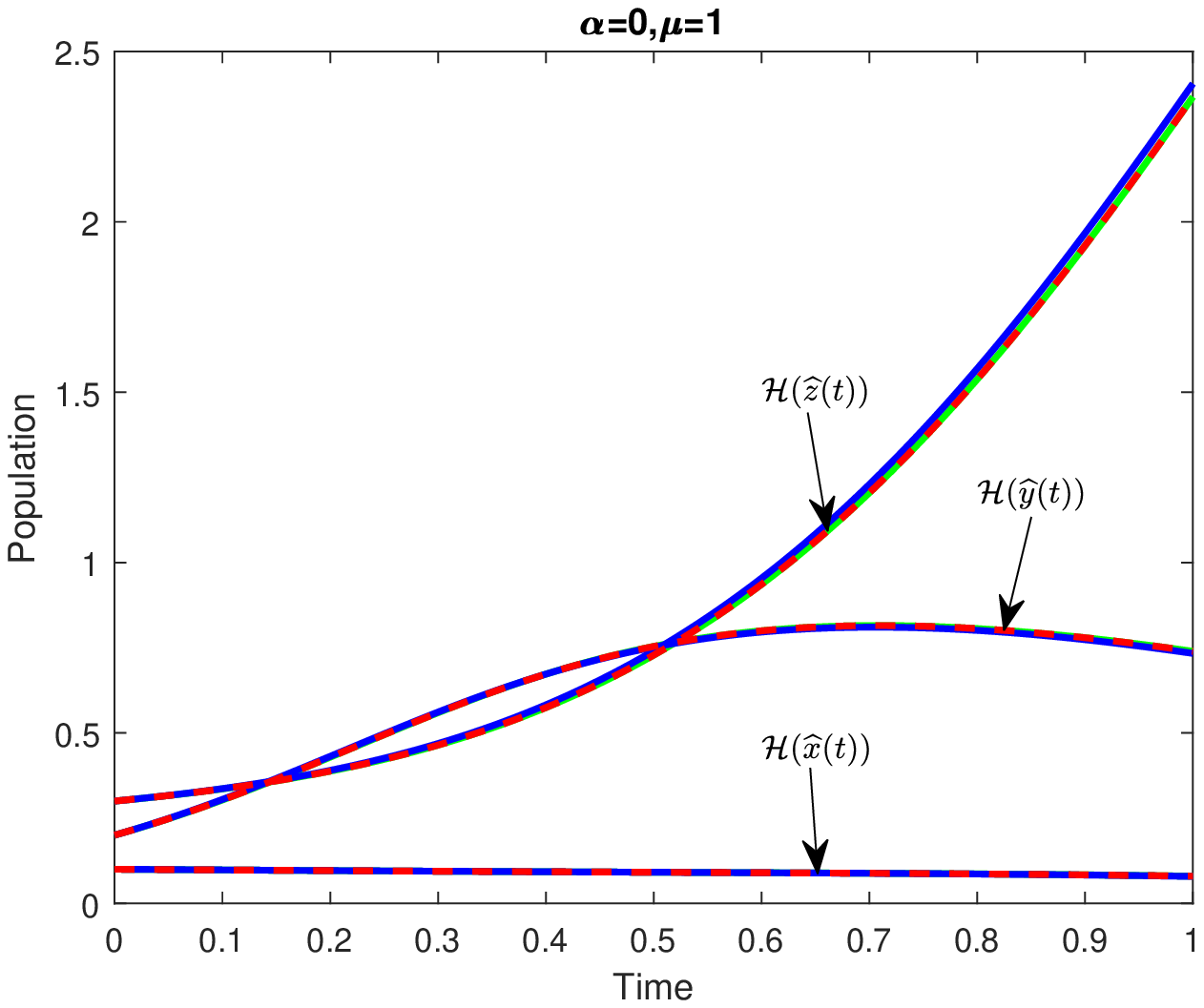}
   \includegraphics[scale=0.5]{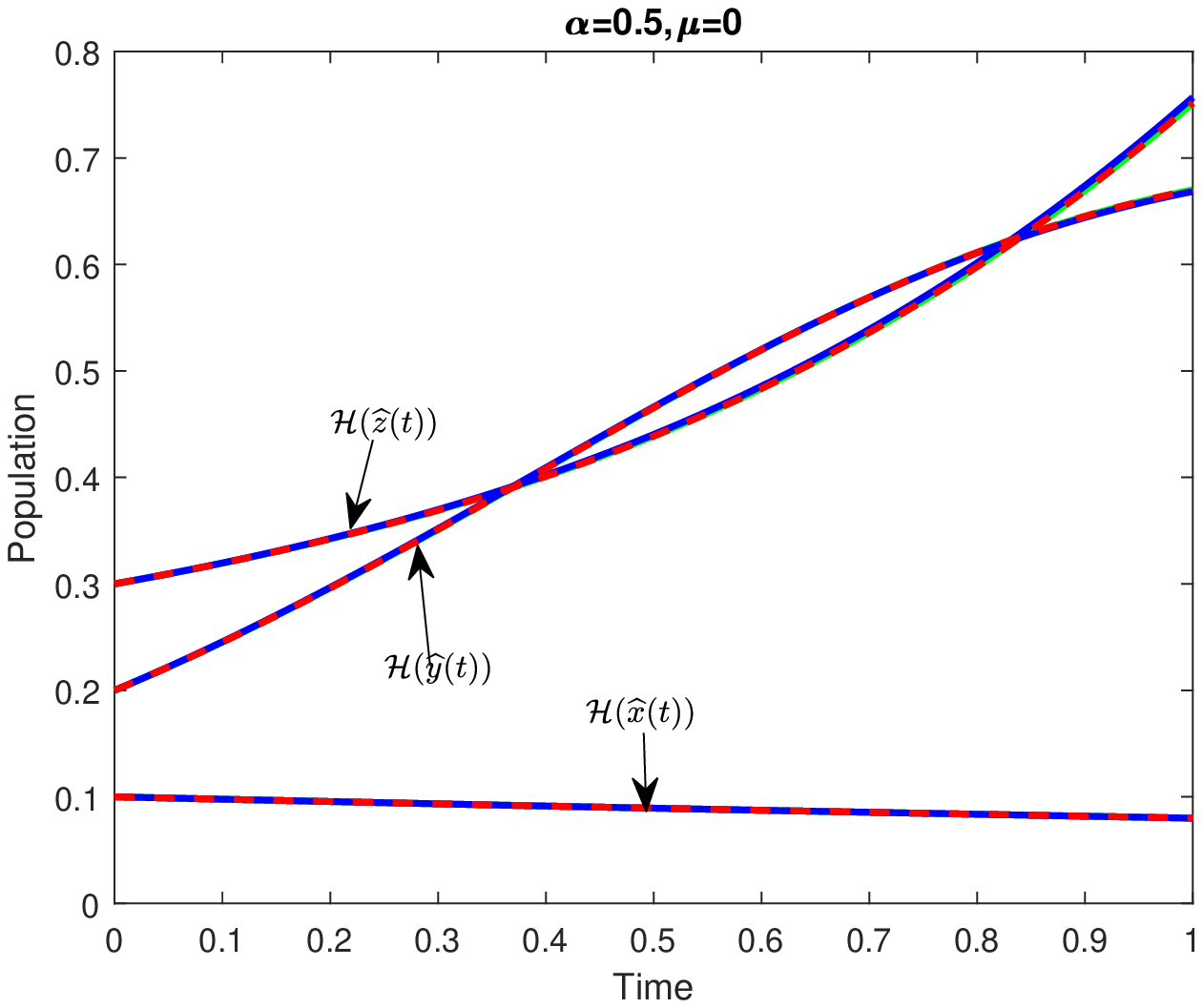}
  \includegraphics[scale=0.5]{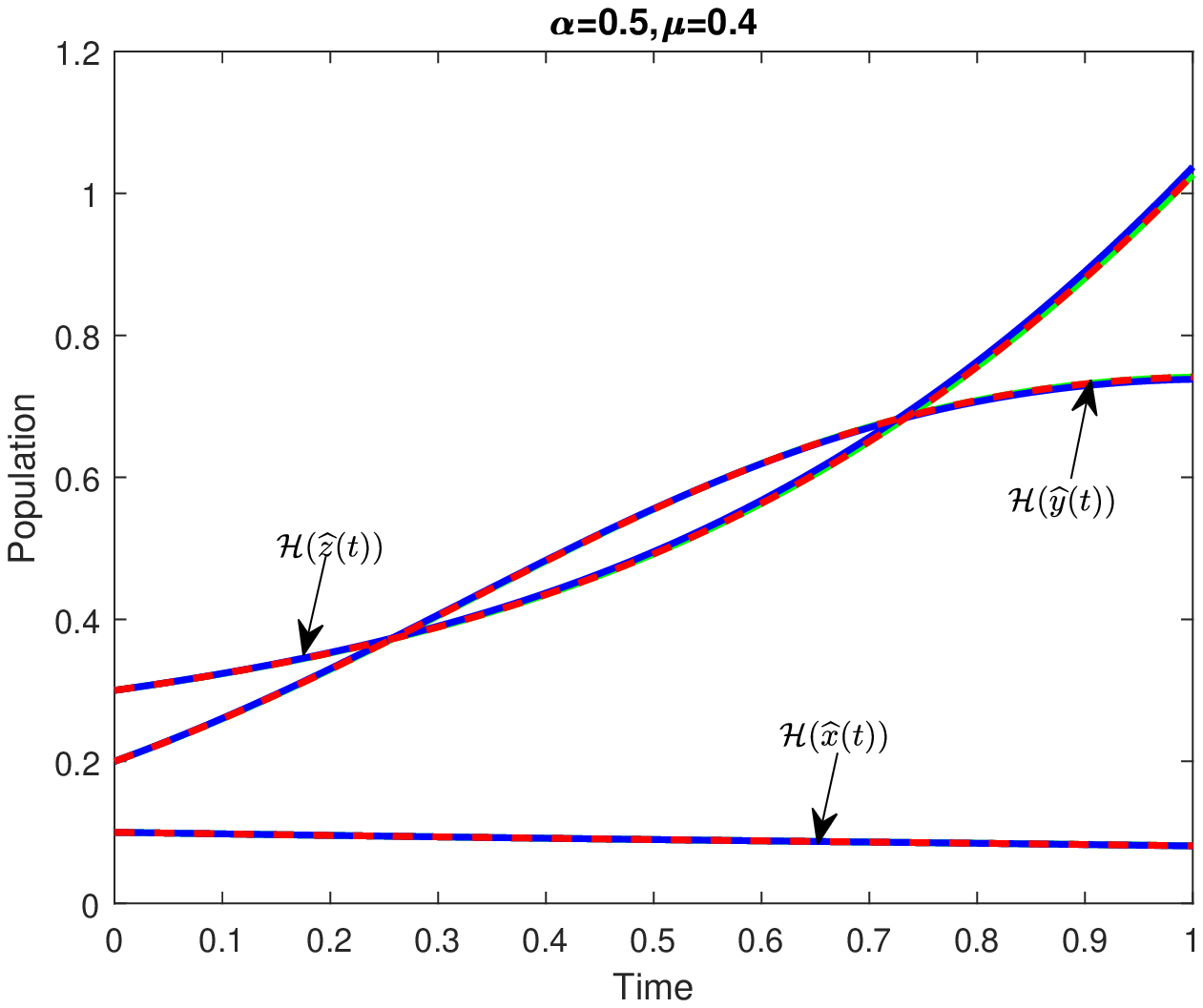}
  \includegraphics[scale=0.5]{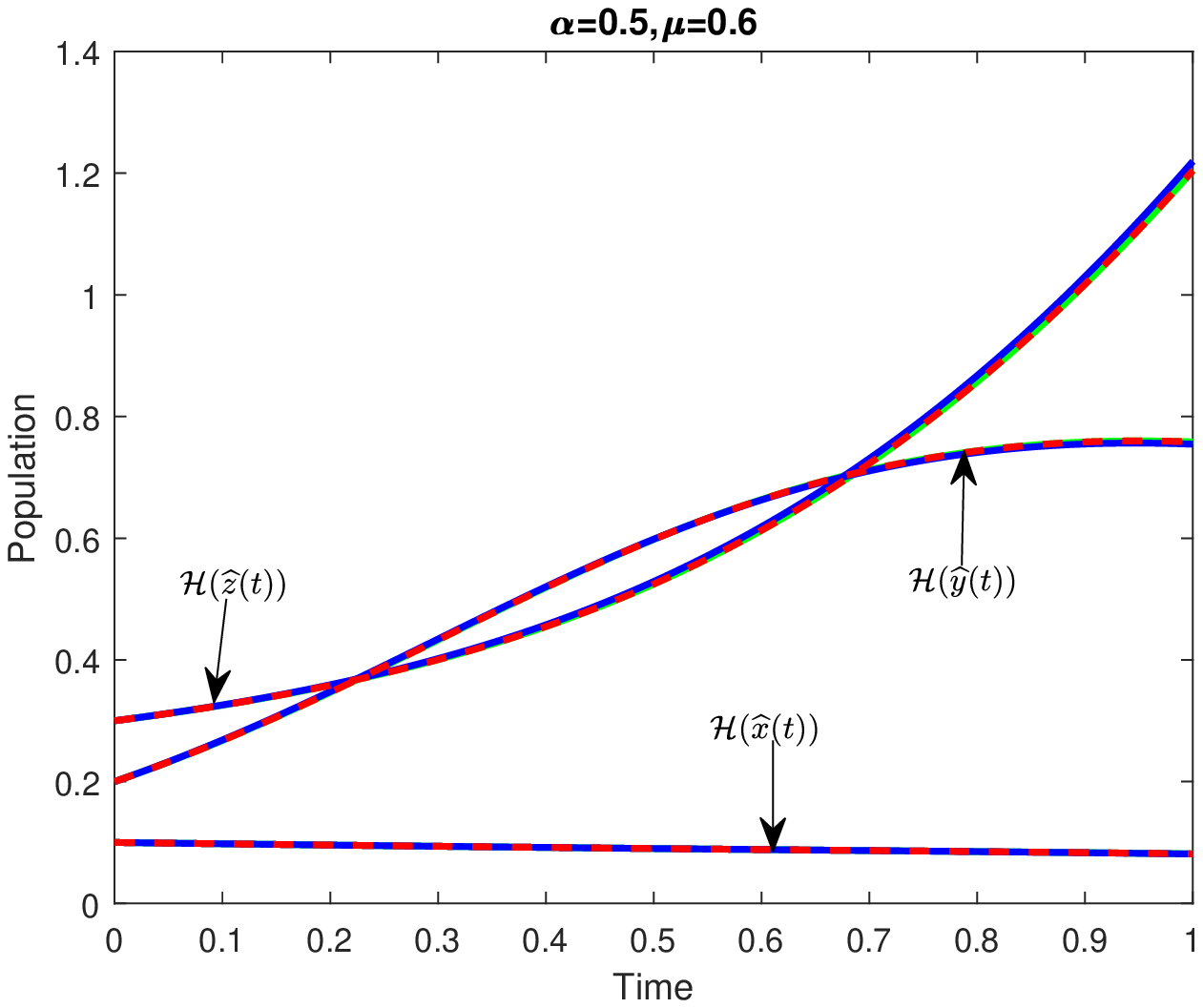}
\includegraphics[scale=0.5]{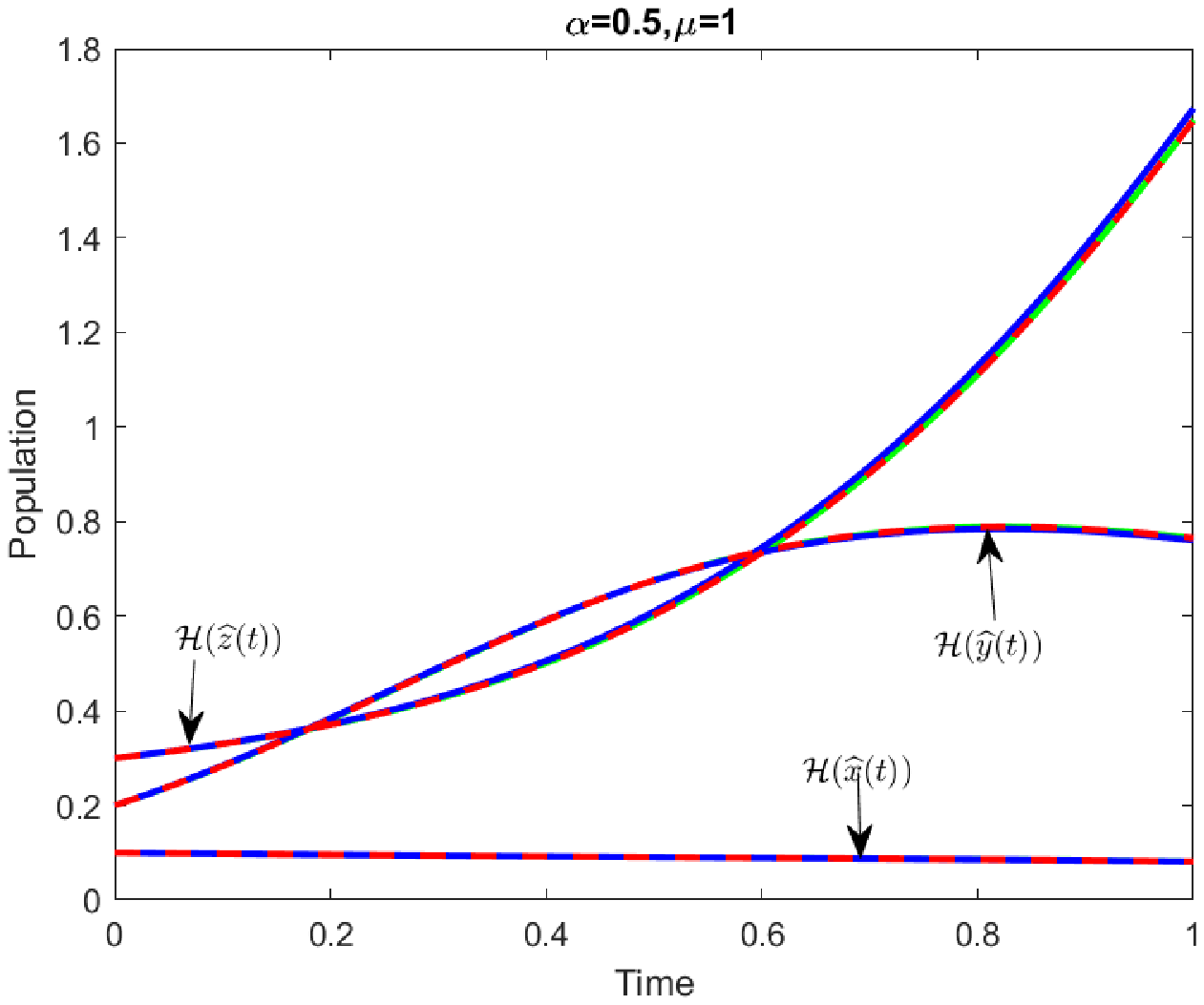}
\end{figure}
\begin{table}[ht]
\centering 
 \begin{adjustbox}{max width=0.90\linewidth}
\begin{tabular}{|c|c |ccc|ccc|ccc|} 
\hline
 $\mu$ &$u_i$ &\multicolumn{3}{|c|}{Exact} &\multicolumn{3}{|c|}{ Euler}
 &\multicolumn{3}{|c|}{FT-Euler mid-point}
\\ 
 & &$\mathbb{K}(\widehat{p}(u))$&$\mathbb{K}(\widehat{q}(u))$&$\mathbb{K}(\widehat{r}(u))$ &$\mathbb{K}(\widehat{p}(u))$& $\mathbb{K}(\widehat{q}(u))$ & $\mathbb{K}(\widehat{r}(u))$ & $\mathbb{K}(\widehat{p}(u))$ & $\mathbb{K}(\widehat{q}(u))$ & $\mathbb{K}(\widehat{r}(u))$ \\
\hline 
 & 0&0.1000&0.2000&0.3000&0.1000&0.2000&0.3000 & 0.1000&0.2000&0.3000 \\
 & 0.2&0.0954&	0.2573&	0.3309 &0.0954&	0.2571&	0.3307& 0.0954&	0.2573&	0.3309\\
{0} & 0.4& 0.0910&	0.3210&	0.3681&0.0910&	0.3207&	0.3677 & 0.0910&	0.3210&	0.3681 \\
 &0.6&0.0867&	0.3871&0.4135 &0.0867&0.3868&0.4128 &0.0867&	0.3872&0.4135 \\
 &0.8 &0.0825&	0.4505&	0.4689 &0.0825&0.4504&	0.4679 & 0.0825&	0.4506&	0.4690\\
  & 1&0.0785&0.5060&0.5362&0.0785&	0.5062&	0.5347&0.0784&0.5060&0.5363\\
  \hline
  & 0&0.1000&0.2000&0.3000&0.1000&	0.2000&	0.3000 & 0.1000&0.2000&0.3000 \\
 & 0.2&0.0956&0.3218&0.3505 &0.0956&0.3209&0.3499& 0.0956&0.3218&0.3505 \\
{0.4} & 0.4& 0.0916	&0.4646&0.4275&0.0916	&0.4635&0.4257 &0.0916	&0.4646&0.4275 \\
 & 0.6&0.0878&0.5958&0.5449&0.0878&0.5958&0.5412&0.0878&0.5958&0.5449 \\
 & 0.8 &0.0843&	0.6861&	0.7174& 0.0844&	0.6876&	0.7111&0.0843&	0.6861&	0.7174\\
  & 1&0.0806 &0.7258&0.9572 &0.0807	&0.7285&0.9478 & 0.0806 &0.7258&0.9572 \\
  \hline
  & 0&0.1000&	0.2000&	0.3000 &0.1000&	0.2000&	0.3000 & 0.1000&0.2000&0.3000 \\
 & 0.2&0.0957&	0.3569	&0.3636 & 0.0957&	0.3555	&0.3612& 0.0957&	0.3570	&0.3622\\
{0.6} & 0.4&0.0919&0.5383&	0.4723&	0.0919&0.5373&	0.4656&	0.0919&0.5385&	0.4687\\
 &0.6 &0.0884&	0.6823&	0.6541 &0.0884&	0.6823&	0.6404 & 0.0885&	0.6830&	0.6472\\
 &0. 8&0.0850&	0.7534&0.9379 &0.0852&	0.7577&0.9138 & 0.0851&	0.7548&0.9256\\
  & 1&0.0809&	0.7565&	1.3412 & 0.0813&	0.7630&	1.3042 & 0.0811&	0.7590&	1.3214\\
  \hline
  & 0&0.1000&	0.2000&	0.3000 &0.1000&	0.2000&	0.3000& 0.1000&0.2000&0.3000 \\
 & 0.2&0.0960&0.4313&0.3901 &0.0960&0.4289&0.3879 & 0.0960&0.4313&0.3901\\
{1} & 0.4& 0.0926&	0.6726&	0.5823 &0.0926&	0.6726&	0.5744& 0.0926&	0.6727&	0.5822 \\
 & 0.6&0.0897&	0.7961&	0.9541&0.0898&	0.7997&	0.9361 & 0.0897&	0.7963&	0.9538 \\
 & 0.8&0.0860&	0.8010&	1.5690 &0.0862&	0.8062&	1.5379&0.0860&	0.8012&	1.5685\\
  & 1&0.0789&	0.7335&	2.4055 &0.0796&	0.7395&	2.3667& 0.0789&	0.7342&	2.4046 \\
\hline 
\end{tabular}
\end{adjustbox}
\caption{Comparison of numerical results of Example \ref{grex5} for $\alpha=0$}
\label{grtab:5}
\end{table}
\begin{table}[ht]
\centering 
 \begin{adjustbox}{max width=0.90\linewidth}
\begin{tabular}{|c|c |ccc|ccc| ccc|} 
\hline
 $\mu$ &$u_i$ &\multicolumn{3}{|c|}{Exact} &\multicolumn{3}{|c|}{ Euler}
 &\multicolumn{3}{|c|}{FT-Euler mid-point}
\\ 
 & &$\mathbb{K}(\widehat{p}(u))$&$\mathbb{K}(\widehat{q}(u))$&$\mathbb{K}(\widehat{r}(u))$ &$\mathbb{K}(\widehat{p}(u))$& $\mathbb{K}(\widehat{q}(u))$ & $\mathbb{K}(\widehat{r}(u))$ & $\mathbb{K}(\widehat{p}(u))$ & $\mathbb{K}(\widehat{q}(u))$ & $\mathbb{K}(\widehat{r}(u))$ \\
\hline 
 & 0&0.1000&0.2000&0.3000&0.1000&	0.2000&	0.3000 & 0.1000&0.2000&0.3000 \\
 & 0.2&0.0956&	0.2967&	0.3426&0.0956&	0.2961&	0.3422 & 0.0956&	0.2967&	0.3426 \\
{0} & 0.4& 0.0913&	0.4092&	0.4020&0.0913&	0.4083&	0.4009&0.0913&	0.4092&	0.4020\\
 & 0.6&0.0874&	0.5202&	0.4857&0.0874&	0.5198&	0.4835 & 0.0874&	0.5203&	0.4857 \\
 & 0.8&0.0836&	0.6107&	0.6015&0.0836&	0.6113&	0.5979 & 0.0836&	0.6108&	0.6015 \\
  & 1&0.0799&	0.6685&	0.7569&	0.0800&	0.6702&	0.7515	&0.0799&	0.6685&	0.7569 \\
  \hline
  & 0&0.1000&0.2000&0.3000&0.1000&	0.2000&	0.3000&0.1000&0.2000&0.3000 \\
 & 0.2&0.0957&	0.3305&0.3533 &0.0957&0.3294&	0.3526 & 0.0957&	0.3304&0.3533\\
{0.4} & 0.4& 0.0916&	0.4832&	0.4370&0.0916&	0.4821&	0.4349& 0.0916&	0.4832&	0.4370 \\
 & 0.6&0.0880	&0.6193&	0.5678 &0.0880&	0.6194&	0.5634&0.0880	&0.6193&	0.5677\\
 & 0.8&0.0845&	0.7067&	0.7633&0.0846&	0.7086&	0.7559& 0.0845&	0.7068&	0.7632\\
  & 1&0.0808&	0.7382&	1.0374 &0.0809&	0.7414&	1.02627& 0.0808&	0.7384&	1.0372 \\
  \hline
  & 0&0.1000&	0.2000&	0.3000 &0.1000&	0.2000&	0.3000 & 0.1000&0.2000&0.3000 \\
 &0. 2&0.0957&	0.3480&	0.3592&0.0957&	0.3467&	0.3583& 0.0957&	0.3480&	0.3591\\
{0.6} &0.4&0.0918&	0.5202&	0.4576&	0.0918&	0.5190&	0.4549& 0.0918&	0.5202&	0.4576 \\
 & 0.6&0.0883&	0.6629&	0.6189&0.0883&	0.6636&	0.6130 & 0.0883&	0.6630&	0.6188\\
 & 0.8&0.0850&	0.7409&	0.8672&0.0850&	0.7435&	0.8570 & 0.0850&	0.7410&	0.8672\\
  & 1&0.0810&	0.7545&	1.2193&0.0812&	0.7583&	1.2042 & 0.0810&	0.7546&	1.2192\\
  \hline
 & 0&0.1000&	0.2000&	0.3000 &0.1000&	0.2000&	0.3000 & 0.1000&0.2000&0.3000 \\
 & 0.2&0.0958&	0.3843&	0.3719&0.0958&	0.3825&	0.3706& 0.0958&	0.3843&	0.3719 \\
{1} & 0.4&0.0922&	0.5918&	0.5061 &0.0921	&0.5908&	0.5015& 0.0922&	0.5918&	0.5060 \\
 & 0.6&0.0890&	0.7351&	0.7447&0.0890&	0.7371&	0.7346& 0.0890&	0.7352&	0.7446 \\
 & 0.8&0.0856&	0.7841&	1.1288 &0.0858&0.7880&	1.1111 & 0.0856&	0.7842&	1.1286 \\
  & 1&0.0808&	0.7598&	1.6740 &0.0812&	0.7647&	1.6489& 0.0808&	0.7600&	1.6738 \\
\hline 
\end{tabular}
\end{adjustbox}
\caption{Comparison of numerical results of Example \ref{grex5} for $\alpha=0.5$} 
\label{grtab:6}
\end{table}
\begin{landscape}

\begin{table}[ht]
\centering 
 \begin{adjustbox}{max width=0.65\linewidth}
 \begin{tabular}{|c |rrr|rrr| rrr|} 
\hline
 $u_i$ &\multicolumn{3}{|c|}{Exact} &\multicolumn{3}{|c|}{ Euler}
 &\multicolumn{3}{|c|}{FT-Euler mid-point}
\\ 
&$\mathbb{K}(\widehat{p}(u))$&$\mathbb{K}(\widehat{q}(u))$&$\mathbb{K}(\widehat{r}(u))$ &$\mathbb{K}(\widehat{p}(u))$& $\mathbb{K}(\widehat{q}(u))$ & $\mathbb{K}(\widehat{r}(u))$ & $\mathbb{K}(\widehat{p}(u))$ & $\mathbb{K}(\widehat{q}(u))$ & $\mathbb{K}(\widehat{r}(u))$ \\
\hline 
 0&0.1000&0.2000&0.3000&0.1000&0.2000&0.3000 & 0.1000&0.2000&0.3000 \\
  0.2&0.0957&0.3392&	0.3562 &0.0957&	0.3380&	0.3554& 0.0957&	0.3380&	0.3554 \\
 0.4& 0.0917&	0.5018&	0.4470&0.0917&	0.5006&	0.4446& 0.0917&	0.5018&	0.4470\\
 0.6&0.0882&	0.6417&	0.5924 &0.0882&	0.6421&	0.5873 & 0.0882&	0.6418&	0.5924\\
  0.8&0.0848&	0.7250&	0.8132&0.0848&0.7272&	0.8044 & 0.0848&	0.7250&	0.8132\\
  1&0.0809&	0.7477&	1.1247&0.0811&	0.7512&	1.1117 & 0.0809&	0.7478&	1.1246\\
\hline 
\end{tabular}
\end{adjustbox}
\caption{Comparison of numerical results of Example \ref{grex5} for $\alpha=1,\mu=0,0.4,0.6,1$}
\label{grtab:7}
\end{table}
\begin{table}[ht]
\centering 
 \begin{adjustbox}{max width=1.0\linewidth}
\begin{tabular}{|c|ccc|ccc|ccc|ccc|ccc|ccc|} 
\hline
 $\mu$&\multicolumn{9}{|c|}{ Euler}
 &\multicolumn{9}{|c|}{FT-Euler Mid-point}
\\
& & $\alpha=0$& && $\alpha=0.5$ & & & $\alpha=1$ & &&$\alpha=0$& && $\alpha=0.5$ & & & $\alpha=1$ & \\
\hline 
& &	& & & & & & &	& & & & & & & & & \\
{0} & 4.0000e-3&	2.1213e-4&	8.1035e-4&0.0000&5.7773e-5&5.7773e-5&4.0000e-3&8.7369e-4&	2.8381e-3& 0.0000& 5.7732e-5& 0.0000&8.1650e-5 &1.8317e-3 &6.7928e-3 & 0.0000&4.9497e-4 & 3.2914e-4\\
 & &	& & & & & & &	& & & & & & & & & \\
  \hline
  & &	& & & & & & &	& & & & & & & & & \\
 {0.4} & 5.7732e-5 & 1.3880e-3&4.7097e-3 & 0.0000& 0.0000&0.0000 &5.7732e-05&	1.6467e-3&5.8375e-3&0.0000&	1.0000e-4&	1.0000e-4& 8.1650e-5 &1.8317e-3 &6.7928e-3 & 0.0000&4.9497e-4 & 3.2914e-4 \\
 & &	& & & & & & &	& & & & & & & & & \\
 \hline
 & &	& & & & & & &	& & & & & & & & & \\
{0.6} &1.8257e-4&	3.2583e-3&1.6367e-2&1.0000e-4&	1.2076e-3&	3.1054e-3 & 8.1652e-5	&2.0339e-3 &7.9053e-3	& 0.0000&1.0000e-4	&1.2910e-4	 & 8.1650e-5 &1.8317e-3 &6.7928e-3 & 0.0000&4.9497e-4 & 3.2914e-4 \\
& &	& & & & & & &	& & & & & & & & & \\
  \hline
 & &	& & & & & & &	& & & & & & & & & \\
{1} & 3.0000e-4&3.6914e-03&	2.1848e-2&	0.0000 &3.2145e-4&	4.3970e-4&1.8708e-4	 & 2.7404e-3&1.3332e-2	& 0.0000&1.2247e-4 &1.2247e-4	&	8.1650e-5 &1.8317e-3 &6.7928e-3 & 0.0000&4.9497e-4 & 3.2914e-4 \\
& &	& & & & & & &	& & & & & & & & & \\
\hline 
\end{tabular}
\end{adjustbox}
\caption{RMS error for Example \ref{grex5} for different values of $\alpha,\mu$} 
\label{grtab:8}
\end{table}
\end{landscape}
\section{Conclusion}
In this contribution, we have introduced and studied the concepts of granular $F$-transforms to enrich the theory of $F$-transforms and explore new applications. Accordingly, we have initiated the said theory, formulated a fuzzy prey-predator model consisting of two prey and one predator team, and discussed the equilibrium points and their stability for this model. Moreover, we have established a numerical method based on the granular $F$-transform to find the numerical solution to the proposed model. Finally, a comparison between two numerical solutions with the exact solution is discussed. In the future, it will be interesting to use the proposed numerical method based on granular $F$-transforms to solve fuzzy fractional differential equations and analyze error estimation. 


\begin{thebibliography}{00}
\bibitem{aziz} M.A. Aziz-Alaoui, M.D. Okiye, Boundedness and global stability for a predator-prey model with modified Leslie-Guwer and Holling-type II schemes, {\em Applied Mathematics Letters}, {\bf 16} (2003) 1069-1075.
\bibitem{bede} B. Bede, S.G. Gal, Generalization of the differentiability of fuzzy-number-valued functions with applications to fuzzy differential equation, {\em Fuzzy Sets and Systems}, {\bf151} (2005) 581-599.
\bibitem{bede1} B. Bede, L. Stefanini, Generalized differentiability of fuzzy-valued functions, {\em Fuzzy Sets and Systems}, {\bf230} (2013) 119-141.
\bibitem{bell} R. Bellman, Stability theory of differential equations, New York: MacGraw-Hill, (1953).
\bibitem{brau} F. Brauer, J.A. Nohel, The qualitative theory of ordinary differential equations: A introduction, New York: Dover Publications Inc, (1969). 
\bibitem{chen} W. Chen, Y. Schen, Approximate solution for a class of second-order ordinary differential equations by the
fuzzy transform, {\em Journal of Intelligent and Fuzzy Systems}, {\bf 27} (2014) 73-82.
\bibitem{cush} J.M. Cushing, Periodic Lotka-Volterra competition equation, {\em Journal of Mathematical Biology}, {\bf 24} (1986) 381-403.
\bibitem{mar} F. Di Martino, V. Loia, I. Perfilieva, S. Sessa, An image coding/decoding method based on direct and inverse fuzzy transforms, {\em International Journal of Approximate Reasoning}, {\bf 48} (2008) 110-131.
\bibitem{mar1} F. Di Martino, V. Loia, S. Sessa, A segmentation method for images compressed by fuzzy transforms, {\em Fuzzy Sets and Systems}, {\bf 161} (2010) 56-74.
\bibitem{mar2} F. Di Martino, V. Loia, S. Sessa, Fuzzy transforms method in prediction data analysis, {\em Fuzzy Sets and Systems}, {\bf 180} (2011) 146-163.
\bibitem{elet} M.F. Elettreby, Two-prey one-preydator model, {\em Chaos Solitons and Fractals}, {\bf39} (2009) 2018-2027.
\bibitem{gak} S. Gakkhar, B. Singh, R.K. Naji, Dynamical behaviour of two predators competing over a single prey, {\em Biosystems}, {\bf 90} (2007) 808-817.
\bibitem{gakk} S. Gakkhar, B. Singh, The dynamics of a food web consisting of two preys and a harvesting predator, {\em Chaos Solitons and Fractals}, {bf 34} (2007) 1345-1356.
\bibitem{hoa} N.V. Hoa, Fuzzy fractional functional differential equations under Caputo gH-differentiability, {\em Communications in Nonlinear Science and Numerical Simulation}, {\bf 22} (2015) 1134-1157.
\bibitem{hoa1} N.V. Hoa, V. Lupulescuc D. O'Regand, Solving interval-valued fractional initial value problems under Caputo gH-fractional differentiability, {\em Fuzzy Sets and systems}, {\bf 309} (2017) 1-34.
\bibitem{holc} M. Hol\v{c}apek, L. Nguyen, Trend-cycle estimation using fuzzy transform of higher degree, {\em Iranian Journal of Fuzzy Systems}, {\bf15} (2018) 23-54.
\bibitem{hut} P. Hurt\'{\i}k, S. Tomasiello, A review on the application of fuzzy transform in data and image compression, {\em Soft Computing}, {\bf 23} (2019) 12641-12653.
\bibitem{hsu} S.B. Hsu, T.W. Hwang, Y. Kuang, Global analysis of the Michaelis-Menten-type ratio-dependent predator-prey system, {\em Journal of Mathematical Biology}, {\bf42} (2001) 489-506.
\bibitem{ali} A. Khastan, Z. Alijani, I. Perfilieva, Fuzzy transform to approximate solution of two-point boundary value problems, {\em Mathematical Methods in the Applied Sciences}, {\bf40} (2017) 6147-6154.
{\bibitem{kh1} A. Khastan, A new representation for inverse fuzzy transform and its application, {\em Soft Computing}, {\bf 21} (2017) 3503-3512.}
\bibitem{kh} A. Khastan, I. Perfilieva, Z. Alijani, A new fuzzy approximation method to Cauchy problems by fuzzy transform, {\em Fuzzy Sets and Systems}, {\bf 288} (2016), 75-95.
\bibitem{land} M. Landowski, Usage of RDM interval arithmetic for solving cubic interval equation, in: {\em  Advances in Fuzzy Logic and Technology}, {\bf 3} (2017) 382–391.
\bibitem{li} M. Liu, D. Chen, C. Wu, H. Li, Approximation theorem of the fuzzy transform in fuzzy reasoning and its application to the scheduling problem, {\em Computers and Mathematics with Applications}, {\bf 51} (2006) 515-526.
\bibitem{liu} B. Liu, Z. Teng, L. Chen, Analysis of predator-prey model with Holling II functional response concerning impulsive control strategy, {\em Journal of Computational and Applied Mathematics}, {\bf 193} (2006) 1147-1162.
\bibitem{liou} L.P. Liou, K.S. Cheng, Global stability of a predator-prey system, {\em Journal of Mathematical Biology}, {\bf26} (1988) 65-71.
\bibitem{long} H.V. Long, N.T.K. Son, H.T.T. Tam, J.C. Yao, Ulam stability for fractional partial integro-differential equation with uncertainty, {\em Acta Mathematica Vietnamica}, {\bf 42} (2017) 675-700.
\bibitem{long1} H.V. Long, N.T.K. Son, H.T.T. Tam, The solvability of fuzzy fractional partial differential equations under Caputo gH-differentiability, {\em Fuzzy Sets and Systems}, {\bf 309} (2017) 35-63.
\bibitem{long2} H.V. Long, N.T.K. Son, N.V. Hoa, Fuzzy fractional partial differential equations in partialy ordered metic spaces, {Iranian Journal of Fuzzy Systems}, {\bf 14} (2017) 107-126.

\bibitem{maz} M. Mazandaeani, M. Najariyan, Differentiability of type-2 fuzzy number-valued functions, {\em Communications in Nonlinear Science and Numerical Simulation}, {\bf19} (2014) 710-725.
\bibitem{mazan} M. Mazandaeani, N. Parizz, A.V. Kamyad, Granular differentiability of fuzzy-number-valued functions, {\em IEEE Transitions on Fuzzy Systems}, {\bf26} (2018) 310-323.
\bibitem{jir} J. Mo\v{c}ko\v{r}, $F$-transforms and semimodule homomorphisms, {\em Soft Computing}, {\bf 23} (2019) 7603-7619.
\bibitem{mo} J. Mo\v{c}ko\v{r}, M. Hol\v{c}apek, Fuzzy objects in spaces with fuzzy partitions, {\em Soft Computing}, {\bf 21} (2016) 7268-7284.
\bibitem{mocko} J. Mo\v{c}ko\v{r}, P. Hurt\'{\i}k, Lattice-valued $F$-transforms and similarity relations, {\em Fuzzy Sets and Systems}, {\bf 342} (2018) 67-89.
\bibitem{murray} J.D. Murray, Mathematical Biology: An introduction, {\em Springer}, New Delhi, (2002).
\bibitem{naj}  M. Najariyan, Fuzzy fractional quadratic regulator problem under graular fuzzy fractional derivative, {\em IEEE Transitions on Fuzzy Systems}, {\bf26} (2018) 2273-2288.
\bibitem{naja}  M. Najariyan,Y. Zhao, On the stability of fuzzy linear dynamical systems, {\em Journal of the Franklin Institute}, {\bf 357} (2020) 5502-5522.
\bibitem{vil} V. Nov\'{a}k, I. Perfilieva, M. Hol\v{c}apek, V. Kreinovich, {Filtering out high frequencies in time series using F-transform}, {\em Information Sciences}, {\bf 274} (2014) 192-209.
\bibitem{oat} A. Oaten, W.W. Murdoch, Functional response and stability in predator-prey systems, {\em The American Naturalist}, {\bf 109} (1975) 289-298.
\bibitem{piegat2} A. Piegat, M. Landowski, Fuzzy arithmetic type 1 with horizontal membership functions, in: {\em Uncertainty Modeling}, {\bf 683} (2017) 233–250.
\bibitem{past} J. Paster, Mathematical ecology of poplations and ecosystems, {\em West Sussex: A John Wiley and Sons Ltd Publication}, (2008).
\bibitem{piegat} A. Piegat, M. Pluci\'{n}ski, Some advantages of the RDM-arithmetic of intervally-precisiated Values, {\em International Journal of Computational Intelligence Systems}, {\bf 8} (2015) 1192-1209.
\bibitem{piegat1} A. Piegat, M. Pluci\'{n}ski, Fuzzy number division and the multi-granularity phenomenon, {\em Bulletin of the Polish Academy of Sciences Technical Sciences}, {\bf 65} (2017) 497-511.
\bibitem{per} I. Perfilieva, $F$-transforms: theory and its applications, \emph{Fuzzy Sets and Systems}, {\bf157} (2006) 993-1023.
\bibitem{p} I. Perfilieva, Fuzzy transform: application to the reef growth problem. In: {\em Fuzzy Logic in Geology}; Demicco, R.V., Klir, G.J., Eds., Academic Press: Amsterdam, The Netherlands, (2003) 275-300.

\bibitem{irin1} I. Perfilieva, Fuzzy transforms: a challenge to conventional  transforms, {\em Advances in Image and Electron Physics}, {\bf 147} (2007) 137-196.
\bibitem{no} I. Perfilieva, V. Nov\'{a}k, A. {Dvo\v{r}\'{a}k}, {Fuzzy transforms in the analysis of data}, {\em International Journal of Approximate Reasoning}, {\bf 48} (2008) 36-46.
\bibitem{anan} I. Perfilieva, A.P. Singh, S.P. Tiwari, On the relationship among $F$-transform, fuzzy rough sets and fuzzy topology, \emph{Soft Computing}, {\bf21} (2017) 3513-3523.
\bibitem{stev} I. Perfilieva, P. \v{S}tevuli\'{a}kova R. Val\'{a}\v{s}ek, $F$-Transform for numerical solution of two-point boundary value problem, {\em Iranian Journal of Fuzzy Systems}, {\bf 14} (2017) 1-13.
\bibitem{spt1} I. Perfilieva, S.P. Tiwari, A.P. Singh, Lattice-valued $F$-transforms as interior operators of $L$-fuzzy pretopological spaces, {\em Communications in Computer and Information Science}, {\bf 854} (2018) 163-174.
\bibitem{Ir} I. Perfilieva, R. Valasek, {Fuzzy transforms in removing noise}, {\em Advances in Soft Computing}, {\bf 2} (2005) 221-230.
\bibitem{rus} C. Russo, Quantale modules and their operators, with application, {\em Journal of Logic and Computation}, {\bf20} (2010) 917-946.
\bibitem{roh} S. B. Roh, S. K. Oh, J. H. Yoon, K. Seo, {Design of face recognition system based on fuzzy transform and radial basis fnction neural networks}, {Soft Computing}, {\bf 23} (2019) 4969-4985.
\bibitem{son} N.T.K Son, N.P. Dong, L.H. Son, H.V. Long, Towards granular calculus of single-valued neutrosophic functions under granular computing, {\em Multimedia Tools and Applications}, {\bf 79} (2020) 16845-16881. 
\bibitem{st} L. Stefanini, $F$-transform with parametric generalized fuzzy partitions, {\em Fuzzy Sets and Systems}, {\bf 180} (2011) 98-120.
\bibitem {ste} M. {\v{S}t\v{e}pni\v{c}ka}, O. {Polakovi\v{c}}, {A neural network approach to the fuzzy transform}, {\em Fuzzy Sets and Systems}, {\bf 160} (2009) 1037-1047.
\bibitem{step} M. {\v{S}t\v{e}pni\v{c}ka}, R. {Val\'{a}\v{s}ek}, Fuzzy transforms and their application to wave equation, {\em Journal of Electrical Engineering}, {\bf 55} (2004) 7-10.
\bibitem{spt} S.P. Tiwari, I. Perfilieva, A.P. Singh, Generalized residuated lattices based $F$-transform, {\em Iranian Journal of Fuzzy Systems}, {\bf15} (2018) 63-182.
\bibitem{tom} S. Tomasiello, An alternative use of fuzzy transform with application to a class of delay differential equations, {\em International Journal of Computational Mathematics}, {\bf 94} (2017) 1719-1726.
\bibitem{jp} J.P. Tripathi, S, Abbas, M. Thakur, Local and global stability analysis of a two prey one predator model with help, {\em Communications in Nonlinear Science and Numerical Simulation}, {\bf 19} (2014) 3284-3297.
\bibitem{tri} A. Tripathi, S.P. Tiwari, I. Perfilieva, $F$-transforms determined by implicators, {\em Iranian Journal of Fuzzy Systems}, {\bf 18} (2021) 19-36.
\bibitem{to} L. Troiano, P. Kriplani, {Supporting trading strategies by inverse fuzzy transform}, {\em Fuzzy Sets and Systems}, {\bf 180} (2011) 121-145.
\bibitem{sun} S. Zhang and J. Sun, Stability of fuzzy differential equations with the second type
of Hukuhara derivative, {\em IEEE Transactions on Fuzzy Systems}, {\bf 23} (2015) 1323-1328.
\bibitem{sun1} S. Zhang and J. Sun, Practical stability of fuzzy differential equations with the second type of Hukuhara derivative, {\em Journal of Intelligent and Fuzzy Systems}, {\bf29} (2015) 307-313.

\end{thebibliography}
\end{document}